%% file: paper.tex
\documentclass{article}
\input{command}
\newcommand*\xbar[1]{%
  \hbox{%
    \vbox{%
      \hrule height 0.5pt 
      \kern0.4ex
      \hbox{%
        \kern-0.05em
        \ensuremath{#1}%
        \kern-0.00em
      }%
    }%
  }%
}

\usepackage[colorinlistoftodos,prependcaption,textsize=tiny]{todonotes}

\newcommand{\red}[1]{{{#1}}}

\begin{document}


\title{Virtual finite element and hyperbolic problems: the PAMPA algorithm.}
\author{R\'{e}mi Abgrall$^\star$, Walter Boscheri$^\dagger$, and Yongle Liu$^\star$\\
$(\star)$ Institute of Mathematics, University of Z\"urich, Switzerland\\
$(\dagger$) Laboratoire de Math\'{e}matiques UMR 5127 CNRS, Universit\'{e} Savoie Mont Blanc, \\France}
\date{}
\maketitle
\begin{abstract}
In this paper, we explore the use of the Virtual Element Method (VEM) concepts to solve scalar and system hyperbolic problems on general polygonal grids. The new schemes stem from the Active Flux approach \cite{AF1}, which combines the usage of point values at the element boundaries with an additional degree of freedom representing the average of the solution within each control volume. Along the lines of the family of residual distribution schemes introduced in \cite{Abgrall_AF,abgrall2023activefluxtriangularmeshes} that integrate the Active Flux technique, we devise novel third order accurate methods that rely on the VEM technology to discretize gradients of the numerical solution by means of a polynomial-free approximation, by adopting a virtual basis that is locally defined for each element. The obtained discretization is globally continuous, and for nonlinear problems it needs a stabilization which is provided by a monolithic convex limiting strategy extended from \cite{Abgrall_BP_PAMPA}. This is applied to both point and average values of the discrete solution. We show applications to scalar problems, as well as to the acoustics and Euler equations in two dimension. The accuracy and the robustness of the proposed schemes are assessed against a suite of benchmarks involving smooth solutions, shock waves and other discontinuities.
\end{abstract}

\input{intro}

\section{Meshes and Approximation space}
\subsection{Meshes}
Given $\Omega\subset \R^2$ that is assumed to be polygonal, we first start by constructing a triangular mesh using GMSH \cite{gmsh}. GMSH can also consider quadrangular meshes. This is also doable in 3D. If we consider only triangular or quadrilateral meshes, we fit the formalism outlined in \cite{Abgrall_AF}. Otherwise, arbitrary polygonal meshes must be faced, and we consider the following options for generating the computational mesh.
\begin{enumerate}
\item The centroids, i.e. the barycenters, of the GMSH elements is connected with the mid points of the edges, hence obtaining a dual mesh with respect to the original one \cite{maire2011,HTCLag_boscheri2024}.
\item A genuinely Voronoi mesh can be constructed from the vertices of the GMSH mesh by connecting the circumcenters of the GMSH elements which share a common vertex.
\item A more regular polygonal grid is built starting from the vertices of the GMSH mesh by connecting the barycenters of the GMSH elements which share a common vertex. This is no longer a Voronoi mesh, but typically it yields more regular hexagonal polygons \cite{CWENO_BGK_boscheri2020}.
\end{enumerate}
In all cases, the physical boundary of the domain is preserved, thus the polygonal boundary elements are modified accordingly. Option 3 is mostly adopted in this paper, because it leads to a higher quality of the element shape, and our numerical method does not need any orthogonality property which would otherwise imply the usage of a Voronoi grid \textit{sensu stricto} \cite{ader_fse_boscheri2020}.

\subsection{Approximation space}
The approximation space is the same of the one adopted by the VEM \cite{hitch}.  We first introduce some notations, following closely \cite{hitch}. The computational domain $\Omega$ is covered by a set of non-empty and non-overlapping polygons that are denoted by $P$. The notation $\bbx_{P}$ represents the centroid of $P$. The element $P$ can be of very general shape (convex or nonconvex polygons), but they are assumed to be star-shaped, with respect to a point $\bbx^\star$ (that may be different of $\bbx_{P}$), for the sake of designing a low-order scheme in section \ref{sec32}. For $P\subset \Omega$, the $L^2$ inner product between two functions in $L^2(P)$ is $\langle u,v\rangle_{P}$. When there is no ambiguity on $P$, we omit the subscript $P$. For $\bsa=(\alpha_1, \alpha_2)$, the scaled monomial of degree $\vert \bsa\vert=\alpha_1+\alpha_2$ is defined by
\begin{equation}\label{scaled_monomial}
 m_{\bsa}=\big ( \dfrac{\bbx-\bbx_{{P}}}{h_{{P}}}\big )^{\bsa}=\big ( \dfrac{x-x_{P}}{h_{P}}\big )^{\alpha_1} \big ( \dfrac{y-y_{P}}{h_{P}}\big )^{\alpha_2}, \qquad \bbx=(x,y),
\end{equation}
where $h_P$ denotes the diameter of $P$. Other choices for the polynomial basis on polygons are possible. For instance, one may use a set of orthogonal polynomials\footnote{This set should contain the constant function $1$. This kind of choice allows for a better conditioning of the linear systems.}, as discussed in \cite{orthVEM2017}. The set of scaled monomials of degree $\vert\bsa\vert\leq k$ is a basis for $\P^k(P)$, which denotes the vector space of polynomials of degree less than or equal to $k$, defined on $P$. Similar definitions and comments can be done in three space dimensions. The scaled monomials are invariant by homothety: if $\bbx=\lambda \widehat{\bbx}$, then 
\begin{equation*}
 \widehat{m}_\bsa(\widehat{\bbx})=\big ( \dfrac{\widehat{\bbx}-\bbx_{\widehat{P}} }{h_{\widehat{P}}}\big )^{\bsa}=m_{\bsa}(\bbx)=\big ( \dfrac{\bbx-\bbx_{P}}{h_{P}}\big )^{\bsa},
\end{equation*}
because $\widehat{P}$ is the image of $P$ by this mapping.
 
Now we introduce the DoFs of the approximation space and define the local virtual space $V_k(P)$ for each $P$. For $k\geq 1$, a function $v_h\in V_k(P)$ is uniquely defined by the following setup:
\begin{enumerate}
 \item $v_h$ is a polynomial of degree $k$ on each edge $e$ of ${P}$, that is $(v_h)_{|e}\in \P^k(e)$;
 \item $\Delta v_h\in \P^{k-2}(P)$ (with the convention $\P^{-1}(P)=\{0\}$ when needed).
\end{enumerate}

The second condition has nothing to do with any PDE problem we could have in mind, we will comment later on its usefulness. Notice that a polynomial of degree $k$ satisfies the aforementioned conditions so that $\P^k(P)\subset V_k(P)$, \red{and a function\footnote{Roughly speaking, $V_k(P)$ contains all polynomials of degree $k$ (which is essential for convergence) plus other functions whose restriction on an edge is still a polynomial of degree $k$.} }of $V_k({P})$ is uniquely defined by the DoFs given by:
\begin{enumerate}
 \item The value of $v_h$ at the vertices of $P$,
 \item On each edge of $P$, the value of $v_h$ at the $k-1$ internal points of the $k+1$ Gauss--Lobatto points on this edge,
 \item The moments up to order $k-2$ of $v_h$ in ${P}$,
   \begin{equation}\label{moments}
      m_{\bsa}(v_h):=\dfrac{1}{\vert P\vert }\int_{P} v_h m_{\bsa}(\bbx)\; {\rm d}\bbx, \qquad |\bsa|\leq k-2.
    \end{equation}
\end{enumerate}
The dimension of $V_k(P)$ is 
\begin{equation*}
  \dim V_k(P)={N_V\cdot } k+\frac{k(k-1)}{2},
\end{equation*}
where $N_V$ represents the number of vertices of $P$. The total number of DoFs in $P$ {is} then referred to as $N_{\text{DoFs}}:=\dim V_k(P)$. Let $\{\varphi_i\}_{i=1}^{N_{\text{DoFs}}}$ be the canonical basis for $V_k(P)$ and $\text{DoF}_i$ be the operator from $V_k(P)$ to $\R$ as
\begin{equation*}
  \text{DoF}_i(v_h)=i-{\rm th} \text{ degree of freedom of }v_h,\quad i=1,\ldots,N_{\text{DoFs}}.
\end{equation*}
We can then represent each $v_h \in V_k(P)$ in terms of its DoFs by means of:
\begin{equation*}
  	v_h = \sum_{i=1}^{N_{\text{DoFs}}} \text{DoF}_i(v_h)\varphi_i.
\end{equation*}
For this basis, the usual interpolation property holds true:
\begin{equation*}
  	\text{DoF}_i(\varphi_j) = \delta_{ij}, \qquad i,j=1, \ldots, N_{\text{DoFs}}.
\end{equation*}

Following the VEM literature,  we recall that the explicit computation of the basis functions $\varphi_i$ is actually not needed. The first step is to construct a projector $\pi^\nabla$ from $V_k(P)$ onto $\P^k(P)$. It is defined by two sets of properties. First, for any $v_h\in V_k({P})$, the orthogonality condition
\begin{equation}\label{projo:1}
 \big \langle \nabla p_k, \nabla \big ( \pi^\nabla v_h-v_h\big )\big \rangle =0, \quad \forall p_k\in {\P^k}(P),
\end{equation}
has to hold true, which is defined up to the projection onto constants $P_0: V_k(P)\rightarrow \P^0(P)$, that can be fixed as follows:
\begin{itemize}
 \item if $k=1$,
 \begin{equation*}
 {\mathcal{P}}_0(v_h)=\frac{1}{N_V}\sum_{i=1}^{N_V} v_h(\bbx_i),
 \end{equation*}
 \item if $k\geq 2$,
  \begin{equation*}
  {\mathcal{P}}_0(v_h)=\frac{1}{\vert P\vert }\int_{P} v_h \; {\rm d}\bbx =m_{(00)}(v_h).
  \end{equation*}
\end{itemize}
Then, we ask that 
\begin{equation}\label{projo:2}
 {\mathcal{P}}_0( \pi^\nabla v_h-v_h)=0.
\end{equation}
We can \textit{explicitly} compute the projector by using only the DoFs previously introduced: using the third condition of the definition of {$V_k(P)$}, it can be seen that ($\bbn$ is the outward unit normal to $\partial P$)
\begin{equation*}
\int_P\nabla p_k\cdot \nabla v_h\; {\rm d}\bbx=-\int_P \Delta p_k \cdot v_h\; {\rm d}\bbx+\int_{\partial P} \nabla p_k\cdot \bbn v_h\; {\rm d}\gamma, \quad \forall p_k\in\P^k(P),
\end{equation*}
so that $\int_P\nabla p_k\cdot\nabla v_h\; {\rm d}\bbx$ is computable from the DoFs only. Since $\pi^\nabla v_h\in {\P^k}(P)$, we can find $\dim {\P^k}(P)=\frac{{(k+2)}(k+1)}{2}$ real numbers $s_{\bsa}$ such that 
 \begin{equation*}
 \pi^\nabla v_h=\sum_{\vert\bsa\vert=1}^{\dim \P^k({P})} s_{\bsa} m_{\bsa},
\end{equation*} 
and then
\begin{equation}\label{proj:1}
 \mathcal{P}_0(\pi^\nabla v_h)=\sum_{\bsa, |\bsa|\leq k} s_{\bsa}\cdot {\mathcal{P}}_0(m_{\bsa})={\mathcal{P}}_0(v_h).
\end{equation}
Moreover, in \eqref{projo:1} we let $p_k$ vary only for the scaled monomial basis defined in \eqref{scaled_monomial} and we denote it by $m_{\bsb}$. Integration by parts yields
\begin{equation*}
 \sum_{\bsa, |\bsa|\leq k} s_{\bsa} \langle \nabla m_{\bsa}, \nabla m_{\bsb}\rangle =\int_{P}\nabla m_{\bsb} \nabla v_h\; {\rm d}\bbx=-\int_{P} \Delta m_{\bsb} v_h\; {\rm d}\bbx+\int_{\partial P} v_h\; \nabla m_{\bsb} \cdot \bbn  \; {\rm d}\gamma,
\end{equation*}
 and since $\Delta m_{\bsb}\in \P^{k-2}(P)$, one can find a family of real coefficients $d_{\bm\delta}(m_{\bsb})$ such that 
 $$\Delta m_{\bsb}=\sum_{\bm\delta, |\bm\delta|\leq k-2} d_{{\bm\delta}}\big (m_{\bsb}\big ) m_{\bm\delta},$$
 the volume integral defined above can be computed by
\begin{equation}\label{proj:2}
 \int_{P} v_h\Delta m_{\bsb} \; {\rm d}\bbx=\sum_{\bm\delta, |\bm\delta|\leq k-2} d_{\bm\delta}\big (m_{\bsb}\big )\;\int_{P} m_{\bm\delta} v_h\; {\rm d}\bbx\stackrel{\eqref{moments}}{=}
 \vert {P} \vert \sum_{\bm\delta, |\bm\delta|\leq k-2} d_{\bm\delta}\big (m_{\bsb}\big ) m_{\bm\delta}(v_h).
\end{equation}
Similarly, $\int_{\partial P} v_h\; \nabla m_{\bsb} \cdot \bbn  \; {\rm d}\gamma$ is computable because $v_h$ is a polynomial of degree $k$ on $\partial {P}$. Note that, since the point values at the Gauss--Lobatto points are known, no additional computation is needed.
 
By gathering the information contained in \eqref{proj:1} and \eqref{proj:2}, we obtain a linear system 
 $$G_P \mathbf{s}=\mathbf{c}$$
 with $\mathbf{s}$ the vector of components $s_{\bsa}$ and $\mathbf{c}$ the vector of components
 $$c_{ \bsb}=-\int_{P} {\Delta m_{\bsb} \; v_h}\; {\rm d}\bbx+\int_{\partial P} \nabla m_{\bsb} \cdot \bbn \, v_h \; {\rm d}\gamma.$$
The matrix $G$ is 
$$G_P=\begin{pmatrix}
A_P\\B_P
\end{pmatrix},$$
where
$A_P$ is the vector containing the coefficients ${\mathcal P}_0(m_{\bsa})$ for $\vert\bsa\vert \leq k$ and $B_P$ is the {$(\dim \P^k(P)-1) \times \dim \P^k(P)$} ``mass matrix'':\begin{equation*}
  B_P=\begin{pmatrix}\langle \nabla m_{\bsa},\nabla m_{\bsb}\rangle\end{pmatrix}_{0{\leq}\vert \bsa\vert,\vert\bsb\vert\leq k}.
\end{equation*}
{It is invertible because from \eqref{projo:1} its kernel contains constant polynomials only, and from \eqref{projo:2} this can be only $0$. {The matrix $G_{{P}}$ can be conveniently computed, inverted and could be  stored once and for all in the pre-processing step, thus improving the efficiency of the overall algorithm. What we have chosen to do instead is to store the coefficients needed to evaluate the gradients at the Lagrange points on the boundary of the polygons, and what is needed to evaluate, at these points, the operator $\mathcal{D}_\sigma$ of \eqref{scheme:stab}. All this is described in section \ref{sec31}.}}

 \bigskip
{ Now we describe the approximation setting on $\Omega$. \red{On} any polygonal domain of $\R^2$ that is covered by non overlapping polygons $P_i$,
$$\Omega=\cup_{i=1}^{n_P} P_i,$$ one can construct a globally continuous approximation of $u$ on $\Omega$ by setting, for all $P_i$
$u_{\mid P_i}\in $ defined by the following degrees of freedom
\begin{enumerate}
    \item on vertex of the polygon, the value of $u$ at \red{this} vertex
    \item on each edge, the values of $u$ on the $k-1$ internal points of the $(k+1)$ Gauss-Lobatto quadrature rule on this edge,
    \item for each polygon, the moments up to order $k-2$ of $u$ in $P$:
    $$\frac{1}{\vert P_i\vert}\int_{P_i} u(\bbx) m_{\bm \alpha}(\bbx)\; d\bbx, \quad \vert \bm \alpha\vert \leq k-2.$$
    \end{enumerate}
    This provides a globally continuous approximation of $u$ on $\Omega$. See \cite{hitch} for the functional analysis details.
}
 


\section{Numerical schemes}\label{sec3}
The mathematical model is given by 
\begin{equation}
\label{scheme:fv}\dpar{\bbu}{t}+{\rm div~}\bbf(\bbu)=0,
\end{equation} 
where $\bbu{(t,\bbx)}\in \mathcal{D}\subset \R^{{m}}$ is the vector of conserved variables, $\mathcal{D}$ is the {convex invariant} domain where $\bbu$ and the flux tensor $\bbf=(f_1, \ldots , f_d)$ are defined. The functions  $f_j$ for every $j=1,\ldots,d$ are assumed to be defined and $C^1$ on $\mathcal{D}$. {The convex invariant domain will be specified according to the problem studied.} This system, at least for smooth solutions, can be rewritten in a non-conservative form as
\begin{equation*}
\dpar{\bbu}{t}+\bbA\cdot\nabla \bbu=0,
\end{equation*}
with the notation
$$\bbA\cdot\nabla \bbu=\sum_{j=1}^d A_j \dpar{\bbu}{x_j} \text{ with }A_j=\dpar{f_j}{\bbu}.$$
The governing equations are assumed to be hyperbolic, i.e. for any $\bbn=(n_1, \ldots , n_d)\in \R^d$, the system matrix
$$\bbA\cdot \bbn=\sum_{j} A_j n_j$$
is diagonalisable in $\R$.

The canonical example of such a system is that of the Euler equations where, if $\rho$ is the density, $\bbv$ the velocity and $E$ the total energy, we have $\bbu=(\rho, \rho \bbv, E)^\top$. The total energy is the sum of the internal energy $\epsilon$ and of the kinetic energy $\tfrac{1}{2}\rho \bbv^2$. The invariant domain $\mathcal{D}$ is
\begin{equation}\label{Euler_invariant}
  \mathcal{D}:=\{ \bbu{=(\rho,\rho\bbv,E)^\top\in\R^4}, \text{ with } \rho>0, \epsilon{=E-\frac{1}{2}\rho\Vert\bbv\Vert^2}> 0\}.
\end{equation}
{\red{The invariant domain  }can be rewritten as 
\begin{equation}\label{Euler_GQL}
   \mathcal{D}=\{\bbu=(\rho, \rho\bbv, E)^\top\in\R^4 \text{ such that for all } \bbn_*\in \mathcal{N}, \bbu^\top\bbn_*>0\}, 
\end{equation}
where, through the GQL (Geometric Quasilinearization) approach \cite{wu2023geometric}, the \red{set} $\mathcal{N}$ is 
\begin{equation*}
     \mathcal{N}=\left\{\begin{pmatrix} 1\\ \mathbf{0}_d\\0\end{pmatrix}, \mathbf{0}_d\in\R^d\right\}\cup\left\{ \begin{pmatrix}\frac{\Vert\bm\nu\Vert^2}{2} \\ -\bm\nu\\ 1\end{pmatrix}, \bm\nu\in \R^d\right\}.
\end{equation*}}
\red{The advantage of describing  the invariant domain \eqref{Euler_invariant} by  \eqref{Euler_GQL} is the following: instead of one non linear relation to get the internal energy from the conserved variable, we have an infinite number of linear relations that will be shown, in section \ref{convex_limiting},  to be much easier to handle.}

The fluxes write
$$\bbf(\bbu)=\begin{pmatrix}\rho \bbv \\ \rho \bbv\otimes\bbv+p\text{Id}_{d} \\ (E+p)\bbv \end{pmatrix}$$
where $\text{Id}_{d}$ is the $d\times d$ identity matrix and we have introduced the pressure $p=p(\rho, \epsilon)$. In all the examples, the system is closed by the perfect gas equation of state, that is $p=(\gamma-1) e$, where the ratio of specific heats $\gamma$ is constant.

Other examples that will be considered {are the cases of scalar conservation laws with $\bbf$ linear taking the form of $\bbf(\bbu,\bbx)=\bba(\bbx)\bbu$ or nonlinear. In that case, $\bbu$ is a function with values in $\R$ and $\mathcal{D}=\R$ but the solution must stay in $\big[\min\limits_{\bbx\in \R^d}\bbu(t=0,\bbx), \max\limits_{\bbx\in \R^d}\bbu(t=0,\bbx)\big]$, due to Kruzhkov's theory. In these particular examples, we set the invariant domain as $\mathcal{D}=\big[\min\limits_{\bbx\in \R^d}\bbu(t=0,\bbx), \max\limits_{\bbx\in \R^d}\bbu(t=0,\bbx)\big]:=[\mathring{\bbu}_{\min},\mathring{\bbu}_{\max}]$. The following two cases are taken into account:} 
\begin{itemize}
\item A convection scalar problem
where $u\in \R$ and 
$$\dpar{u}{t}+\text{div }\big (\bba \; u)=0,$$
with the advection speed $\bba$ possibly depending on the spatial coordinate.
\item A nonlinear example
$$\dpar{u}{t}+\dpar{\big (\sin u\big )}{x}+\dpar{\big (\cos u\big )}{y}=0.
$$
\end{itemize}

The time evolution will be carried out by \red{means} of SSP Runge-Kutta schemes, so we only describe the first order forward Euler scheme. {We only describe the construction of the numerical schemes on general 2D polygons below and would like to use notations $f_x\equiv f_1$, $f_y\equiv f_2$, $A_x\equiv A_1$, and $A_y\equiv A_2$ instead.}
\subsection{High order schemes}\label{sec31} 
\red{To discretise in time, we use the method of lines with a SSP Runge-Kutta method. Since each step of this algorithm is obtained by a linear combination of Euler forward methods, we describe here only what we do for an Euler forward time stepping method.}

In our scheme, on a generic polygon $P$, the numerical solution $\bbu$ is represented by point values at the Gauss--Lobatto points of the edges of $P$ and the average. The update of the average is simply done applying the divergence theorem to \eqref{scheme:fv}:
\begin{equation*}
 	\vert P \vert \dfrac{{\rm d}\xbar\bbu_P}{{\rm d}t}+\oint_{\partial P}\bbf(\bbu)\cdot \bbn\; {\rm d}{\gamma}=0,
\end{equation*}
where $\bbn$ is the outward pointing unit normal at almost each point on the boundary $\partial P$. \red{Here, the symbol $\oint$ indicates that the volume integral is computed by means of a quadrature formula, the same convention is employed fr the boundary integrals.}Using Gauss--Lobatto quadrature rule to approximate the surface integral and the first order forward Euler method to discretize in time, we have
\begin{equation*}
  \xbar{\bbu}_P^{n+1}=\xbar{\bbu}_P^n-\frac{\dt}{\vert P\vert}\sum_{e \text{ edge of }{P}}\vert e\vert\; \Hbbf_{\bbn_e}(\bbu_P^n),\quad \bbu_P^n=\{\bbu_\sigma^n\}_{\sigma\in P},
\end{equation*}
where {$\vert e\vert$ is the measure of edge $e$ and}
\begin{equation}\label{flux}
  \Hbbf_{\bbn_e}(\bbu_P^n)=\sum_{\sigma\in e} \omega_\sigma \bbf(\bbu^n_\sigma)\cdot \bbn_e.
\end{equation}
Here, $\{\omega_\sigma\}$ are the weights of the Gauss--Lobatto points and  $\bbn_e$ is the outward normal unit to the edge $e$. 
Note that we have used the global continuity property of the approximation, and then no numerical flux is needed.

The update of the boundary values is more involved and we describe an extension of what \red{has} been proposed in several papers.
To update $\bbu_\sigma$, we consider a semi-discrete scheme of the form:
\begin{equation}\label{scheme:point}
  \dpar{\bbu_\sigma}{t}+\sum_{P, \sigma\in P}\bbPhi_\sigma^{P,\rm HO}(\bbu)=0,
\end{equation}
that again will be approximated by an Euler forward time stepping:
\begin{equation}\label{tfirstorder}
{\bbu_{\sigma_i}^{n+1}=\bbu_{\sigma_i}^n-\dt\sum_{P, \sigma_i\in P}\bbPhi_{\sigma_i}^{P, \rm HO}(\bbu).}
\end{equation}

\begin{figure}[!h]
\centerline{\subfigure[case of a vertex]{\includegraphics[trim=0.01cm 0.005cm 0.005cm 0.02cm,clip,width=4.5cm]{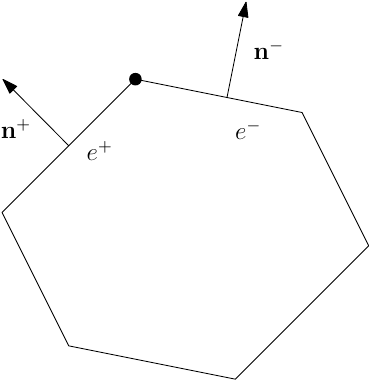}}\hspace*{2.5cm}\subfigure[case of a non vertex]{\includegraphics[trim=0.01cm 0.005cm 0.005cm 0.02cm,clip,width=4.5cm]{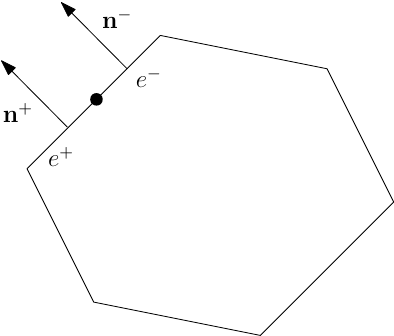}}}
\caption{\label{normals} Definitions of $e^\pm$ and $\bbn^\pm$.}
\end{figure}

The quantities $\bbPhi_\sigma^{P,{\rm HO}}(\bbu)$ are defined such that if the solution is linear and the problem linear with constant Jacobians. We have
\begin{equation*}
  	\sum_{P, \sigma\in P}\bbPhi_\sigma^{P,{\rm HO}}(\bbu)=\bbA\cdot \nabla \bbu(\sigma).
\end{equation*}	
There is still a lot of freedom in the definition of $\bbPhi_\sigma^{P,{\rm HO}}(\bbu)$. Inspired by Residual Distribution schemes, and in particular looking at the LDA scheme \cite{DeconinckMario}, we consider
\begin{subequations}
\label{residu:HO}
\begin{equation}
\label{residu_sigma}
\bbPhi_\sigma^{P,{\rm HO}}={N_\sigma} K_\sigma^+\bigg ( \bbA(\bbu_\sigma)\cdot\nabla \pi^\nabla\bbu(\sigma)\bigg ),
\end{equation}
where 
\begin{equation}
\label{N} N_\sigma^{-1}=\sum_{P,\sigma \in P}K^+_\sigma,
\end{equation}  
\end{subequations}
with the following definitions:
\begin{itemize}
\item $\bbA(\bbu_\sigma)$ is the vector of Jacobian matrices evaluated for the state $\bbu_\sigma$.
\item For any $\sigma\in P$, we define a {scaled} normal $\bbn_\sigma$ for the polygon $P$ and the DoFs $\sigma$ as follows: {if $\sigma$ is a vertex,} the vector $\bbn_\sigma$ is the sum of the {scaled outward} normals of the faces $e^+$ and $e^-$ sharing $\sigma$. {Otherwise, it is the half sum of the corresponding normals. We refer to} Figure \ref{normals} for a definition of the $e^\pm$ faces and normals $\bbn^\pm$. {To lighten the notation, we omit that these normals are referred to the element $P$, meaning that we have $\bbn_{\sigma,P}$.}
\item For any $\bbn_\sigma=(n_x,n_y)$, ${K_\sigma}:={\bbA(\bbu_\sigma)}\cdot \bbn_\sigma=A_xn_x+A_y n_y$ with $A_x=\dpar{{f_x}}{\bbu}$ and $A_y=\dpar{{f_y}}{\bbu}$.
\item Since the problem is hyperbolic, the matrix ${K_\sigma}$ is diagonalisable in $\R$, and we can take its positive part
\end{itemize}

Assuming that $\bbu$ is linear and the problem linear, we see that $\pi^\nabla \bbu=\bbu$ because $\bbu$ is linear and using the fact that $\bbA$ does not depend on $\bbu$, we have 
\begin{equation*}
\begin{split}
\sum_{P, \sigma\in P} \bbPhi_\sigma^{P,{\rm HO}}&=\sum_{P, \sigma\in P}N_\sigma K_\sigma^+\big ( \bbA\cdot\nabla \pi^\nabla\bbu(\sigma)\big )\\
& =N_\sigma \bigg ( \sum_{P, \sigma\in P} K_\sigma^+\bigg ) \; \bbA\cdot\nabla \bbu(\sigma)\\
& =\bbA\cdot\nabla \bbu(\sigma),
\end{split} \end{equation*}
if the condition \eqref{N} is met. The only thing is to show that the matrix $\sum_{{P},\sigma \in {P}}K^+_\sigma$ is invertible. It turns out that this is true for a hyperbolic system that is symmetrizable; see \cite{Abgrall99} for a proof. For the sake of completeness, we reproduce the proof and give precise assumptions in Appendix \ref{append:N}, as well as a precise statement.

We can also define the matrix $N$ from
\begin{equation*}
	N^{-1}=\sum_{P, \sigma\in P }\text{sign}\;{K_\sigma}.
\end{equation*}
This defines a different scheme, that provides (qualitatively) the same results.

The scheme is formally third order accurate in space because
{
\begin{equation*}
\begin{split}
\sum_{P, \sigma\in P} \bbPhi_\sigma^{P,{\rm HO}}-\bbA(\bbu_\sigma)\cdot \nabla \bbu(\bbx_\sigma)&=\sum_{P, \sigma\in P}\bbN_\sigma K_\sigma^+\bbA(\bbu_\sigma)\cdot \big [\nabla \pi^\nabla\bbu(\bbx_\sigma)-\nabla \bbu(\bbx_\sigma) \big ],
\end{split}
\end{equation*}
so that for any matrix norm,
\begin{equation*}
  \lVert \sum_{P, \sigma\in P} \bbPhi_\sigma^{P,{\rm HO}}-\bbA(\bbx_\sigma)\cdot \nabla \bbu(\bbx_\sigma)\rVert \leq \sum_{P, \sigma\in P} \lVert \bbN_\sigma K_\sigma^+\rVert \; {\lVert\bbA(\bbx_\sigma) \rVert}\;\lVert \nabla \pi^\nabla\bbu(\bbx_\sigma)-\nabla \bbu(\bbx_\sigma) \rVert,
\end{equation*}
and then
\begin{equation*}
 \sum_{P, \sigma\in P} \bbPhi_\sigma^{P,{\rm HO}}-\bbA(\bbx_\sigma)\cdot \nabla \bbu(\bbx_\sigma)=\red{O(h^2)}.
\end{equation*}
}
{The exponant $2$ comes from the following.
It is easy, at least for convex polygons, using the technique of \cite{CiarletRaviart}, to show that
for $\bbu\in C^{k+1}(K)$, the VEM approximation $\pi_h(\bbu)$ satisfies, in the $L^q(K)$ norm, $0<q\leq +\infty$, that
for any $\bbx\in K$ that
$$\Vert D^p\pi_h(\bbu)-D^p\bbu(\bbx)\Vert_{q}\leq C(K)h^{k+1-p}\Vert D^{p+1}\bbu\Vert_q$$
where the constant $C(K)$ depends on the polygon $K$. We have denoted the $p$-th derivative of $\bbu$ by $D^p\bbu$\footnote{\red{so that for $p=1$, $D^1\bbu=\nabla \bbu$.}}. We believe that, using the shape regularity assumptions of \cite{chinois}, this constant is independent of $K$, provided that $K$ belongs to the class defined in \cite{chinois}. In practice, this means that the polygons have a "nice" aspect ratio.

Now, in practice, it seems that we have better than second order accuracy. This is shown in the numerical section. The explanation of this fact, at least for linear hyperbolic problems, is given in the section 3 of  \cite{AbgrallOeffnerLiuDG}, where a link with the discontinous Galerkin method is made.  For fairness, we must mention the reference \cite{WassilijPetrov} published independently  on Arxiv, 2 days before.  Improving this point will be the topic of a future publication.
}

\medskip
{
It turns out that the scheme \eqref{scheme:point}-\eqref{residu:HO} is not fully satisfactory. Simulations on the vortex problem described in Section \ref{vortex}, with the scheme for point values defined by \eqref{residu:HO},   show that spurious modes exist. They are not damped and  do not amplify. We interpret this as a loss of information while going from $\bbu$ to $\pi^\nabla\bbu$: the dimension of $V_k(P)$ is always larger than that of $\P^k(P)$ for any $P$ and any $k$. This problem does not exist with  the approximation described in \cite{abgrall2023activefluxtriangularmeshes}: we do not need any projector on triangular meshes, since the gradient can be explicitly computed. Because of that fact, we need to 
modify \eqref{residu:HO} by adding a term that will damp out the spurious modes that are suspected to come from the mismatch between the VEM approximation $\bbu$ and its projection $\pi^\nabla\bbu$. This can be achieved if one add to \eqref{residu:HO} a term that is strictly dissipative when $\bbu-\pi^\nabla\bbu\neq 0$.

Inspired by \cite{hitch}, such a term can be 
\begin{equation}
\label{scheme:stab}\mathfrak{D}_\sigma=\frac{\alpha_{P}}{\sqrt{h}}\sum_{r=1}^{N_{{\text{DoFs}}}} \text{ DoF}_r(\bbu-\pi^\nabla \bbu)\text{ DoF}_r(\varphi_\sigma-\pi^\nabla\varphi_\sigma),
\end{equation}
where $\alpha_{P}$ {is} the spectral radius of $\bbA(\bbu)\cdot \bbn$ {(i.e., $A_xn_x+A_yn_y$)} in $P$.  In \cite{hitch}, this is the term that is added to the approximation 
of $$\int_P\nabla \bbu\cdot \nabla \varphi_\sigma\; d\bbx$$ by
$$\int_P\nabla \big (\pi^\nabla\bbu\big )\cdot \nabla \big (\pi^\nabla \varphi_\sigma\big )\; d\bbx,$$
i.e.
$$\int_P\nabla \bbu\cdot \nabla \varphi_\sigma\; d\bbx\approx \int_P\nabla \pi^\nabla\bbu\cdot \nabla\big ( \pi^\nabla\varphi_\sigma\big ) \; d\bbx+
\sum_{r=1}^{N_{{\text{DoFs}}}} \text{ DoF}_r(\bbu-\pi^\nabla \bbu)\text{ DoF}_r(\varphi_\sigma-\pi^\nabla\varphi_\sigma).$$

We note that  for all $\sigma$ on the boundary,  we have 
$$\int_{P} \varphi_\sigma\; {\rm d}\bbx=\int_{P}\pi^\nabla \varphi_\sigma\; {\rm d}\bbx=0,$$
so that 
\begin{equation*}
  \mathfrak{D}_\sigma=\frac{\alpha_P}{\sqrt{h}} \sum_{r{\in \{1,2,\cdots,N_{{\text{DoFs}}}\}\backslash\{\iota\}}}\text{ DoF}_r(\bbu-\pi^\nabla \bbu)\text{ DoF}_r(\varphi_\sigma-\pi^\nabla \varphi_\sigma),
\end{equation*}
{where $\iota$ corresponds to the average degree of freedom.}

We also note that
$$\sum_{\sigma\in \partial P}\text{DoF}_\sigma(\bbu)=\sum_{r{\in \{1,2,\cdots,N_{{\text{DoFs}}}\}\backslash\{\iota\}}}\text{ DoF}_r(\bbu-\pi^\nabla u)^2>0,$$
if $\bbu$ is not a polynomial. 
The residual for the point values \eqref{residu_sigma} becomes
\begin{equation}
\label{residu_sigma_2}
\bbPhi_\sigma^{P,{\rm HO}}={N_\sigma} K_\sigma^+\bigg ( \bbA(\bbu_\sigma)\cdot\nabla \pi^\nabla\bbu(\sigma)\bigg )+\mathfrak{D}_\sigma.
\end{equation}
We also have that $\text{ DoF}_r(\bbu-\pi^\nabla u)=O(h^{k+1})$, and the same applies for the $\varphi_\sigma$, so that let us notice that $\mathfrak{D}_\sigma=O(h^{2k+1/2})$: the stabilization term does not spoil the accuracy.
}
\subsection{Low order schemes}\label{sec32}
The update of the average value is carried out as follows:
\begin{equation}\label{low:average}
\xbar\bbu^{n+1}_{P}=\xbar \bbu_{P}^n-\frac{\Delta t}{\vert {P}\vert} \sum_{e \text{ edge of }{P}}\vert e\vert \lbbf_{\bbn_e}(\xbar \bbu_{P}^n, \xbar \bbu_{P^-}^n),
\end{equation}
where ${P}^-$ is the polygon sitting on the other side of $e$. If $\lbbf_{\bbn_e}$ is a monotone flux, then, the scheme will be stable. In the numerical implementation, we have used the first order Local Lax--Friedrichs (or Rusanov) flux:
\begin{equation}\label{flux:LO}
\lbbf_{\bbn_e}(\xbar \bbu_{P}^n, \xbar \bbu_{P^-}^n)=\frac{\big(\bbf(\xbar{\bbu}_P^n)+\bbf(\xbar{\bbu}_{P^-}^n)\big)\cdot\bbn_e}{2}-\frac{\alpha_e}{2}(\xbar{\bbu}_{P^-}^n-\xbar{\bbu}_P^n),
\end{equation}
where $\alpha_e=\alpha_e(\xbar{\bbu}_{P^-}^n,\xbar{\bbu}_P^n,\bbn_e)$ is the maximum speed obtained from the Riemann problem between the states $\xbar{\bbu}_P^n$ and $\xbar{\bbu}_{P^-}^n$ in the direction $\bbn_e$ as evaluated in \cite{GuermondPopovFast}. It can be shown that under condition $\alpha_e\geq\max_{\bbw\in I(\xbar{\bbu}_{P^-}^n,\xbar{\bbu}_P^n)}\big(\vert\bbf'(\bbw)\cdot\bbn_e\vert\big)$, where $I(a,b)=[\min(a,b),\max(a,b)]$, the numerical flux \eqref{flux:LO} is a monotone flux.
The main drawback of this approach is that there is no longer any coupling between the point values and the average values in the update \eqref{low:average}. This might be a problem, but we have not found any concrete example where this approach fails.

\begin{figure}[!h]
\begin{center}
\subfigure[Sub-triangulation]{\includegraphics[trim=0.01cm 0.005cm 0.01cm 0.02cm,clip,width=5.5cm]{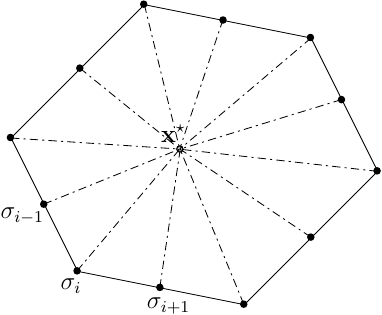}}\hspace*{2.5cm}
\subfigure[Normals]{\includegraphics[trim=0.01cm 0.06cm 0.01cm 0.02cm,clip,width=4.5cm]{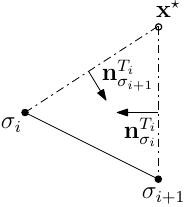}}
\end{center}
\caption{\label{subtriangulation} Sub-triangulation and normals that used for the low order scheme.}
\end{figure}

Again, the update of the point values is a bit more subtle. For notations and graphical illustration, we refer to Figure \ref{subtriangulation}.
By assumption, there exists a point, ${\bbx^\star}$, such that $P$ is star-shaped with respect to this point. In practice, we have always taken the centroid, because all the polygons we have considered are convex. Then, as drawn in Figure \ref{subtriangulation}--(a), we connect this point to the DoFs on $\partial P$, and this creates a sub-triangulation of $P$. Taking a {counter-clockwise} orientation of $\partial P$, we denote the DoFs on $\partial P$ as $\{\sigma_i\}_{i=1, \ldots , N_P}$ with $N_P$ representing the number of point values DoFs and $\sigma_{N_P+1}=\sigma_1$, so that the sub-triangles will be $T_i=\{\sigma_i,\sigma_{i+1},\bbx^\star\}$ for $i=1, \ldots, N_P$. The vertex $\sigma_i$ is shared by the triangles $\{\sigma_i, \sigma_{i+1},\bbx^\star\}$ and $\{\sigma_{i-1},\sigma_{i},\bbx^\star\}$ and the list of sub-triangles in $P$ is denoted by $\mathcal{T}_{P}$. {For each sub-element $T_i$, we can also get the scaled normals for the DoFs $\sigma_i$ and $\sigma_{i+1}$ as shown in Figure \ref{subtriangulation}--(b).} We identify the average value $\xbar \bbu_P$ with an approximation of $\bbu$ at $\bbx^\star$. This has no impact on the accuracy since we are looking for a first order scheme.

The forward Euler update of the $i$-th point value {$\bbu_{\sigma_i}$ is given by} \eqref{tfirstorder}. Again inspired by what has been done in the framework of Residual Distribution schemes, the so-called residuals in \eqref{tfirstorder} are given by
\begin{equation}\label{low_residu}
{\bbPhi_{\sigma_i}^{P,{\rm LO}}=\sum\limits_{T_j\in \mathcal{T}_P, \sigma_i\in T_j}  \Psi_{\sigma_i}^{T_j}=\Psi_{\sigma_i}^{T_{i-1}}+\Psi_{\sigma_i}^{T_i}},
\end{equation} with
\begin{subequations}
\begin{equation}\label{low_residua}
{
\begin{aligned}
\vert C_\sigmai\vert\Psi_{\sigma_i}^{T_{i-1}}&=\frac{1}{3}\ \oint_{T_{i-1}} \bbA\cdot \nabla\bbu\; {\rm d}\bbx+\alpha_{T_{i-1}}\big ( \bbu_{\sigma_i}^n-\xbar{\bbu}_{T_{i-1}}\big)\\
&=\frac{1}{6}\bigg ( \big (\bbf(\bbu_{\sigma_{i-1}}^n)-\bbf(\xbar{\bbu}_P^n)\big )\cdot \bbn_{\sigma_{i-1}}^{T_{i-1}}+
  \big ( \bbf(\bbu_{\sigma_i}^n)-\bbf(\xbar{\bbu}_P^n)\big )\cdot \bbn_{\sigma_i}^{T_{i-1}}\bigg )+\frac{\alpha_{\sigma_{i}}}{3}\sum_{\substack{j\in T_{i-1}\\j\neq i}}\big ( \bbu_{\sigma_i}^n-\bbu_{\sigma_j}^n\big)
\end{aligned}}
\end{equation}
and
{
\begin{equation}\label{low_residub}
\begin{aligned}
\vert C_\sigmai\vert\Psi_{\sigma_i}^{T_{i}}&=\frac{1}{3}\ \oint_{T_{i}} \bbA\cdot \nabla\bbu\; {\rm d}\bbx+\alpha_{T_{i}}\big ( \bbu_{\sigma_i}^n-\xbar{\bbu}_{T_{i}}\big)\\
&=\frac{1}{6}\bigg ( \big (\bbf(\bbu_{\sigma_{i}}^n)-\bbf(\xbar{\bbu}_P^n)\big )\cdot \bbn_{\sigma_{i}}^{T_{i}}+
  \big ( \bbf(\bbu_{\sigma_{i+1}}^n)-\bbf(\xbar{\bbu}_P^n)\big )\cdot \bbn_{\sigma_{i+1}}^{T_{i}}\bigg )+\frac{\alpha_{\sigma_i}}{3}\sum_{\substack{j\in T_{i}\\j\neq i}}\big ( \bbu_{\sigma_i}^n-\bbu_{\sigma_j}^n\big),
\end{aligned}
\end{equation}
where {$C_{\sigma_i}$} is the dual cell of area
\begin{equation}\label{area}
{\vert C_{\sigma_i}\vert =\sum\limits_{P, \sigma_i\in P}\sum\limits_{T_j\in \mathcal{T}_P, \sigma_i\in T_j} \frac{\vert T_j\vert }{3}\xlongequal[]{\text{Figure \ref{subtriangulation}--(a)}}\sum\limits_{P, \sigma_i\in P}\frac{1}{3}\big(\vert T_{i-1}\vert+\vert T_i\vert\big)}, 
\end{equation}$\alpha_{\sigmai}=\max(\alpha_{T_i},\alpha_{T_{i-1}})$ and $\alpha_{T_j}$ is an upper bound of the spectral radius of $\bbA$ evaluated from the two point values of $T_j$ and the average value $\xbar \bbu_P$ multiplied by an upper bound of the scaled normals. The scheme is stable under} {a CFL like condition of the form}
\begin{equation}
	{\Delta t \leq \text{CFL} \, \min_\Omega \frac{\vert C_{\sigma_j}\vert }{\vert \partial C_{\sigma_j} \vert \, \vert \alpha_{T} \vert},}
	\label{eqn.cfl}
\end{equation}
{with $\vert \partial C_{\sigma_j} \vert$ being the perimeter of the associated dual cell and $\text{CFL}<1$.}
\end{subequations}


\subsection{{Limiting strategies}}
It is well known that high order scheme will not be oscillation free without any form of {limiting strategies}. When the solution develops discontinuities, these schemes will be prone to numerical oscillations and possibly non-physical solutions. \red{For instance, in gas dynamics, the  density or  the internal energy, or both  might become negative.} To address this issue, {limiting} is required. Below, we introduce two techniques.

\subsubsection{{A-posteriori limiting strategy}}
{A straightforward nonlinear stabilization method is to use the {a-posteriori limiting} as it was done in \cite{Abgrall_AF,abgrall2023activefluxtriangularmeshes, Liu2024,Abgrall2024_WBAF}.} Since point values and averages are independent variables, we need to test both against a set of admissibility criteria. The idea is to work with several schemes ranging from order $k={3}$ to $k=1$, with the lowest order one able to provide results staying in the invariance domain $\mathcal{D}$. For the element $P$, we write the scheme for $\xbar\bbu_P$ as $S_P(k)$ and for $\bbu_\sigma$ as $S_\sigma(k)$, namely for average and point values, respectively.

We denote by $\xbar\bbu^n_P$ and $\bbu_\sigma^n$ the solution at the time $t_n$. After each Runge--Kutta cycle, the updated solution is denoted by $\widetilde{\xbar\bbu}_P$, $\widetilde\bbu_\sigma$. 

We first run the high order scheme, for each Runge-Kutta cycle.  Concerning the average $\xbar\bbu_P$, we store the flux {$\oint_e\bbf(\bbu)\cdot\bbn \; {\rm d}\gamma$} for reasons of conservation. Then, for the average solution, we proceed as follows:
\begin{enumerate}
\item {Computer Admissible Detector (CAD):} we check if $\widetilde{\xbar\bbu}_P$ is a valid vector with real components, namely we verify that each component is not NaN. If this is not the case, we flag the element and go to the next one in the list.
\item {Physically Admissible Detector (PAD):} we check if $\widetilde{\xbar\bbu}_P\in \mathcal{D}$. If this is false, the element is flagged, and we proceed the control on the next element.
\item 
Then we check if at $t^n$, the solution is not constant in the elements used in the numerical stencil (so we check and compare between the average and point values $\bbu_\sigma$ in $P$ with those in the elements sharing a face with $P$). This is done  in order to avoid to detect a wrong maximum principle. If the solution is declared locally constant {up to an error of $10^{-10}$}, we move to the next element and we let the polygon $P$ unflagged.
\item Discrete Maximum Principle (DMP): we check if $\widetilde{\xbar\bbu}_P$ is a local extremum. If we are dealing with the Euler equations, we compute the density and the pressure and perform this test on these two quantities only, even though for a system this is not really meaningful. We denote by $\xi$ the variable on which we perform the test (i.e. $\bbu$ itself for scalar problems, and the density/pressure for the Euler system). Let $\mathcal{V}(P)$ be the set of elements $F$ that share a face or a vertex with $P$, excluding $P$ itself. We say that we have a potential extremum if 
$${\xi}_P^{n+1}\not\in \big[ \min\limits_{F\in \mathcal{V}(P)}\xi_F^n{-}\varepsilon_P^n, \max\limits_{F\in \mathcal{V}(P)}\xi_F^n{+}\varepsilon_P^n\big],$$
where {$\varepsilon_P^n=\max\Big(10^{-4},10^{-3}\big(\max\limits_{\text{all~}P}\xi_P^n-\min\limits_{\text{all~}P}\xi_P^n\big)\Big)$}.
If the above test is true, the element is flagged. 
\end{enumerate}
If an element is flagged, then each of its faces are flagged, and we recompute the flux of the flagged faces using the first order scheme. 

For the point values, the procedure is similar, and then for the flagged DoFs, we recompute the residuals
$$\sum\limits_{P, \sigma\in P}\bbPhi_\sigma^P(\bbu^n)$$ with the first order scheme.

\subsubsection{A monolithic convex limiting method}\label{convex_limiting}
In an alternative direction, we extend the convex limiting method developed for the one-dimensional case in \cite{Abgrall_BP_PAMPA} to the two-dimensional case. In both cases, the source of inspiration is \cite{Vilar_DGFV}. Our goal is to find a convex limiting approach between the high-order and low-order (parachute) fluxes, operators, or schemes so that the resulted method is invariant domain preserving (IDP). For simplicity and clarity, we first illustrate the blending procedure with the focus on the scalar conservation laws. The extension to Euler equations will be addressed later. To simplify the text, we do not consider boundary conditions, which will be described afterwards.

\paragraph{Evolving  cell averages.} Considering $P$, it may have two types of faces: the faces that are interior to $\Omega$, and those which intersect $\Gamma=\partial \Omega$. In the following,  $\FF_P$ is the set of faces interior to $\Omega$, and the set of boundary faces is $\FF_P^b$. We may have $\FF_P^b=\emptyset$. We can rewrite the forward Euler update of cell averages as follows:
\begin{subequations}\label{ave:scheme}
\begin{equation}\label{ave:scheme:1}
    \xbar{\bbu}_P^{n+1}=\xbar{\bbu}_P^n-\frac{\dt}{\vert P\vert}\bigg ( \sum_{e\in \FF_P }\vert e\vert \; \hbbf_{\bbn_e}({\bbu}_P^n, \xbar{\bbu}_P^n,\xbar{\bbu}_{P^-}^n)
    \bigg ),
\end{equation}
where the flux $\hbbf_{\bbn_e}$ is written as
\begin{equation}\label{ave:scheme:2}
\hbbf_{\bbn_e}(\bbu_P^n,\xbar{\bbu}_P^n,\xbar{\bbu}_{P^-}^n)=\lbbf_{\bbn_e}(\xbar{\bbu}_P^n,\xbar{\bbu}_{P^-}^n)+\eta_{e}\Delta \hbbf_{\bbn_e}(\bbu_P^n,\xbar{\bbu}_P^n,\xbar{\bbu}_{P^-}^n)
\end{equation}
with $\eta_e\in[0,1]$, $\lbbf_{\bbn_e}(\xbar{\bbu}_P^n,\xbar{\bbu}_{P^-}^n)$ given by \eqref{flux:LO}, and
\begin{equation}\label{ave:dflux}
  \Delta \hbbf_{\bbn_e}(\bbu_P^n,\xbar{\bbu}_P^n,\xbar{\bbu}_{P^-}^n)=\Hbbf_{\bbn_e}(\bbu_P^n)-\lbbf_{\bbn_e}(\xbar{\bbu}_P^n,\xbar{\bbu}_{P^-}^n),
\end{equation}
where $\Hbbf_{\bbn_e}(\bbu_P^n)$ is given by \eqref{flux}. 
\end{subequations}

Next, we rewrite the Euler forward update as a convex combination of quantities defined at the previous time step:
\begin{equation}\label{eq:ave_convex}
\begin{split}
\xbar{\bbu}_P^{n+1}&=\xbar{\bbu}_P^n-\frac{\dt}{\vert P\vert}\bigg ( \sum_{e\in \FF_P }\vert e\vert \; \hbbf_{\bbn_e}(\bbu_P^n,\xbar{\bbu}_P^n,\xbar{\bbu}_{P^-}^n)
\bigg )\\
 &=\xbar{\bbu}_P^n-\frac{\dt}{\vert P\vert }\sum_{e\in \FF_P}\vert e\vert 
 \hbbf_{\bbn_e}(\bbu_P^n,\xbar{\bbu}_P^n, \xbar{\bbu}_{P^-}^n)+\frac{\dt}{\vert P\vert} \sum_{e\in \FF_P} \vert e\vert\big ( \alpha_{e}\xbar{\bbu}_P^n - \alpha_{e}\xbar{\bbu}_{P}^n+\bbf(\xbar{\bbu}_P^n)\cdot \bbn_e \big )\\
 %
 &=\Bigg ( 1-\frac{\dt}{\vert P\vert }\Big (\sum_{e\in \FF^P}\vert e\vert \alpha_{e}
 \Big )\Bigg ) \xbar{\bbu}_P^n+
 \frac{\dt}{\vert P\vert} \sum_{e\in \FF_P} \vert e\vert \alpha_{e}\bigg ( \xbar{\bbu}_P^n-\dfrac{\hbbf_{\bbn_e} (\bbu_P^n,\xbar{\bbu}_P^n,\xbar{\bbu}_{P^-}^n)-\bbf(\xbar{\bbu}_P^n)\cdot \bbn_e}{\alpha_{e}}\bigg ).
\end{split}
\end{equation}
Then, defining the blended Riemann intermediate states
\begin{equation*}
  \widetilde{\xbar{\bbu}}_P^e=\xbar{\bbu}_P^n-\dfrac{ \hbbf_{\bbn_e} (\bbu_P^n,\xbar{\bbu}_P^n,\xbar{\bbu}_{P^-}^n)-\bbf(\xbar{\bbu}_P^n)\cdot \bbn_e}{\alpha_{e}},
\end{equation*}
the previous expression of forward Euler update can finally be recast into the following convex form:
\begin{equation}\label{ave:convex_form}
  \xbar{\bbu}_P^{n+1}=\Bigg ( 1-\frac{\dt}{\vert P\vert }\Big (\sum_{e\in \FF^P}\vert e\vert \alpha_{e}
  \Big )\Bigg ) \xbar{\bbu}_P^n
 +\frac{\dt}{\vert P\vert} \sum_{e\in \FF_P} \vert e\vert \alpha_{e}\widetilde{\xbar{\bbu}}_P^e
\end{equation}
provided that the standard CFL condition 
\begin{equation*}
\frac{\dt}{\vert P\vert }\bigg ( \sum_{e\in \FF_P}\vert e\vert \alpha_{e}
\bigg )\leq 1
\end{equation*}
is satisfied.

By means of the blended numerical flux in \eqref{ave:scheme:2}, \eqref{ave:dflux}, and the low order monotone numerical flux given by \eqref{flux:LO}, the blended Riemann intermediate state of the interior part $\widetilde{\xbar{\bbu}}_P^e$ can be rewritten into the following form:
\begin{equation*}
  \widetilde{\xbar{\bbu}}_P^e=\xbar{\bbu}_P^n-\dfrac{ \hbbf_{\bbn_e} (\bbu_P^n,\xbar{\bbu}_P^n,\xbar{\bbu}_{P^-}^n)-\bbf(\xbar{\bbu}_P^n)\cdot \bbn_e}{\alpha_{e}}
  =\xbar{\bbu}_P^{e,\star}-\eta_e\frac{\Delta\hbbf_{\bbn_e}(\bbu_P^n,\xbar{\bbu}_P^n,\xbar{\bbu}_{P^-}^n)}{\alpha_e},
\end{equation*}
where $\xbar{\bbu}_P^{e,\star}$ is nothing but the first order finite volume Riemann intermediate state given by
\begin{equation*}
  \xbar{\bbu}_P^{e,\star}=\frac{\xbar{\bbu}_P^n+\xbar{\bbu}_{P^-}^n}{2}-\frac{\big(\bbf(\xbar{\bbu}_{P^-}^n)-\bbf(\xbar{\bbu}_P^n)\big)\cdot\bbn_e}{2\alpha_e}.
\end{equation*}
Since $\alpha_e\geq\max_{\bbw\in I(\xbar{\bbu}_{P^-}^n,\xbar{\bbu}_{P}^n)}\big(\vert\bbf'(\bbw)\cdot\bbn_e\vert\big)$, it follows that $\xbar{\bbu}_P^{e,\star}\in I(\xbar{\bbu}_{P^-}^n,\xbar{\bbu}_{P}^n)\subset\mathcal{D}$. 

We now introduce the definition of the blending coefficient $\eta_e$ to ensure that if the numerical initial solution for scalar conservation laws lies in $\mathcal{D}=[\mathring{\bbu}_{\min},\mathring{\bbu}_{\max}]$, the solution $\xbar{u}_P$ remains in $\mathcal{D}$ throughout the entire calculation. To guarantee this property, it is sufficient, by convexity of relation \eqref{ave:convex_form}, that the blended Riemann intermediate state $\widetilde{\xbar{\bbu}}_P^e$ also remains in $\mathcal{D}=[\mathring{\bbu}_{\min},\mathring{\bbu}_{\max}]$. Since $\xbar{\bbu}_P^{e,\star}$ does, a sufficient condition is then to set $\eta_e$ such that
\begin{equation}\label{eq:theta_average}
  \eta_e\leq\min\Big(1,\frac{\alpha_e}{\big\vert\Delta\hbbf_{\bbn_e}(\bbu_P^n,\xbar{\bbu}_P^n,\xbar{\bbu}_{P^-}^n)\big\vert}\min\big(\mathring{\bbu}_{\max}-\xbar{\bbu}_P^{e,\star},
  \xbar{\bbu}_P^{e,\star}-\mathring{\bbu}_{\min}\big)\Big).
\end{equation}

\paragraph{Evolving point values.} We blend the residuals for the update of any point values $\bbu_\sigmai$ as
\begin{equation}\label{point:blend_residu}
  \bbPhi_\sigmai^P=\bbPhi_\sigmai^{P,{\rm LO}}+\theta_\sigmai\Delta\bbPhi_\sigmai,
\end{equation}
where
\begin{equation}\label{point:dresidu}
  \Delta\bbPhi_\sigmai=\bbPhi_\sigmai^{P,{\rm HO}}-\bbPhi_\sigmai^{P,{\rm LO}},
\end{equation}
with $\bbPhi_\sigmai^{P,{\rm HO}}$ extracted from \eqref{residu_sigma}, and $\bbPhi_\sigmai^{P,{\rm LO}}$ given by \eqref{low_residu}, \eqref{low_residua}--\eqref{low_residub}. 
The forward Euler time integration for the point value $\bbu_\sigmai$ reads as
\begin{equation}\label{scheme:residual:o}
\begin{split}
\bbu_\sigmai^{n+1}&=\bbu_\sigmai^n-\dt\sum_{\substack{P, \sigmai \in P
}} \bbPhi_\sigmai^P(\bbu)
=\bbu_\sigmai^n-\dt\sum_{\substack{P, \sigmai \in P
}}\bigg( \bbPhi_\sigmai^P(\bbu)+\sum_{T_j\in \mathcal{T}_P, \sigma_i\in T_j}\frac{\alpha_{\sigma_i}}{\vert C_\sigmai\vert}\big(\bbu_\sigmai^n-\bbu_\sigmai^n\big)\bigg)\\
&=\Bigg (1-\dt \bigg ( \frac{1}{\vert C_\sigmai\vert}\sum_{\substack{P, \sigmai \in P
}}\sum_{T_j\in \mathcal{T}_P, \sigma_i\notin T_j}\alpha_{\sigma_i}
 \bigg )\Bigg ) \bbu_\sigmai^n  
+\frac{\dt}{\vert C_\sigmai\vert}\sum_{\substack{P, \sigmai \in P\\ \sigmai\notin\Gamma}}\alpha_{\sigma_i} \Big( \sum_{T_j\in \mathcal{T}_P, \sigma_i\in T_j} \bbu_\sigmai^n-\vert C_\sigmai\vert\frac{\bbPhi_\sigmai^P}{\alpha_{\sigma_i} }\Big).
\end{split}
\end{equation}
Then, defining the blended Riemann intermediate states
\begin{equation*}
  \widetilde{\bbu}_\sigmai=\sum_{T_j\in \mathcal{T}_P, \sigma_i\in T_j}\bbu_\sigmai^n-\vert C_\sigmai\vert\frac{\bbPhi_\sigmai^P}{\alpha_\sigmai},
\end{equation*}
the previous expression of forward Euler update can finally be recast into the following convex form:
\begin{equation*}
\begin{split}
  {\bbu}_\sigmai^{n+1}&=\Bigg (1-\dt \bigg ( \frac{1}{\vert C_\sigmai\vert}\sum_{\substack{P, \sigmai \in P\\ \sigmai\notin\Gamma}}\sum_{T_j\in \mathcal{T}_P, \sigma_i\notin T_j}\alpha_{\sigma_i}
 \bigg )\Bigg ) \bbu_\sigmai^n
 +\frac{\dt}{\vert C_\sigmai\vert}\sum_{\substack{P, \sigmai \in P
 }} \alpha_\sigmai \widetilde{\bbu}_\sigmai,
\end{split}
\end{equation*}
provided that the  CFL condition 
\begin{equation*}
\dt \bigg ( \frac{1}{\vert C_\sigmai\vert}\sum_{\substack{P, \sigmai \in P
}}\sum_{T_j\in \mathcal{T}_P, \sigma_i\notin T_j}\alpha_{\sigma_i}
 \bigg )\leq 1
\end{equation*}
is satisfied.

Using the blended residuals in \eqref{point:blend_residu}--\eqref{point:dresidu}, and the low order residuals given in \eqref{low_residu}, \eqref{low_residua}--\eqref{low_residub}, the blended Riemann intermediate state of the inner point values $\widetilde{\bbu}_\sigmai$ can be rewritten into the following form:
\begin{equation*}
 \widetilde{\bbu}_\sigmai=\sum_{T_j\in \mathcal{T}_P, \sigma_i\in T_j} \bbu_\sigmai^n-\vert C_\sigmai\vert\frac{\bbPhi_\sigmai^{P,{\rm LO}}}{\alpha_{\sigma_i}}-\theta_\sigmai\vert C_\sigmai\vert\frac{\Delta\bbPhi_\sigmai}{\alpha_{\sigma_i}}=\bbu_\sigmai^{P,\star}-\theta_\sigmai\vert C_\sigmai\vert\frac{\Delta\bbPhi_\sigmai}{\alpha_{\sigma_i}},
\end{equation*}
where 
\begin{equation*}
\begin{split}
  \bbu_\sigmai^{P,\star}&= \frac{1}{3}\Big(\frac{\bbu^n_{\sigma_{i-1}}+\bbu^n_{\sigma_{i+1}}}{2}-\frac{(\bbf(\bbu^n_{\sigma_{i-1}})-\bbf(\bbu^n_{\sigma_{i+1}}))\bbn_{\sigma_{i-1}}^{T_{i-1}}}{2\alpha_\sigmai}
  +\bbu^n_\sigmai+\xbar{\bbu}_P^n-\frac{(\bbf(\bbu^n_\sigmai)-\bbf(\xbar{\bbu}_P^n))}{2\alpha_\sigmai}\big(\bbn_{\sigmai}^{T_{i-1}}+\bbn_{\sigmai}^{T_{i}}\big)\Big)\\
  &\quad +\frac{1}{6}(\bbu^n_{\sigma_{i-1}}+2\bbu^n_{\sigma_i}+2\xbar{\bbu}_P^n+\bbu^n_{\sigma_{i+1}}) 
\end{split}
\end{equation*}
By using the fact that $\bbu^n_{\sigma_{i-1}}, \bbu^n_{\sigma_i}, \bbu^n_{\sigma_{i+1}}, \xbar{\bbu}_P^n \in \mathcal{D}$ and the choice of $\alpha_\sigmai$ in Section \ref{sec32}, we can verify that ${\bbu}_\sigmai^{P,\star}\in \mathcal{D}$. 
 Analogously, to ensure that the solution $\bbu_\sigmai^{n+1}$ remains in the convex invariant domain $\mathcal{D}$, it is sufficient to take the residuals blending coefficient $\theta_\sigmai$ as
\begin{equation}\label{eq:theta_point}
  \theta_\sigmai\leq\min\Big(1,\frac{\alpha_\sigmai}{\vert C_\sigmai\vert\vert\Delta\bbPhi_\sigmai\vert}\min\big(\mathring{\bbu}_{\max}-\bbu_\sigmai^{P,\star},\bbu_\sigmai^{P,\star}-\mathring{\bbu}_{\min}\big)\Big).
\end{equation}

{
\begin{remark}
We have demonstrated the convex limiting approach based on the forward Euler time-stepping method. When extended to the third-order strong stability preserving Runge-Kutta (SSP--RK) method, the key properties are retained, as these high-order time integration scheme can be expressed as convex combinations of several forward Euler steps.

However, one important consideration in practical implementation is the choice of the adaptive time step. When the time step is determined using the standard CFL condition and the solution from the previous time level, it guarantees the validity of the convex combinations in \eqref{eq:ave_convex} and \eqref{scheme:residual:o} only for the first RK stage. If the solution leaves the convex invariant domain during subsequent RK stages, it is typically due to this adaptive time step being too large so that the convex combinations in \eqref{eq:ave_convex} and \eqref{scheme:residual:o} are not validated.

A simple remedy is to recompute the time step using the standard CFL condition, but based on the current solution. With this updated time step, the semi-discrete system is then re-evolved starting from the first RK stage.  

In practical applications, this procedure was never needed, since there were never dramatic change of the time step during one Runge Kutta cycle.
\end{remark}}

\paragraph{Extension to Euler equations.}  The goal here is to define the blending coefficients for the system of Euler equations so that the solution remains in the convex invariant domain $\mathcal{D}_{\bm \nu}$ defined by \eqref{Euler_GQL}. To simplify the text, we avoid to consider boundary conditions. This can be done in the same way as for the scalar case, and details will be given in Section \ref{sec:BCs}.

According to \eqref{eq:theta_average} and \eqref{eq:theta_point}, to guaranty the positivity of the density $\rho$, it is enough to take
\begin{equation*}
  \eta_e^{\rho}\leq\min\Big(1,\frac{\alpha_e \xbar{\rho}_P^{e,\star}}{\big\vert\Delta\hbbf_{\bbn_e}^{\rho}(\bbu_P^n,\xbar{\bbu}_P^n,\xbar{\bbu}_{P^-}^n)\big\vert}\Big)
\end{equation*} 
and
\begin{equation*}
  \theta_\sigmai^\rho\leq\min\Big(1,\frac{\alpha_\sigmai\rho_\sigmai^{P,\star}}{\vert C_\sigmai\vert\vert\Delta\bbPhi_\sigmai\vert}\Big).
\end{equation*}

Note that the constraint in the GQL representation \eqref{Euler_GQL} is linear with respect to $\bbu$. To guarantee the positivity of the internal energy, we therefore only need
\begin{equation*}
  \xbar{\bbu}_P^{e,\star}\cdot\bbn_{*}\pm\eta_e\frac{\Delta\hbbf_{\bbn_e}(\bbu_P^n,\xbar{\bbu}_P^n,\xbar{\bbu}_{P^-}^n)\cdot\bbn_*}{\alpha_e}>0
\end{equation*}
and
\begin{equation*}
  \bbu_\sigmai^{P,\star}\cdot\bbn_{*}-\theta_\sigmai\vert C_\sigmai\vert\frac{\Delta\bbPhi_\sigmai \cdot\bbn_*}{\alpha_{\sigma_i}}>0,
\end{equation*}
where $\bbn_*=(\frac{\Vert \bm\nu\Vert^2}{2}, -\bm\nu, 1)^\top$, $\bm\nu\in\R^d$. To obtain the blending coefficients, we need to minimize
\begin{equation*}
  \eta_e^{\epsilon}=\alpha_e\min_{\bm\nu\in\R^d}\frac{\xbar{\bbu}_P^{e,\star}\cdot\bbn_{*}(\bm\nu)}{\vert \Delta\hbbf_{\bbn_e}(\bbu_P^n,\xbar{\bbu}_P^n,\xbar{\bbu}_{P^-}^n)\cdot\bbn_*(\bm\nu) \vert}
\end{equation*}
and
\begin{equation*}
  \theta_\sigmai^{\epsilon}=\frac{\alpha_\sigmai}{\vert C_\sigmai\vert}\min_{\bm\nu\in\R^d}\frac{ \bbu_\sigmai^{P,\star}\cdot\bbn_{*}(\bm\nu)}{\vert \Delta\bbPhi_\sigmai \cdot\bbn_*(\bm\nu) \vert}.
\end{equation*}
Once the minimization problems is solved, we take
\begin{equation*}
  \eta_e=\min\big(\eta_e^\rho,\eta_e^\epsilon\big) \text{ and } \theta_\sigmai=\min\big(\theta_\sigmai^\rho,\theta_\sigmai^\epsilon\big).
\end{equation*}

Of course, the remaining issue is to evaluate
\begin{equation*}
  \min_{\bm\nu\in\R^d}\frac{\xbar{\bbu}_P^{e,\star}\cdot\bbn_{*}(\bm\nu)}{\vert \Delta\hbbf_{\bbn_e}(\bbu_P^n,\xbar{\bbu}_P^n,\xbar{\bbu}_{P^-}^n)\cdot\bbn_*(\bm\nu) \vert} \text{ and } \min_{\bm\nu\in\R^d}\frac{ \bbu_\sigmai^{P,\star}\cdot\bbn_{*}(\bm\nu)}{\vert \Delta\bbPhi_\sigmai \cdot\bbn_*(\bm\nu) \vert}. 
\end{equation*}
We now discuss this for the flux and for simplicity we denote $\Delta\hbbf_{\bbn_e}:=\Delta\hbbf_{\bbn_e}(\bbu_P^n,\xbar{\bbu}_P^n,\xbar{\bbu}_{P^-}^n)$. It is exactly the same for the residual and thus is omitted here. The first thing to notice is $\xbar{\bbu}_P^{e,\star}\cdot\bbn_{*}(\bm\nu)>0$ so that 
\begin{equation*}
  \min_{\bm\nu\in\R^d}\frac{\xbar{\bbu}_P^{e,\star}\cdot\bbn_{*}(\bm\nu)}{\vert \Delta\hbbf_{\bbn_e}\cdot\bbn_*(\bm\nu) \vert}=\bigg(\max_{\bm\nu\in\R^d}\frac{\vert \Delta\hbbf_{\bbn_e}\cdot\bbn_*(\bm\nu) \vert}{\xbar{\bbu}_P^{e,\star}\cdot\bbn_{*}(\bm\nu)}\bigg)^{-1},
\end{equation*}
where $\xbar{\bbu}_P^{e,\star}=(\xbar{\rho}_P^{e,\star}, \xbar{\rho\bbv}_P^{e,\star},\xbar{E}_P^{e,\star})^\top$ and
\begin{equation*}
  \Delta\hbbf_{\bbn_e}\cdot\bbn_*(\bm\nu)=\Delta\hbbf^{\rho}_{\bbn_e}\frac{\Vert\bm\nu\Vert^2}{2}-\Delta\hbbf^{\rho\bbv}_{\bbn_e}\cdot \bm\nu+\Delta\hbbf^E_{\bbn_e},\quad \xbar{\bbu}_P^{e,\star}\cdot\bbn_{*}(\bm\nu)=\xbar{\rho}_P^{e,\star}\frac{\Vert\bm\nu\Vert^2}{2}-\xbar{\rho\bbv}_P^{e,\star}\cdot \bm\nu+\xbar{E}_P^{e,\star}.
\end{equation*}
We assume that $\bm\nu=\frac{\bbw}{\omega}$ with $(\bbw,\omega)\neq\mathbf{0}$ and obtain
\begin{equation*}
 \frac{ \vert \Delta\hbbf_{\bbn_e}\cdot\bbn_*(\bm\nu) \vert}{\xbar{\bbu}_P^{e,\star}\cdot\bbn_{*}(\bm\nu)}=\frac{\vert\Delta\hbbf^{\rho}_{\bbn_e}\Vert\bbw\Vert^2-2\omega\Delta\hbbf^{\rho\bbv}_{\bbn_e}\cdot \bbw+2\omega^2\Delta\hbbf^E_{\bbn_e}\vert}{\xbar{\rho}_P^{e,\star}\Vert\bbw\Vert^2-2\omega\xbar{\rho\bbv}_P^{e,\star}\cdot \bbw+2\omega^2\xbar{E}_P^{e,\star}}=\frac{\vert\langle \bbz, B\bbz\rangle\vert}{\langle\bbz, C\bbz\rangle},
\end{equation*}
where 
\begin{equation}\label{Rqmatrix}
   \bbz=\left(
          \begin{array}{c}
            \bbw \\
            \omega\\
          \end{array}
        \right)\in\R^{d+1}\backslash\{\mathbf 0\},~
   B=\left(
           \begin{array}{cc}
             \Delta\hbbf^{\rho}_{\bbn_e}\Id_d & -\Delta\hbbf^{\rho\bbv}_{\bbn_e}  \\
             -\big(\Delta\hbbf^{\rho\bbv}_{\bbn_e}\big)^\top & 2\Delta\hbbf^E_{\bbn_e}  \\
           \end{array}
         \right),~\text{and}~C=\left(
           \begin{array}{cc}
             \xbar{\rho}_P^{e,\star}\Id_d & -\xbar{\rho\bbv}_P^{e,\star}  \\
             -\big(\xbar{\rho\bbv}_P^{e,\star}\big)^\top & 2\xbar{E}_P^{e,\star}  \\
           \end{array}
         \right).
   \end{equation}
Notice that $\langle\bbz, C\bbz\rangle=\xbar{\bbu}_P^{e,\star}\cdot\bbn_{*}(\bm\nu)>0$, we know that $C$ is positive definite and the problem reduces to study the Rayleigh quotient
\begin{equation*}
  \max\limits_{\substack{\bbz\in \R^{d+1}\\\bbz\neq\mathbf 0}}\frac{\vert\langle \bbz, B\bbz\rangle\vert}{\langle\bbz, C\bbz\rangle}=\max\limits_{\substack{\bbz\in \R^{d+1}\\\bbz\neq\mathbf 0}}\frac{\vert\langle \bbz, C^{-1/2}BC^{-1/2}\bbz\rangle\vert}{\Vert\bbz\Vert^2}=\max\limits_{\lambda \text{ eigenvalue of }C^{-1/2}BC^{-1/2}}\vert\lambda\vert .
\end{equation*}
The evaluation of eigenvalues of $C^{-1/2}BC^{-1/2}$ can be done by some iterative method. It turns out that in the particular case we consider, we can find a simple analytical formula.
Let us temporarily denote the matrices in \eqref{Rqmatrix} as
\begin{equation*}
  B=\begin{pmatrix}
\beta_0\Id_d & -\bbb\\
-\bbb^\top & 2\beta_{d+1}\end{pmatrix}, \quad C=\begin{pmatrix}
\alpha_0\Id_d & -\bba\\
-\bba^\top & 2\alpha_{d+1}\end{pmatrix}.
\end{equation*}
If $\lambda$ is an eigenvector of $C^{-1/2}BC^{-1/2}$, this means there is a non zero $\bbz=(\bbx, \theta)$, $\bbx\in \R^d$, $\theta\in \R$ such that $C^{-1/2}BC^{-1/2}\bbz=\lambda\bbz$, i.e
\begin{equation*}
  \big ( B-\lambda C\big ) \bbz=0
\end{equation*}
and $\lambda$ must be real since $C^{-1/2}BC^{-1/2}$ is symmetric in $\R$. This means that
\begin{equation}\label{Rq1}
\left\{\begin{array}{ll}
\big (\beta_0-\lambda \alpha_0\big )\bbx&-\big (\bbb-\lambda\bba\big )\theta=0,\\
-\big(\bbb-\lambda\bba\big)^T\bbx&+2\big (\beta_{d+1}-\lambda \alpha_{d+1}\big )\theta=0.
\end{array}\right.
\end{equation}
If $\theta=0$, we must have $\lambda=\tfrac{\beta_0}{\alpha_0}$ (remember that $\alpha_0=\xbar{\rho}_P^{e,\star}>0$) from the first equation in \eqref{Rq1}. If $\theta\neq 0$,
by multiplying the first line by $\big(\bbb-\lambda\bba\big)^T$ and the second one by  $\beta_0-\lambda \alpha_0$, we end up with the condition 
\begin{equation*}
  \theta \Big ( -\big(\bbb-\lambda\bba\big)^T\big(\bbb-\lambda\bba\big)+2\big (\beta_{d+1}-\lambda \alpha_{d+1}\big )\big (\beta_0-\lambda \alpha_0\big )\Big )=0,
\end{equation*}
which is equivalent to
\begin{equation*}
  \big ( \Vert \bba\Vert^2-2\alpha_0\alpha_{d+1}\big )\lambda^2+2\big(\beta_{d+1}\alpha_0+\beta_0\alpha_{d+1}-\bba^\top\bbb\big )\lambda+\big (\Vert \bbb\Vert^2-2\beta_0\beta_{d+1}\big )=0.
\end{equation*}
Since $C$ is positive definite, we must have $\det C=\alpha_0\big (2 \alpha_0\alpha_{d+1}-\Vert \bba\Vert^2\big )>0$ so that $\Vert \bba\Vert^2- 2\alpha_0\alpha_{d+1}<0$. From this we get that the two other eigenvalues are
\begin{equation*}
  \lambda_{\pm}=\dfrac{\big (\bba^T\bbb-\beta_{d+1}\alpha_0-\beta_0\alpha_{d+1}\big )\pm \sqrt{\Delta}}{\Vert \bba\Vert^2-2\alpha_0\alpha_{d+1}}
\end{equation*}
with
\begin{equation*}
  \Delta=\big (\beta_{d+1}\alpha_0+\beta_0\alpha_{d+1}-\bba^T\bbb\big )^2-\big (\Vert \bba\Vert^2-2\alpha_0\alpha_{d+1}\big )\big (\Vert \bbb\Vert^2-2\beta_0\beta_{d+1}\big )\geq 0.
\end{equation*}
These formulae are independent of the dimension, and
\begin{equation*}
  \max\limits_{\lambda \text{ eigenvalue of }C^{-1/2}BC^{-1/2}}\vert\lambda\vert =\max\Big(\frac{\vert\beta_0\vert}{\alpha_0}, \vert \lambda_+\vert, \vert\lambda_-\vert\Big).
\end{equation*} 
Next, we use the values in \eqref{Rqmatrix} to get
\begin{equation*}
  \max\limits_{\lambda \text{ eigenvalue of }C^{-1/2}BC^{-1/2}}\vert\lambda\vert =\max\Big(\frac{\vert\Delta\hbbf^{\rho}_{\bbn_e}\vert}{\xbar{\rho}_P^{e,\star}},\frac{\vert\kappa_1\vert+\sqrt{\Delta}}{-\kappa_0}\Big),
\end{equation*}
where $\kappa_0=\Vert\xbar{\rho\bbv}_P^{e,\star}\Vert^2-2\xbar{\rho}_P^{e,\star}\xbar{E}_P^{e,\star}$, $\kappa_1=\xbar{\rho\bbv}_P^{e,\star}\cdot \Delta\hbbf^{\rho\bbv}_{\bbn_e} -\xbar{\rho}_P^{e,\star}\Delta\hbbf^{E}_{\bbn_e}-\xbar{E}_P^{e,\star}\Delta\hbbf^{\rho}_{\bbn_e}$, and
\begin{equation*}
  \Delta=\kappa_1^2-\kappa_0\big(\Vert\Delta\hbbf^{\rho\bbv}_{\bbn_e}\Vert^2-2\Delta\hbbf^{\rho}_{\bbn_e}\Delta\hbbf^{E}_{\bbn_e}\big).
\end{equation*}
Finally, we take
\begin{equation*}
  \eta_e^{\epsilon}= \alpha_e\bigg(\max\limits_{\lambda \text{ eigenvalue of }C^{-1/2}BC^{-1/2}}\vert\lambda\vert\bigg)^{-1}= \alpha_e\bigg(\max\Big(\frac{\vert\Delta\hbbf^{\rho}_{\bbn_e}\vert}{\xbar{\rho}_P^{e,\star}},\frac{\vert\kappa_1\vert+\sqrt{\Delta}}{-\kappa_0}\Big)\bigg)^{-1}.
\end{equation*}

\begin{remark} This technique is not specific to the schemes developed in this paper but has a larger potential. It has also been used in \cite{GauthierLBM}, for a completely different scheme that uses a kinetic formulation of the Euler equations.
\end{remark}

\subsection{Boundary conditions}\label{sec:BCs}
This section only deals with the Euler equations: all the scalar problems are chosen such that the solution does not change in a neighborhood of the boundary.
We need to consider 3 types of boundary conditions: i) wall, ii) inflow/outflow, iii) and {homogeneous} Neumann type {(i.e., the zero-order extrapolation boundary conditions)}. We describe what is done first for the average update and then the point values.
The average values are evolved by
\begin{equation}\label{scheme:flux:bc}
    \xbar{\bbu}_P^{n+1}=\xbar{\bbu}_P^n-\frac{\dt}{\vert P\vert}\bigg ( \sum_{e\in \FF_P }\vert e\vert \; \hbbf_{\bbn_e}(\bbu_P^n,\xbar{\bbu}_P^n,\xbar{\bbu}_{P^-}^n)
    +\sum_{e^b\in \FF_P^b}\vert e^b\vert \; \hbbf^b_{\bbn_e^b}(\bbu_P^n, \xbar{\bbu}_P^n)
    \bigg ),
\end{equation}
and the point values by:
\begin{equation}\label{scheme:residual:bc}
\bbu_\sigmai^{n+1}=\bbu_\sigmai^n-\dt\sum_{\substack{P, \sigmai \in P\\ \sigmai\notin\Gamma
}} \bbPhi_\sigmai^P(\bbu)-\dt\sum_{\substack{P, \sigmai \in P\\ \sigmai\in\Gamma}}\bbPhi_{\sigmai}^{P,b}(\bbu)
\end{equation}
where the boundary flux $\hbbf^b$ and the boundary residuals $\bbPhi_{\sigmai}^{P,b}$ need to be defined.
Since the structure of the updates are similar to \eqref{ave:scheme} and \eqref{scheme:residual:o}, we can apply the same technique and define blending $\theta_b$ and $\eta_b$.
\subsubsection{Average values}
The problem is to evaluate the flux contribution on the faces $e^b$ of a polygon $P$ which is on the boundary of the computational domain. 
\begin{figure}[!h]
\begin{center}
{\includegraphics[trim=0.25cm 0.001cm 0.25cm 0.001cm,clip,width=8.5cm]{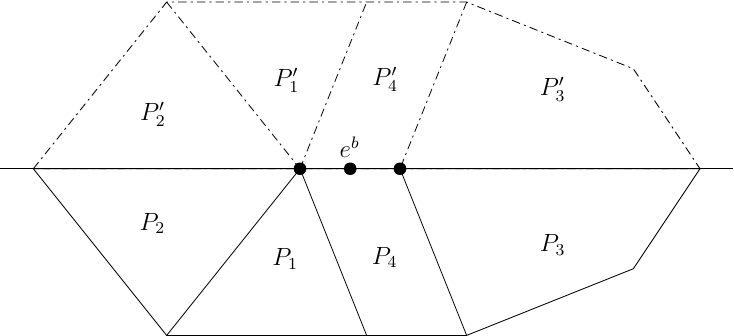}}
\caption{\label{fig:bc} Notation for the boundary conditions.}
\end{center}
\end{figure}
\begin{itemize}
\item Wall. The flux at the wall is obtained from a numerical flux $\hbbf_{\bbn_{e^b}}$ given by
\begin{equation*}
{\int_{e^b} \hbbf_{\bbn_{e^b}}(\bbu_{e^b},\bbu_{e^b}^s)\; {\rm d}\gamma},
\end{equation*}
where $\int$ is a quadrature formula \footnote{We always use the Gauss--Lobatto points since they are our degrees of freedom on the edges.}, $\bbu_{e^b}=(\rho, \rho \bbv, E)^T$ represents the polynomial approximation on $e^b$, and $\bbu_{e^b}^s$ is the state with the same density, the same total energy and the velocity which has been symmetrized with respect to $\bbn_{e^b}$ {using the following symmetric operator:
\begin{equation*}
  s_{\bbn_{e^b}}(\bbv)=\bbv-2\frac{\langle \bbv,\bbn_{e^b}\rangle}{\Vert \bbn_{e^b}\Vert^2}\bbn_{e^b}.
\end{equation*}} 
In practice, {for the high-order case, we take
\begin{equation*}
  \hbbf_{\bbn_{e^b}}=\frac{\bbf(\bbu_{e^b})+\bbf(\bbu_{e^b}^s)}{2}\cdot \bbn_{e^b},
\end{equation*}
while for the low order scheme,} $\lbbf_{\bbn_{e^b}}$ is the Rusanov flux given by \eqref{flux:LO} where the state $\xbar{\bbu}_{P^-}$ is simply replaced with the symmetric state of $\xbar\bbu_{P}$ denoted by $\xbar\bbu_P^s$. 

\item Inflow/outflow. The flux across the boundary edge $e^b$ is evaluated as 
\begin{equation*}
{\int_{e^b} \hbbf_{\bbn_{e^b}}(\bbu_{e^b},\bbu^\infty)\; {\rm d}\gamma},
\end{equation*}
where $\bbu^\infty$ is the state at infinity. To make sure one \red{gets} complete upwinding for supersonic inflow or outflow, the {high-order} numerical flux is evaluated by modified Steger--Warming flux as:   
{
\begin{equation*}
  \hbbf_{\bbn_{e^b}}(\bbu_{e^b},\bbu^\infty)=\left\{
  \begin{array}{ll}
  \big(\bbA(\bbu_{e^b})\cdot\bbn_{e^b}\big)^+\bbu_{e^b}+\big(\bbA(\bbu^\infty)\cdot\bbn_{e^b}\big)^-\bbu^\infty,& \text{ at an inflow boundary},\\
  \big(\bbA(\bbu_{e^b})\cdot\bbn_{e^b}\big)^-\bbu_{e^b}+\big(\bbA(\bbu^\infty)\cdot\bbn_{e^b}\big)^+\bbu^\infty,& \text{ at an outflow boundary},
  \end{array}\right.
\end{equation*}
while the low order numerical flux is the Rusanov one given by \eqref{flux:LO} where the state $\bbu_{P^-}$ is simply $\bbu^\infty$.}

\item {homogeneous} Neumann type. We refer to Figure \ref{fig:bc}. For the edge $e^b$ of $P_4$, we consider the element $P_4'$ obtained by symmetrizing $P_4$ with respect to the outward normal to $e^b$ and we put in $P'_4$ the same state as in $P_4$, then the high- and low-order fluxes in \eqref{flux} and \eqref{flux:LO} are computed. 
\end{itemize}

{Applying the above defined boundary values and high- and low-order fluxes to the numerical flux blending procedure, we can define the corresponding blending coefficient $\eta_{e^b}$ using \eqref{eq:theta_average}.}

\subsubsection{Point values}
The update of the $i$-th point values localized at the $\bullet$ points depicted in Figure \ref{fig:bc} is obtained as 
\begin{equation*}
  \dfrac{{\rm d}\bbu_\sigmai}{{\rm d}t}+\sum_{P, \sigmai\in P}\bbPhi_\sigmai^{P}(\bbu)=0,
\end{equation*}
where in the sum we consider the polygons of the triangulation that share $\sigmai$ and we have also included those obtained by symmetrization: the polygons {$P'_1$, $P'_2$, $P'_3$, and $P'_4$} of Figure \ref{fig:bc}. Here we focus on how we evaluate $\bbPhi_\sigmai^{P}$ for $P$.
\begin{itemize}
\item Wall. In the {polygon $P_i$, $i=1,\ldots, 4$} with state $(\rho, \rho \bbv, E)^T$, we consider the symmetrized state {in the corresponding element $P'_i$} given by {$(\rho, \rho s_{\bbn_{e^b}}(\bbv), E)$}. All the geometrical elements are also mirrored with respect to $e^b$.
\item Inflow/outflow. {We populate the element $P'_i$ with the state $(\rho_\infty, \rho_\infty\bbv, E_\infty)^T$.} 
\item Neumann. {We feed the elements $P'_i$ with the state of $P_i$.}
\end{itemize}
For the implementation, we have to take into account that, for the first order scheme, the area $C_\sigmai$ is again evaluated as in \eqref{area}. {With the obtained symmetrized values, we can define the residual blending coefficients $\eta_\sigmai^b$ using the same procedure for \eqref{eq:theta_point}.}

\section{Results}
{In this section, we demonstrate the robustness and effectiveness of the proposed scheme on a number of classical numerical examples. These results also show that the {monolithic convex limiting} procedure is more accurate than the {a-posteriori limiting} one.}
\subsection{Scalar case}
\subsubsection{Convection case}
\begin{equation}\label{scalar:test}\dpar{u}{t}+\mathbf{a}\cdot \nabla u=0, \qquad u_0(x)=\exp(-20||\bbx-\bbx_0||^2).
\end{equation}
We run this case on the domain $[-2,2]^2$, the velocity field is $\bba=2\pi(-y,x)$ and $\bbx_0=(0,1)$ for one full rotation, i.e. $t_f=1$. The errors for a triangular mesh are given in Table \ref{error:convection:tri}.  We start from a coarse mesh, and {the finer meshes are generated by cutting the triangles into 4 sub-triangles by connecting the edge midpoints with three new edges.} From a given mesh, we construct a polygonal mesh by agglomeration. \red{We have used the high order scheme, without BP stabilisation.}
The errors as well as the convergence order of the new schemes are given in Table \ref{error:convection:poly}, confirming that the formal third-order accuracy is achieved in all cases.
\input{pente}

\subsubsection{KPP case}
This problem has been considered by Kurganov, Popov and Petrova in \cite{zbMATH05380178}. It writes 
   \begin{equation}
\label{eq:KPP}
\dpar{u}{t}+\dpar{\big(\sin \; u\big )}{x}+\dpar{\big (\cos \; u\big )}{y}=0,
\end{equation}
in a domain $[-2,2]^2$
with the initial condition
$$u_0(\bbx)=\left \{\begin{array}{ll}
\frac{7}{2}\pi& \text{ if } \Vert \bbx-(0,0.5)\Vert^2 \leq 1,\\
\frac{\pi}{4}& \text{ else.}
\end{array}
\right .
$$
The problem \eqref{eq:KPP} is non-convex, in the sense that compound waves may exist {, characterized by a shock and a rarefaction wave that are attached together}. Here, we have used the {two limiting strategies---a-posteriori limiting and a monolithic convex limiting---}we have described above. We compute the numerical solution until the final time $t_f=1$. The triangular mesh has $29909$ point value DoFs,  $14794$ {triangles} and $22351$ faces. The polygonal mesh has $38425$ point value DoFs,  $7558$ elements and $22991$ faces. Though the number of point value DoFs is larger for the polygonal mesh, the resolution of the polygonal mesh is similar to that of the triangular one because the total number of DoF is the same for both the meshes which are both very regular. 

\begin{figure}[!ht]
\begin{center}
\subfigure[{Convex limiting, Polygonal mesh}]{\includegraphics[width=0.45\textwidth]{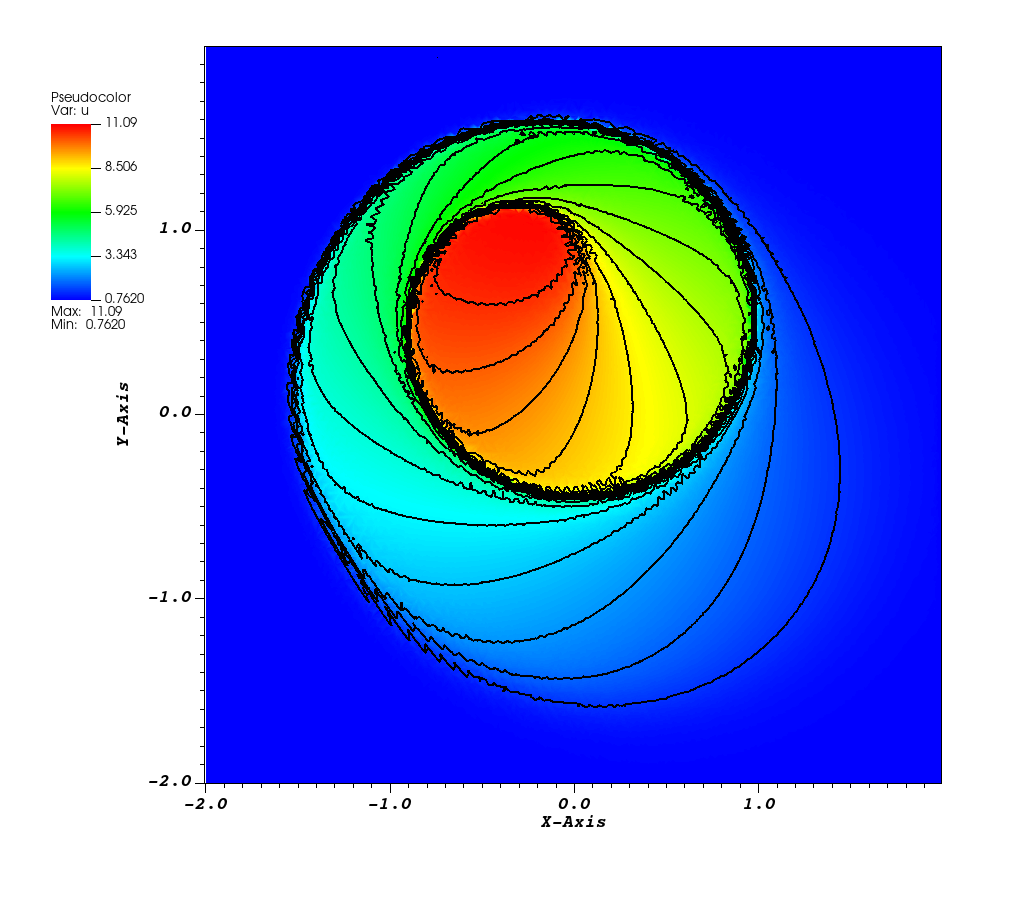}}
\subfigure[{Convex limiting, Triangular mesh}]{\includegraphics[width=0.45\textwidth]{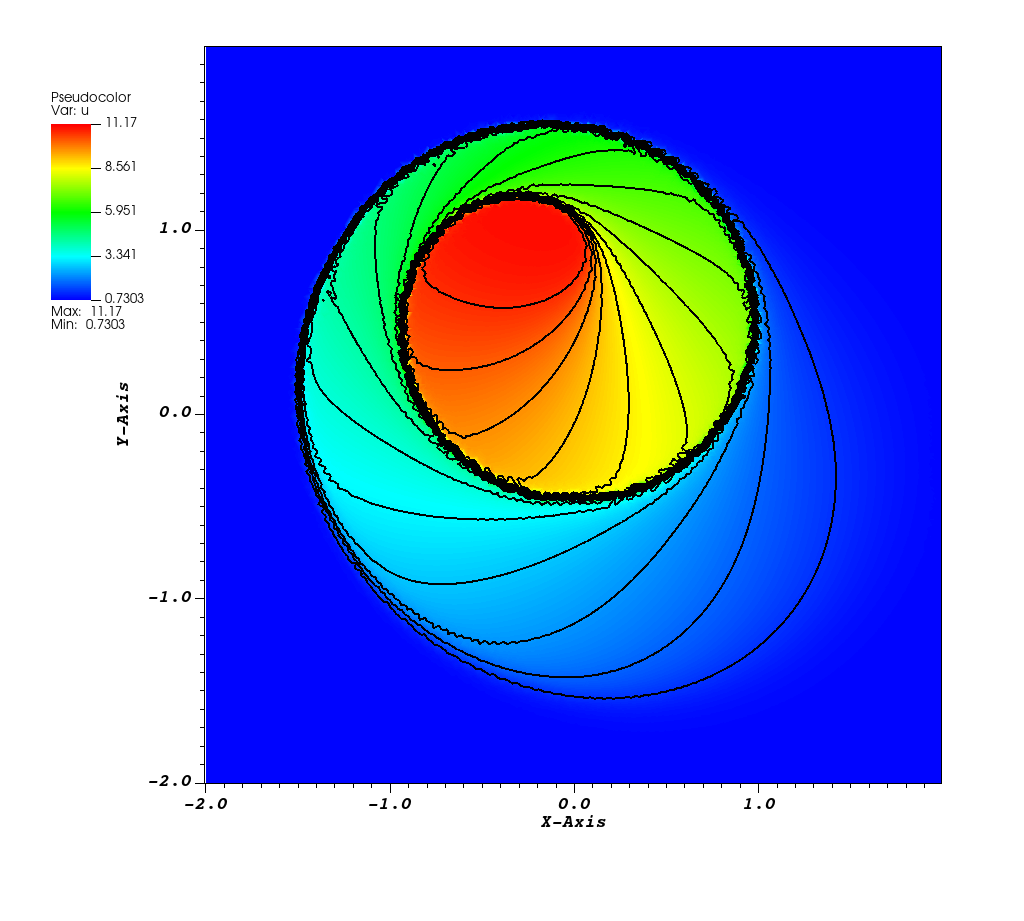}}
\subfigure[{a-posteriori limiting, Polygonal mesh}]{\includegraphics[width=0.45\textwidth]{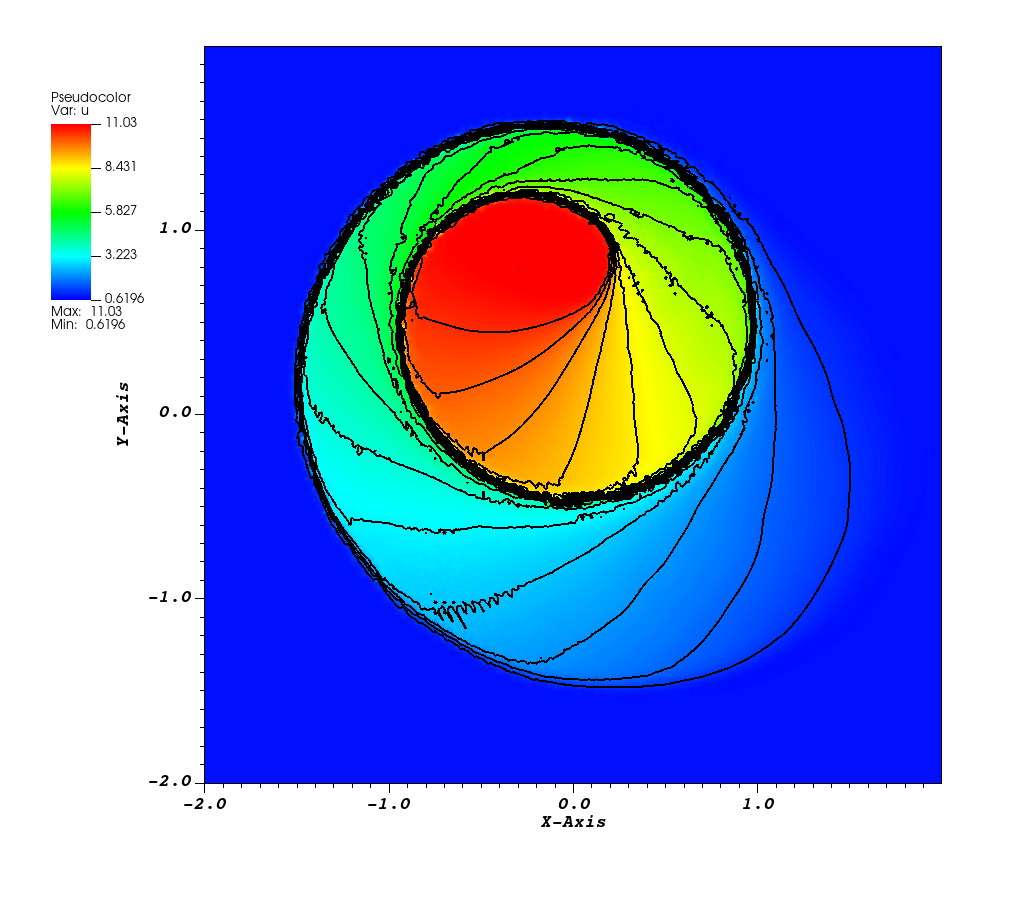}}
\subfigure[{a-posteriori limiting, Triangular mesh}]{\includegraphics[width=0.45\textwidth]{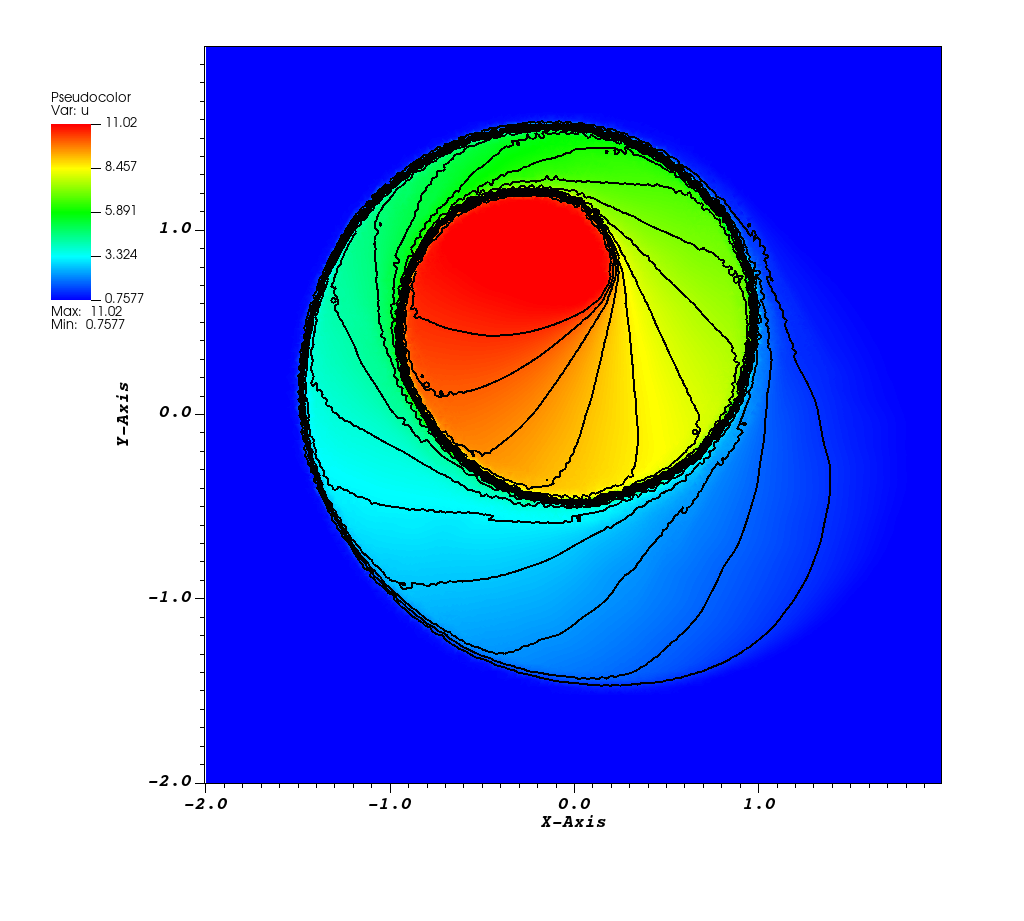}}
\end{center}
\caption{\label{fig:kpp_p2} KPP problem using quadratic approximation. Average values. The CFL is $0.3$.}
\end{figure}

The results we obtain are in good agreement with published results. In particular, we have managed to have a good shock structure. In that respect, the use of Rusanov scheme for the point values and local Lax--Friedrichs flux for the average values, as first order scheme, is important to get the correct shock structure. Comparing the rows (a)-(b) and (c)-(d) of Figure \ref{fig:kpp_p2} indicates that {the scheme with convex limiting strategy} provides smoother solution than those obtained with the {a-posteriori limiting approach}, with crisp shock structure.

\subsection{System case}
\subsubsection{Acoustics}
For $\bbv\in \R^2$ and $p\in \R$, the acoustics system reads 
\begin{equation}\label{eq:acoustics}
\left\{\begin{array}{ll}
\dpar{\bbv}{t}&+\nabla p=0,\\
\dpar{p}{t}&+c^2\text{ div }\bbv=0.
\end{array} \right.
\end{equation}
In \cite{Barsukow}, it was mentioned that Active Flux on Cartesian meshes manages to preserve very well the steady solutions of \eqref{eq:acoustics}.
For that, besides a theoretical analysis, examples of the following problem is shown:
\begin{equation}\label{acoustic:steady}p(\bbx,t=0)=0, \quad \bbv(\bbx,t=0)=\left\{ \begin{array}{ll}
0 &\text{ if }\Vert \bbx\Vert\geq 0.4,\\
(2-5\Vert \bbx\Vert ){\frac{\bbx^\bot}{\Vert \bbx\Vert}} & \text{ if } 0.2\leq \Vert \bbx\Vert\leq 0.4,\\
5{\frac{\bbx^\bot}{\Vert\bbx\Vert}} & \text{ if } 0\leq \Vert \bbx\Vert \leq 0.2,
\end{array}
\right.
\end{equation}
where for $\bbx=(a,b)$, $\bbx^\bot=(-b,a)$. In \eqref{eq:acoustics}, we set $c=1$, and run the scheme until $t_f=100$ in the domain $[-1,1]^2$. The problem \eqref{acoustic:steady} is a steady problem, so the solution should not change. Note that the mesh has no particular symmetry and it is not particularly regular (it has been obtained with the option Delaunay in meshing with GMSH).
\begin{figure}[!h]
\begin{center}
\subfigure[Triangular mesh, average values]{\includegraphics[width=0.45\textwidth]{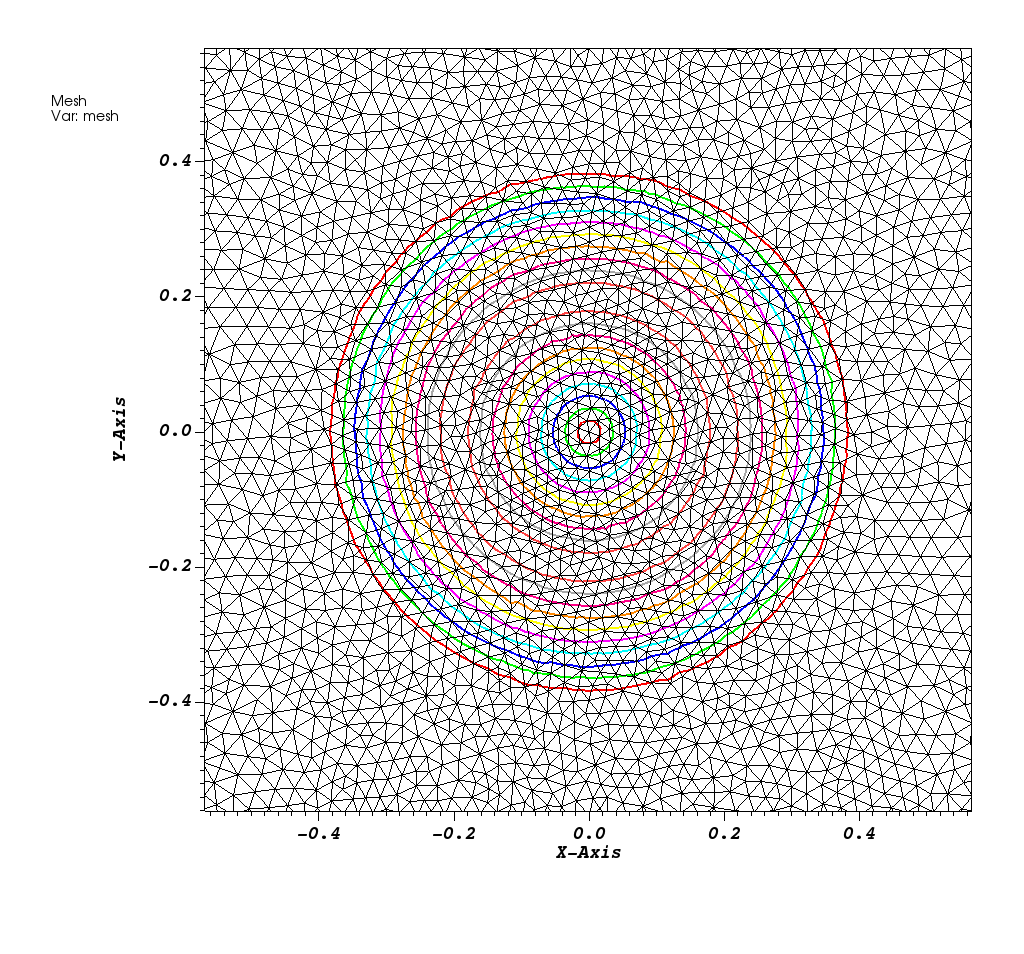}}
\subfigure[Triangular mesh, point values]{\includegraphics[width=0.45\textwidth]{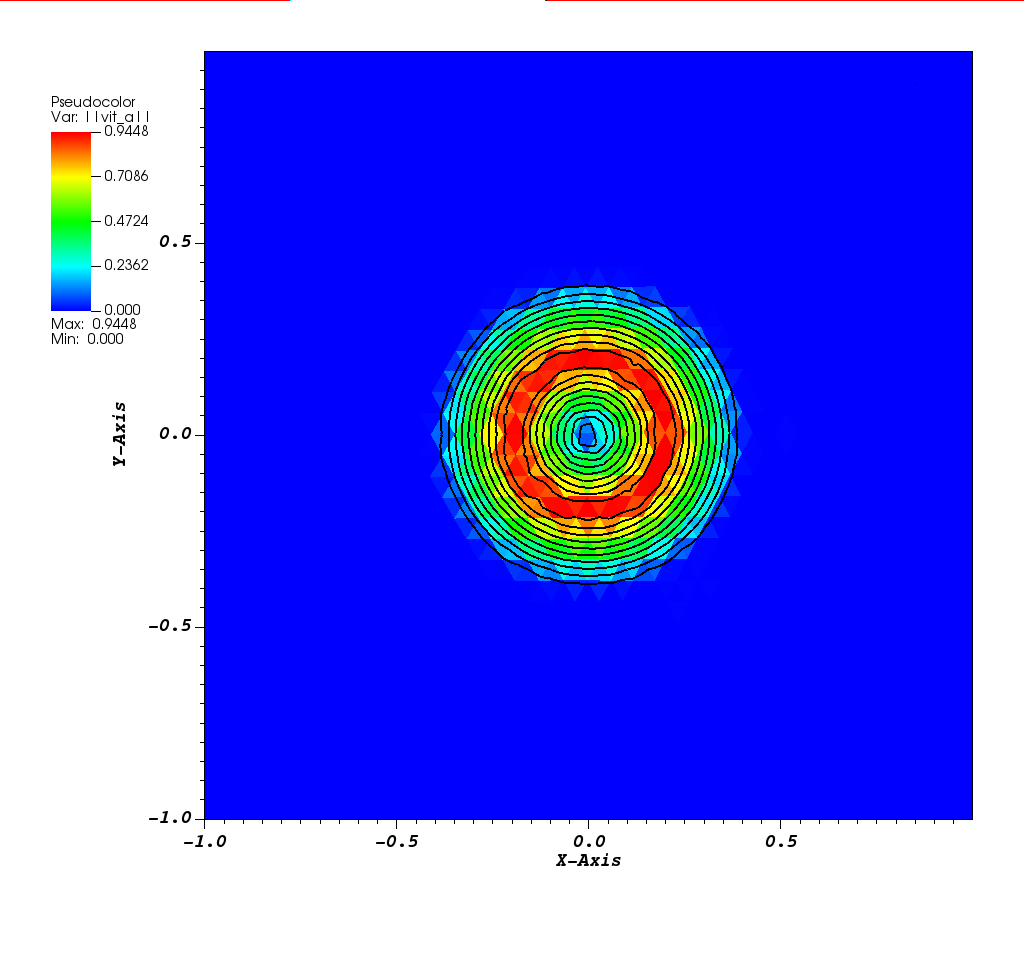}}
\subfigure[Polygonal mesh, average values]{\includegraphics[width=0.45\textwidth]{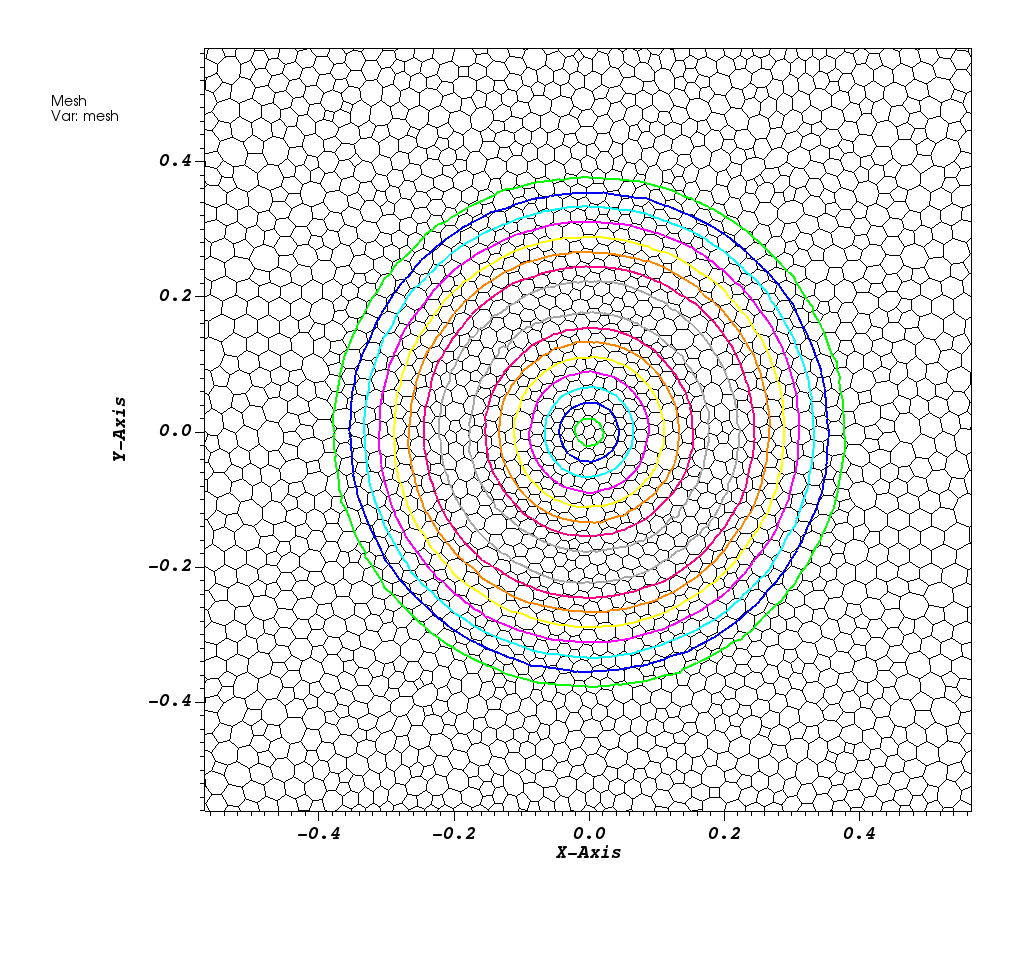}}
\subfigure[Polygonal mesh, point values]{\includegraphics[width=0.45\textwidth]{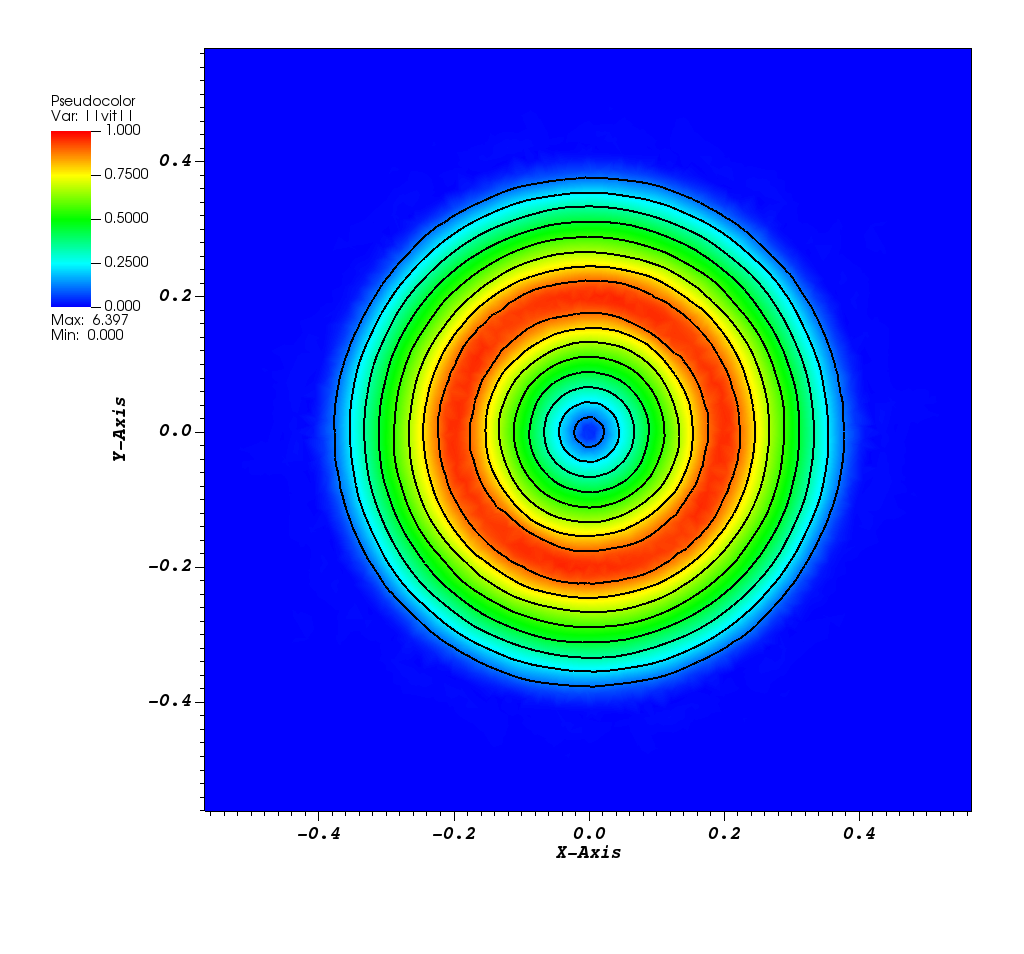}}
\end{center}
\caption{\label{sol:acoustics} Acoustics: plots of the velocity and the computational mesh.}
\end{figure}
\begin{figure}[!h]
\begin{center}
\subfigure[Polygonal mesh]{\includegraphics[width=0.45\textwidth]{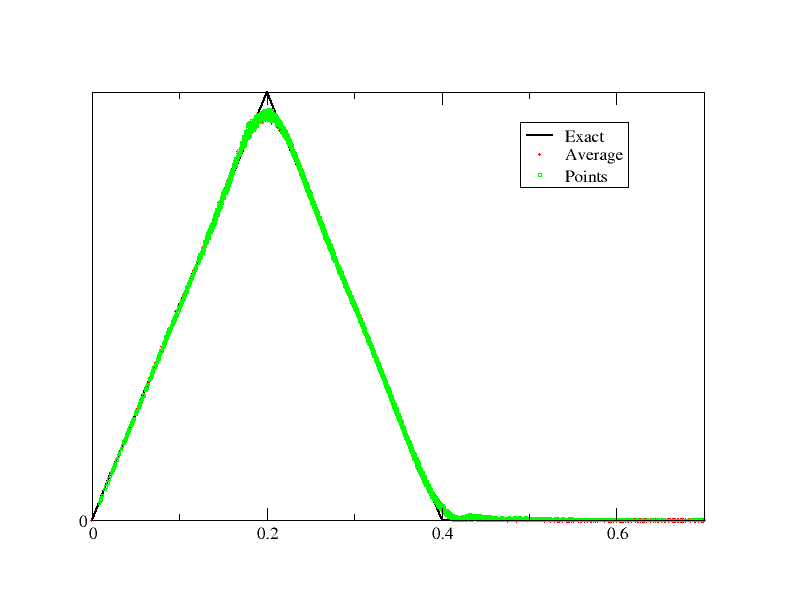}}
\subfigure[Triangular mesh]{\includegraphics[width=0.45\textwidth]{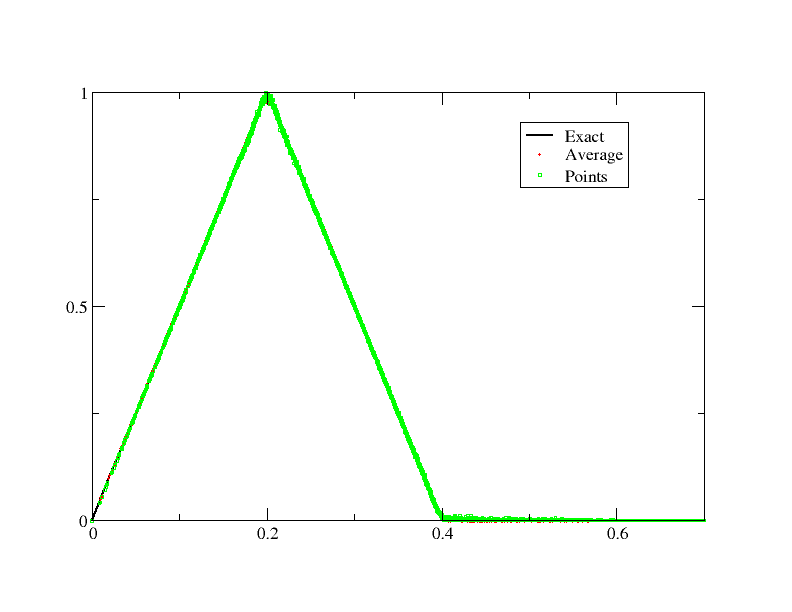}}
\end{center}
\caption{\label{sol:acoustics:sca} Acoustics: scatter plots of the velocity norm, compared to the exact solution.}
\end{figure}
Figure \ref{sol:acoustics} shows the norm of the velocity field at time $t_f=100$, as well as the meshes that have been used. The scatter plots depicted in Figure \ref{sol:acoustics:sca} confirm that the solution has very little dispersion and is almost equal to the initial condition. {The CFL number is set to $0.4$, to reach $T=100$, we need $15,000$ time steps for the triangular mesh and about $9,000$ for the polygonal one which is coarser.} Notice that the meshes are very coarse. This seems to indicate that this kind of schemes also have a very good behavior with respect to irrotational flows, as explained in \cite{Barsukow}. However, the explanation is likely to be different because the schemes, and the meshes are very different.
\subsubsection{Euler equations}
Four cases are considered:  the moving vortex case, two examples of 2D Riemann problems, and the double Mach reflection case.
\paragraph{Moving vortex case.}\label{vortex}
In the first example of nonlinear Euler equations, a moving vortex case is defined in the computational domain $[-20,20]^2$. The final time of the simulation is $t_f=20$, and the initial condition is given by
  {
\begin{equation*}
  \rho=\big ( T_\infty+\Delta T\big )^{\frac{1}{\gamma-1}},\quad \bbv=\bbv_\infty+Me^{\frac{1-R}{2}}\widehat{\bbx},\quad p=\rho^\gamma,
\end{equation*}
where $T_\infty=1$, $M=\frac{5\pi}{2}$, $\bbv_\infty=(1,\frac{\sqrt 2}{2})$, $\widehat{\bbx}=(\bbx^\bot-\bbx_0)/2$ with $\bbx_0=(-10,-10)$, and
\begin{equation*}
  R=\Vert\widehat{\bbx}\Vert^2, \quad \Delta T=-\frac{\gamma-1}{2\gamma}M^2e^{1-R}.
\end{equation*}}

\begin{figure}[h!]
\begin{center}
  \subfigure[Plot of the density with 20 equally spaced isolines]{ \includegraphics[width=0.45\textwidth]{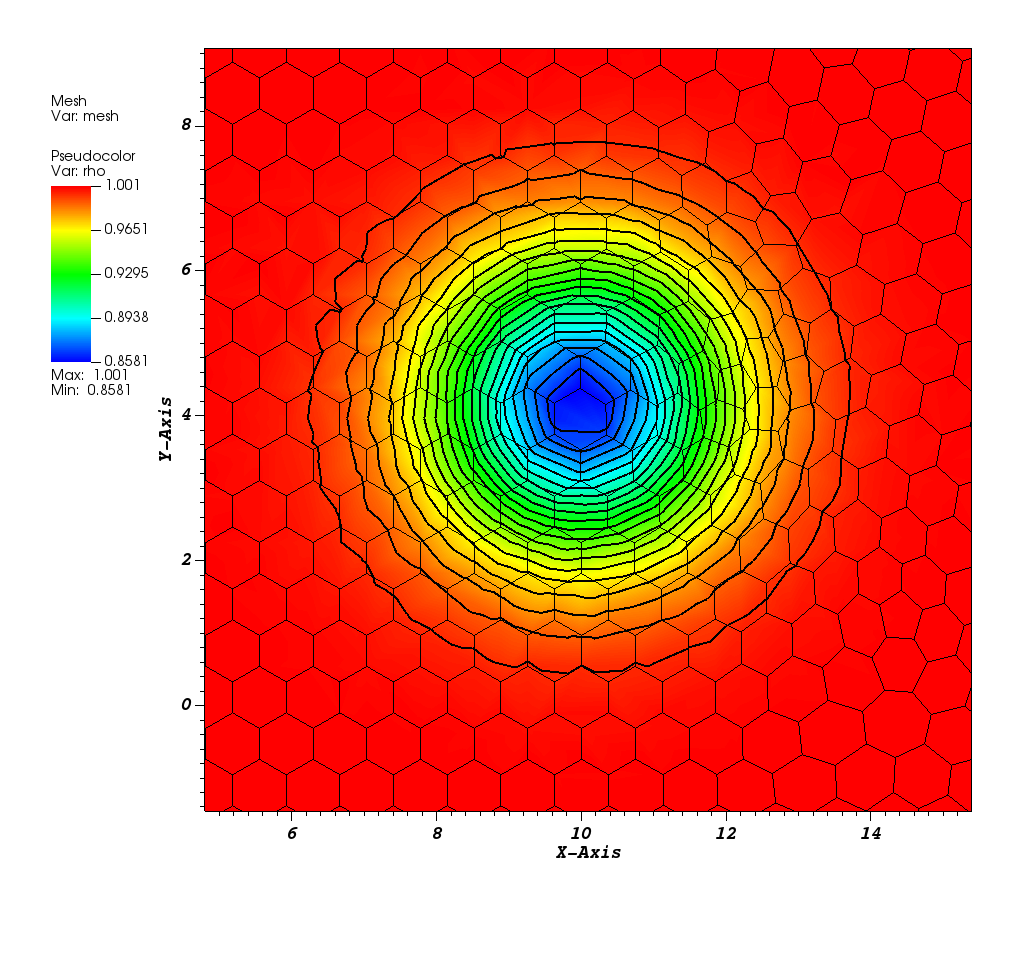}}
    \subfigure[Exact solution (red) versus the numerical one]{ \includegraphics[width=0.45\textwidth]{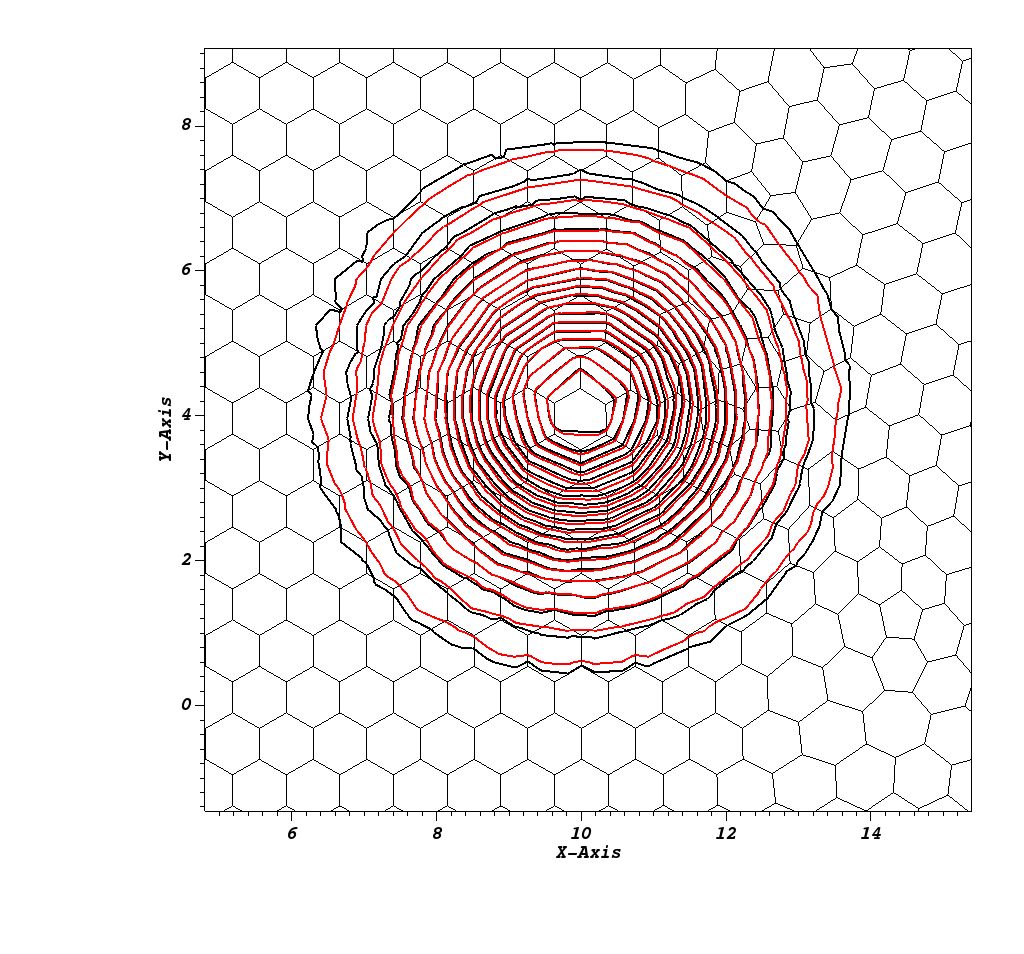}}
\end{center}
\caption{\label{fig:vortex} {Moving vortex. The CFL is $0.2$ and only the high-order schemes introduced in section \ref{sec31} will be activated.}}
\end{figure}

The solution computed on a polygonal mesh with 18252 point DoFs and 3565 elements is shown in Figure \ref{fig:vortex}-(a). Additionally, a comparison with the exact solution is depicted in Figure \ref{fig:vortex}-(b), where the exact solution is obtained by advecting the vortex with the velocity $\bbv_\infty$.

\begin{figure}[!h]
\begin{center}
\subfigure[{a-posteriori limiting, point values}]{\includegraphics[width=0.45\textwidth]{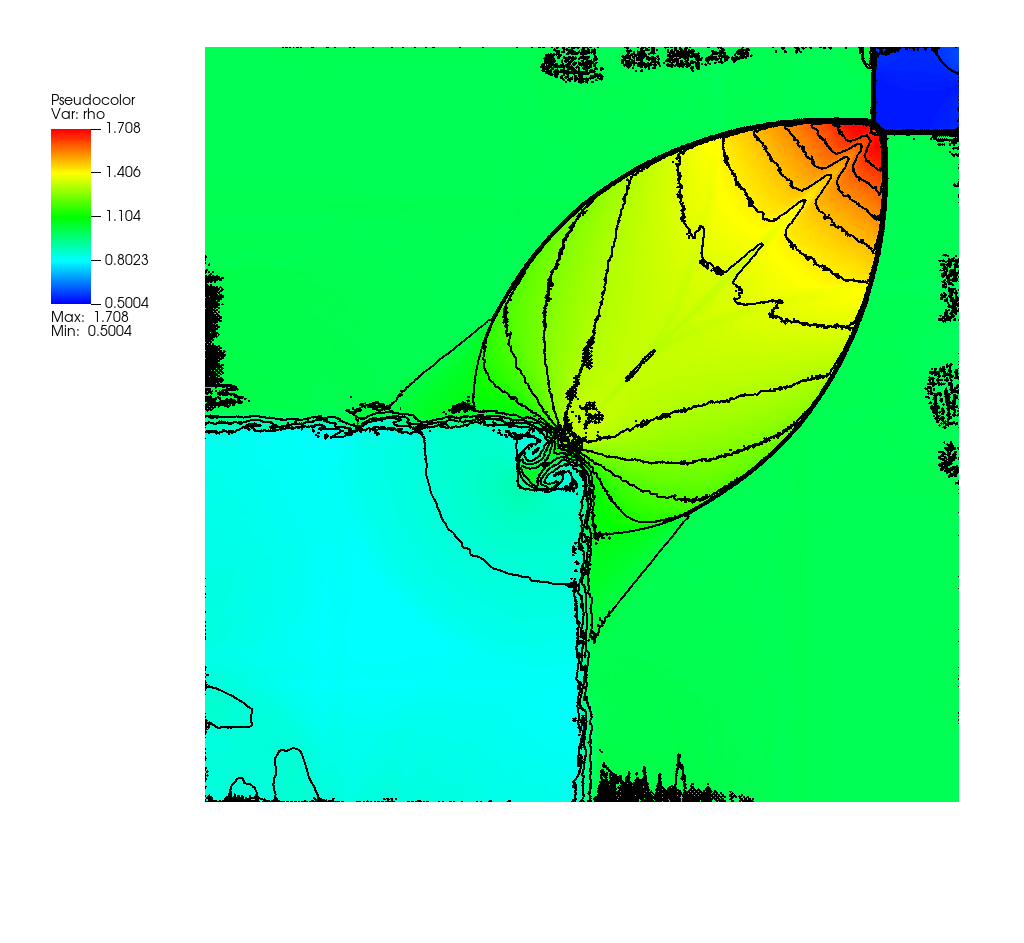}}
\subfigure[{convex limiting, point values}]{\includegraphics[width=0.45\textwidth]{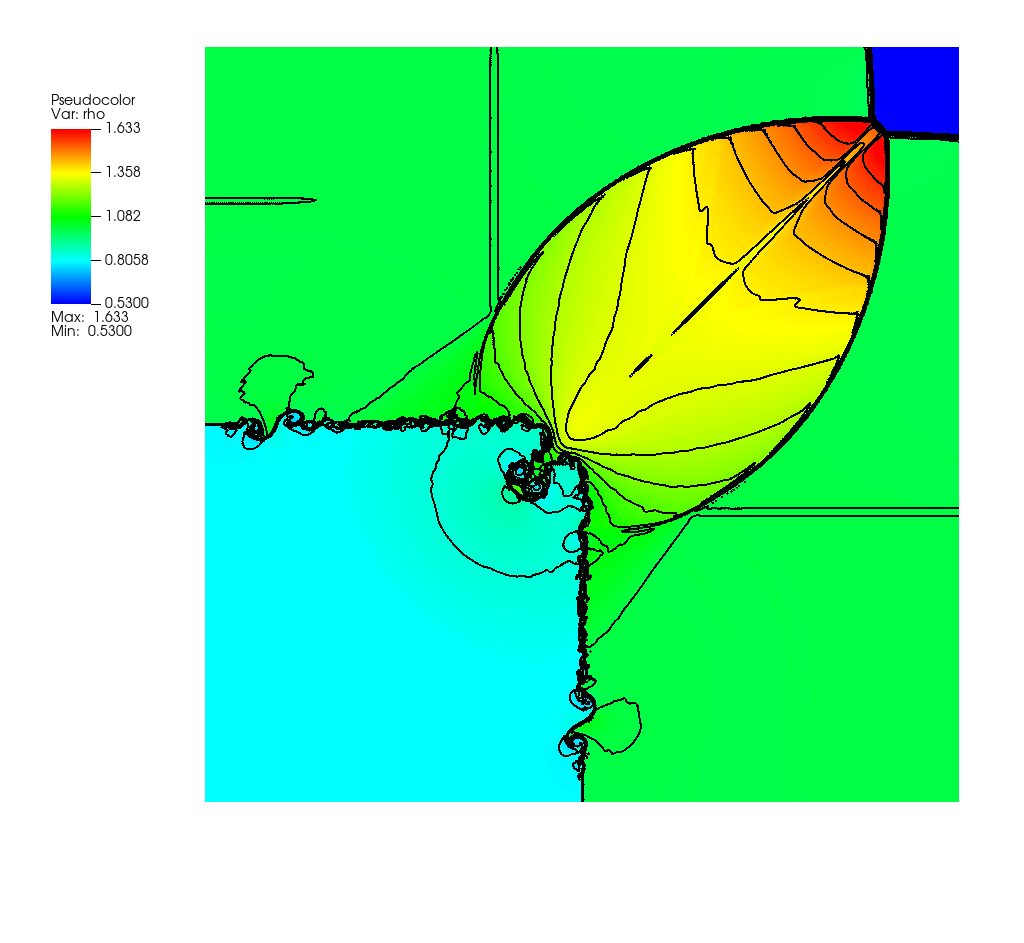}}
\subfigure[{a-posteriori limiting, average values}]{\includegraphics[width=0.45\textwidth]{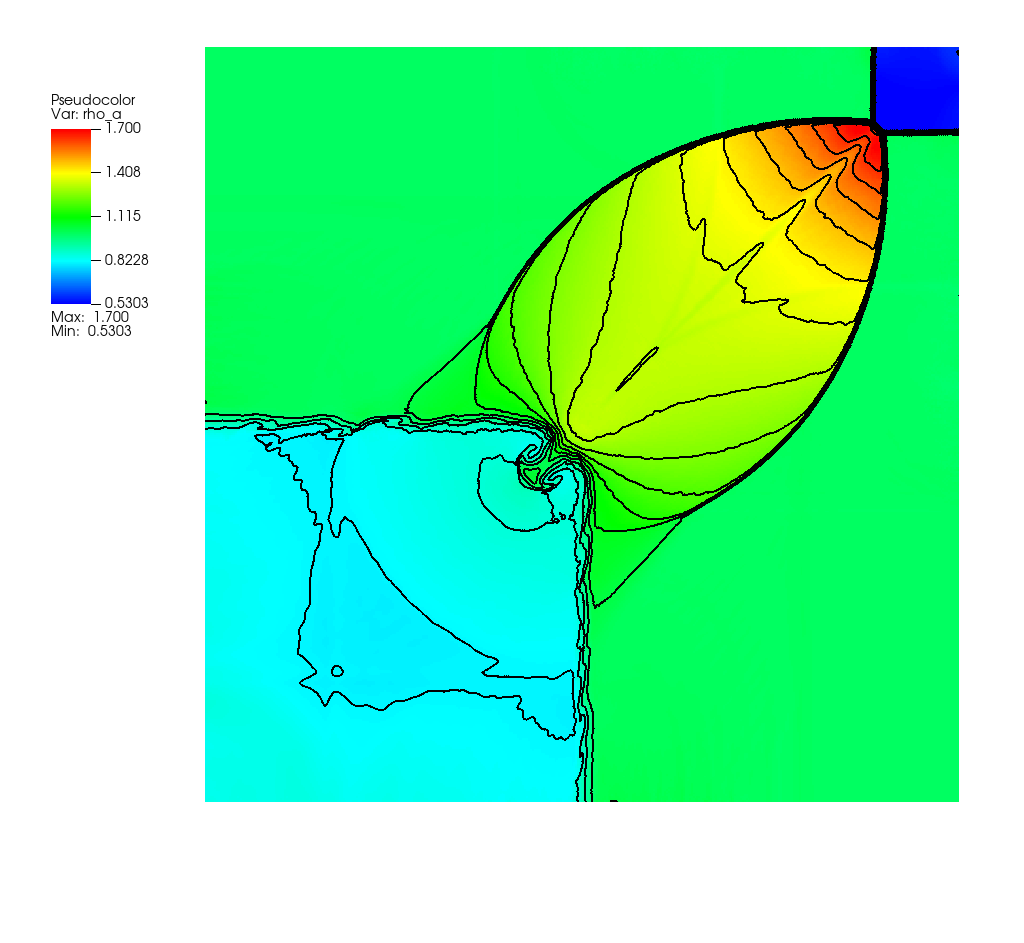}}
\subfigure[{convex limiting, \red{average} values}]{\includegraphics[width=0.45\textwidth]{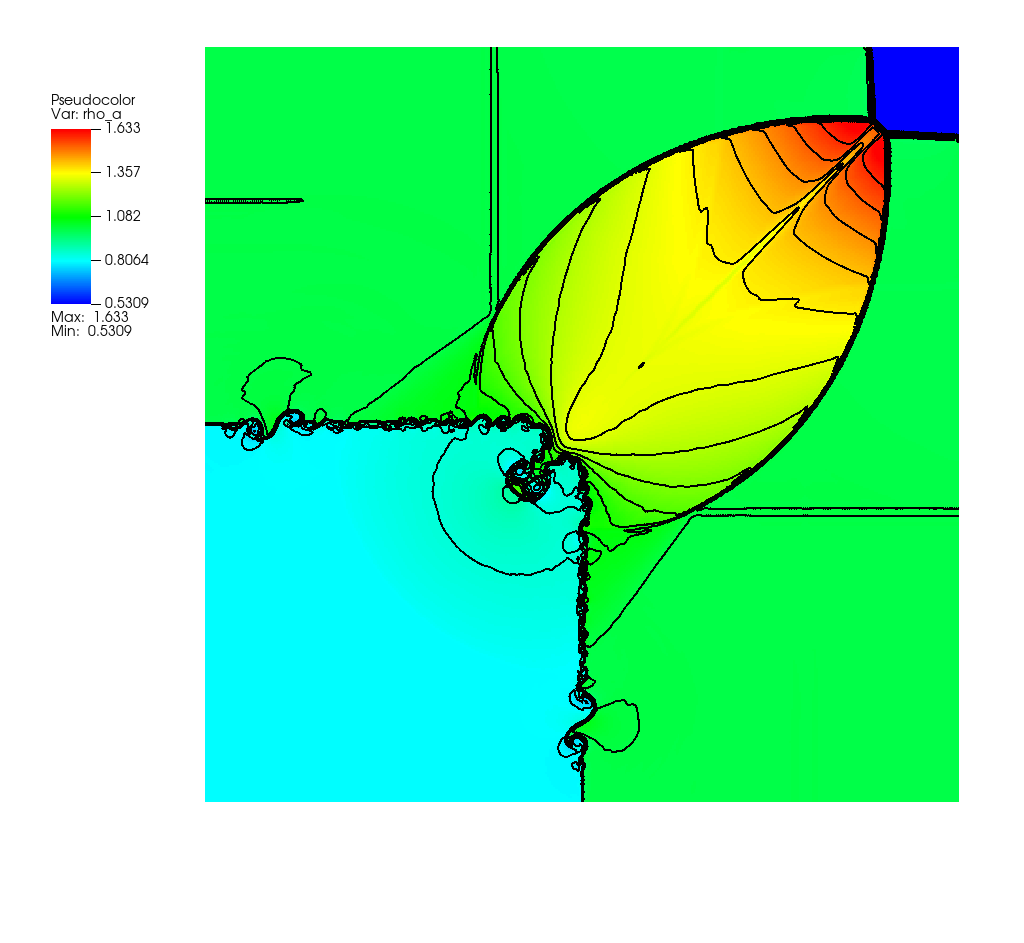}}
\caption{\label{LL13} Lax--Liu test case. Polygonal mesh, 20 isolines CFL=0.3.}
\end{center}
\end{figure}

\paragraph{Lax--Liu problem.}
In the second example of nonlinear Euler equations, we consider the following initial condition{, which corresponds to configuration 13 of \cite{lax1998},}
$$ (\rho, u,v,p)=\left \{\begin{array}{ll}
(\rho_1,u_1,v_1,p_1)=(0.5313,0,0,0.4)& \text{ if } x\geq0\text{ and } y\geq 0,\\
(\rho_2,u_2,v_2,p_2)=(1,0.7276,0,1) & \text{ if } x\leq 0 \text{ and } y\geq 0,\\
(\rho_3,u_3,v_3,p_3)=(0.8,0,0,1)&\text{ if } x\leq 0\text{ and }y\leq 0,\\
(\rho_4,u_4,v_4,p_4)=(1,0, 0.7276,1) &\text{ if } x\geq0\text{ and } y\leq 0,
\end{array}\right .
$$
prescribed in the computational domain $[-2,2]^2$. {All boundary conditions are set to Neuman's.} The states 1 and 2 are separated by a shock. The states 2 and 3 are separated by a slip line. The states 3 and 4 are separated by a steady slip-line. The states  4 and 1 are separated by a shock. {The mesh corresponds to a $400\times400$ cells, i.e. with $476801$ DoFs.} The final time is set to $t_f=1$.  The obtained results of the density field are plotted in Figure \ref{LL13} for the {a-posteriori limiting and convex limiting} versions. The results are qualitatively coherent with those obtained in the literature \cite{lax1998}. However, what can be noted is that the slip lines separating the states 2/3 and 3/4 have a \red{different structure} for the two methods: the {convex limiting} scheme leads to vortices while {while they do not appear when using the a-posteriori limiting procedure}. In our opinion, this is an indication that the {convex limiting} scheme is less dissipative than the {a-posteriori limiting} one: it is very well known that slip lines are not stable structures. This is what can be observed, and the occurrence of this phenomenon in numerical simulation is an indication of a reduced numerical dissipation. We also note that the {convex limiting} solution is smoother than the {a-posteriori limiting} one.

\begin{figure}[!h]
\begin{center}
\subfigure[]{\includegraphics[width=0.45\textwidth]{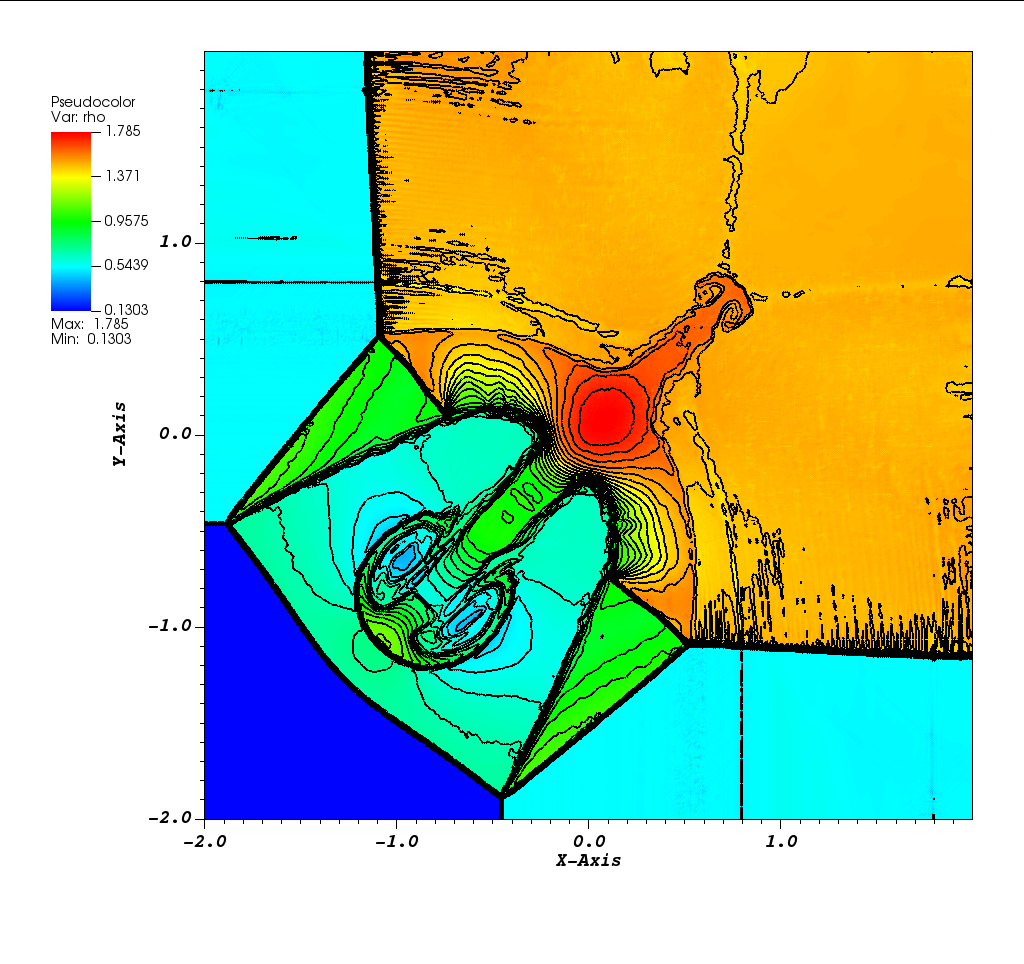}}
\subfigure[]{\includegraphics[width=0.45\textwidth]{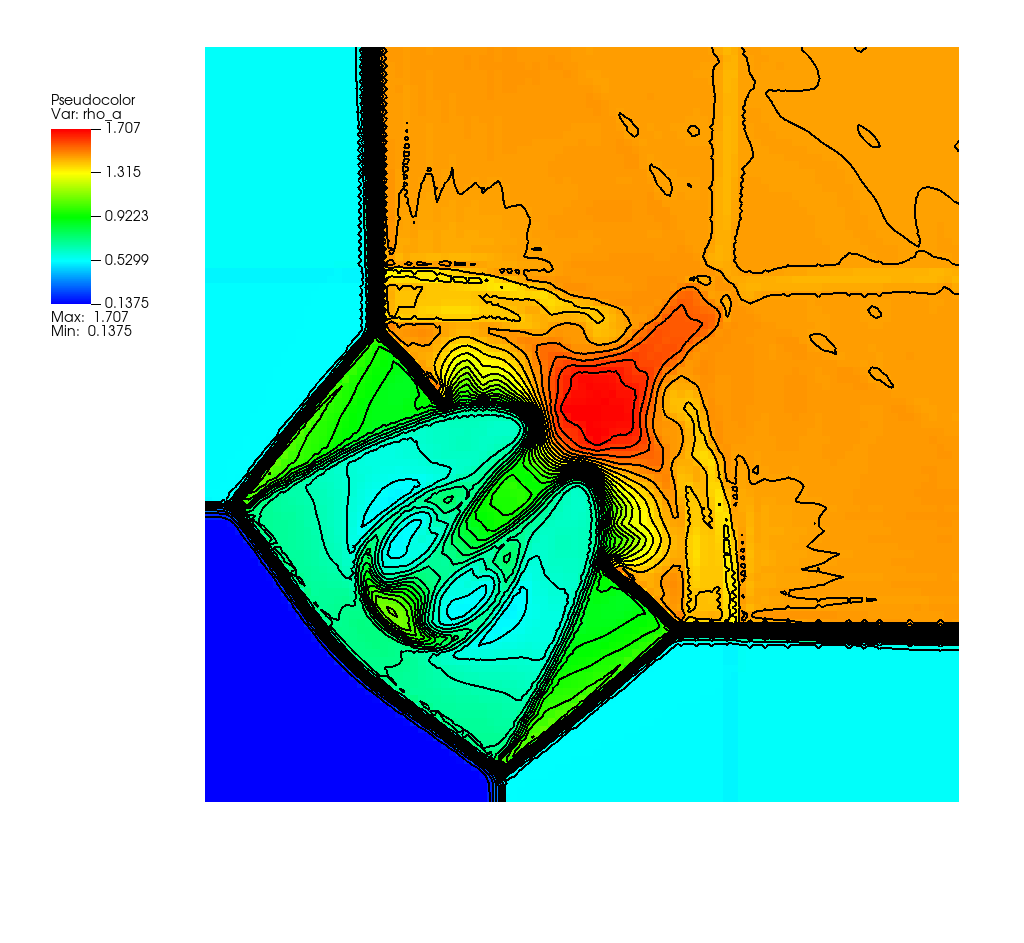}}
\subfigure[]{\includegraphics[width=0.45\textwidth]{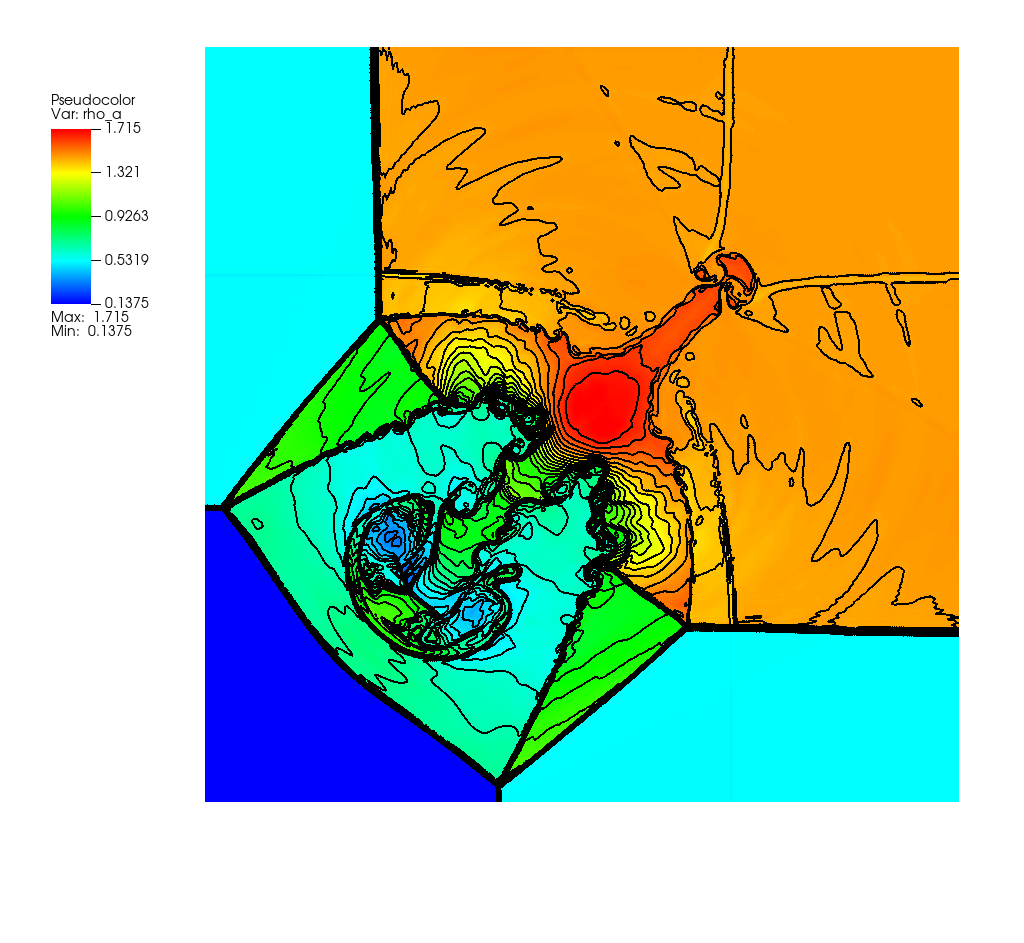}}
\caption{\label{fig:KT} Kurganov-Tadmor test case. Average density iso lines (30 isolines). (a): Mood with $400\times 400$, (b): BP with $100\times 100$, (c): BP with $400\times 400$ CFL=0.3.}
\end{center}
\end{figure}

\paragraph{Kurganov--Tadmor problem.} 
In the third example of nonlinear Euler equations, {which is also known as Configuration 3 of Lax \& Liu in \cite{lax1998}}, the initial condition is 
$$ (\rho, u,v,p)=\left \{\begin{array}{ll}
(\rho_1,u_1,v_1,p_1)=(1.5, 0, 0, 1.5)& \text{ if } x\geq1\text{ and } y\geq 1,\\
(\rho_2,u_2,v_2,p_2)=(0.5323, 1.206, 0, 0.3) & \text{ if } x\leq 1 \text{ and } y\geq 1,\\
(\rho_3,u_3,v_3,p_3)=(0.138, 1.206, 1.206, 0.029)&\text{ if } x\leq 1\text{ and }y\leq 1,\\
(\rho_4,u_4,v_4,p_4)=(0.5323, 0, 1.206, 0.3) &\text{ if } x\leq1\text{ and } y\leq 1.
\end{array}\right .
$$
Here, the four states are separated by shocks. The domain is $[-2,2]^2$. The solution at $t_f=3$ is displayed in Figure \ref{fig:KT}.
{Two meshes are used. One with $100\times 100$ cells, i.e. with $120 312$ DoFs, and a second one with  $400\times400$ cells, i.e. with $481601$ DoFs.  We see, on figure \ref{fig:KT}-(a) and (b),  that the $100\times 100$ BP solutions looks similar to the $400\times 400$ MOOD, but is less noisy than the MOOD one. These solutions look very similar to what was obtained in the literature, see e.g. \cite{abgrall2023activefluxtriangularmeshes,KLZ, GKL, WDKL}. The $400\times 400$ BP solution on figure \ref{fig:KT}-(c) has similar shock structures, but the jet in the middle has a much more complicated one: again this is because the slip lines are unstable.}

From the results we have obtained, it seems clear that the MOOD technique, at least in the way we have implemented it, is surpassed by our BP method. For this reason, we will only consider the Bound preserving scheme for the next test case.

\begin{figure}[!h]
\begin{center}
\subfigure[]{\includegraphics[width=0.45\textwidth]{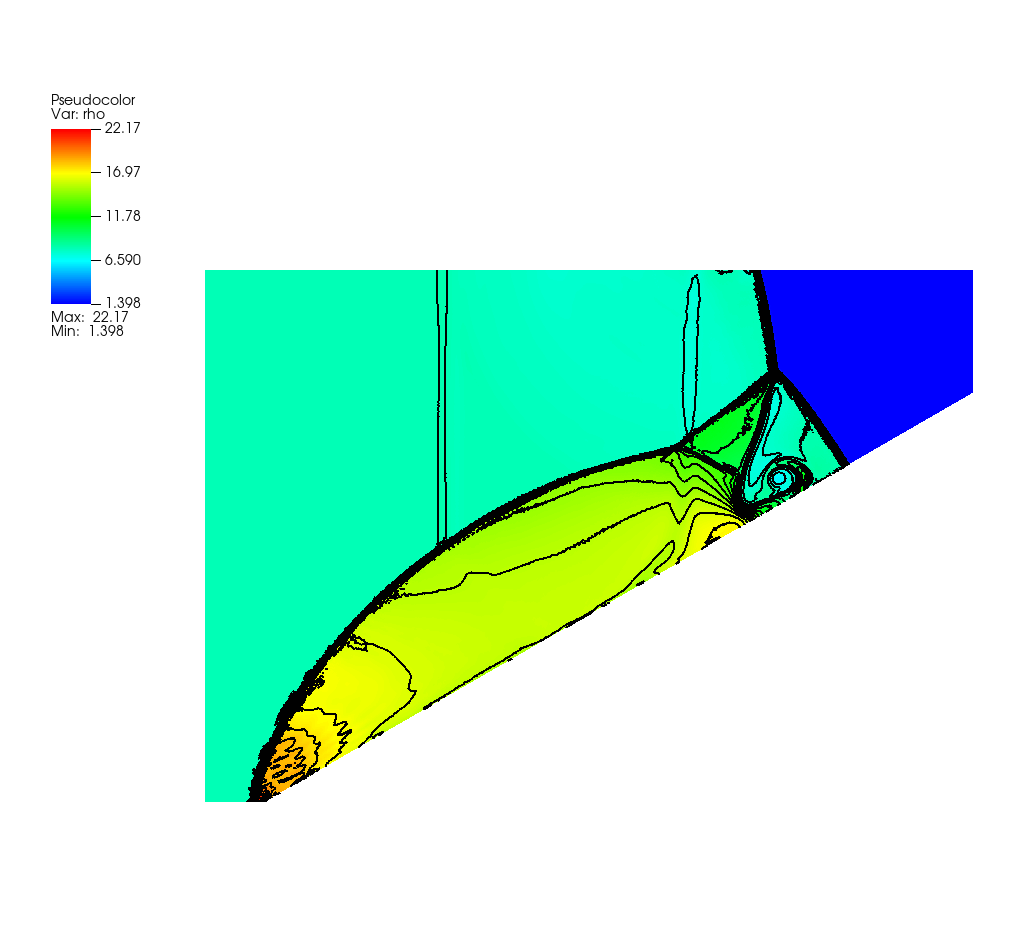}}
\subfigure[]{\includegraphics[width=0.45\textwidth]{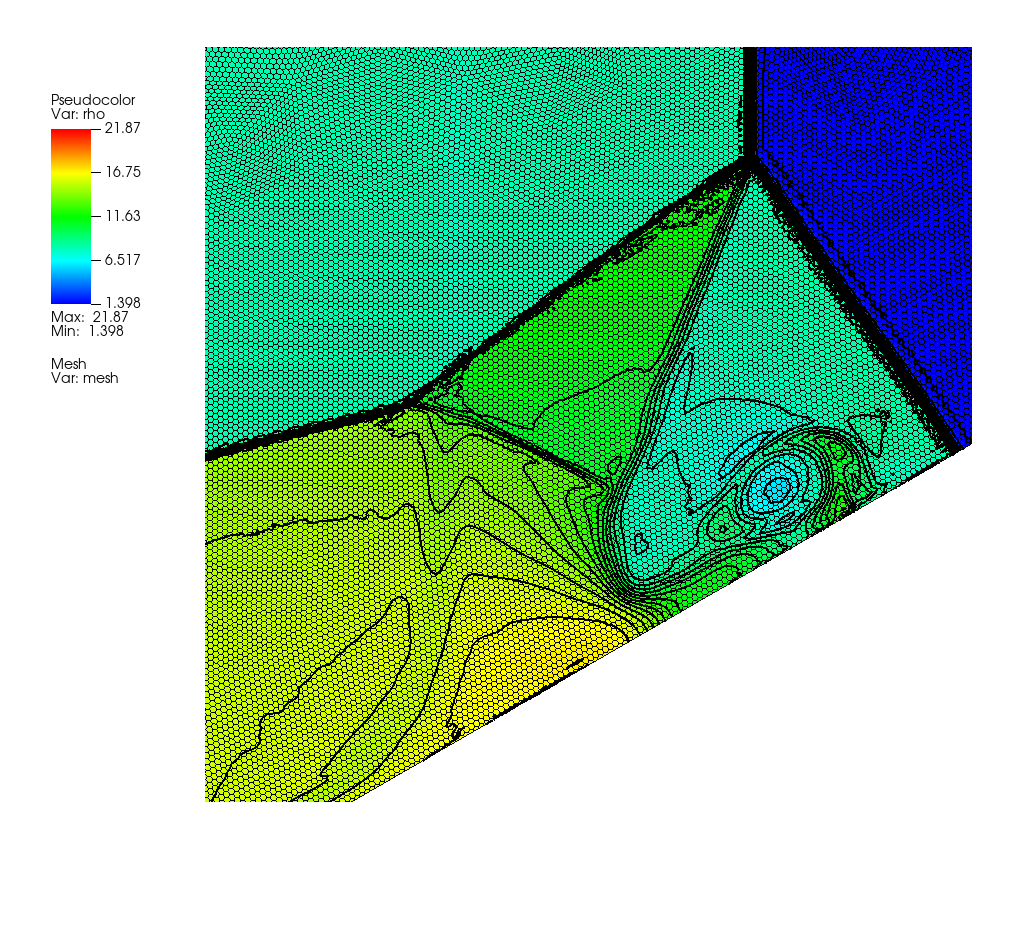}}
\subfigure[]{\includegraphics[width=0.45\textwidth]{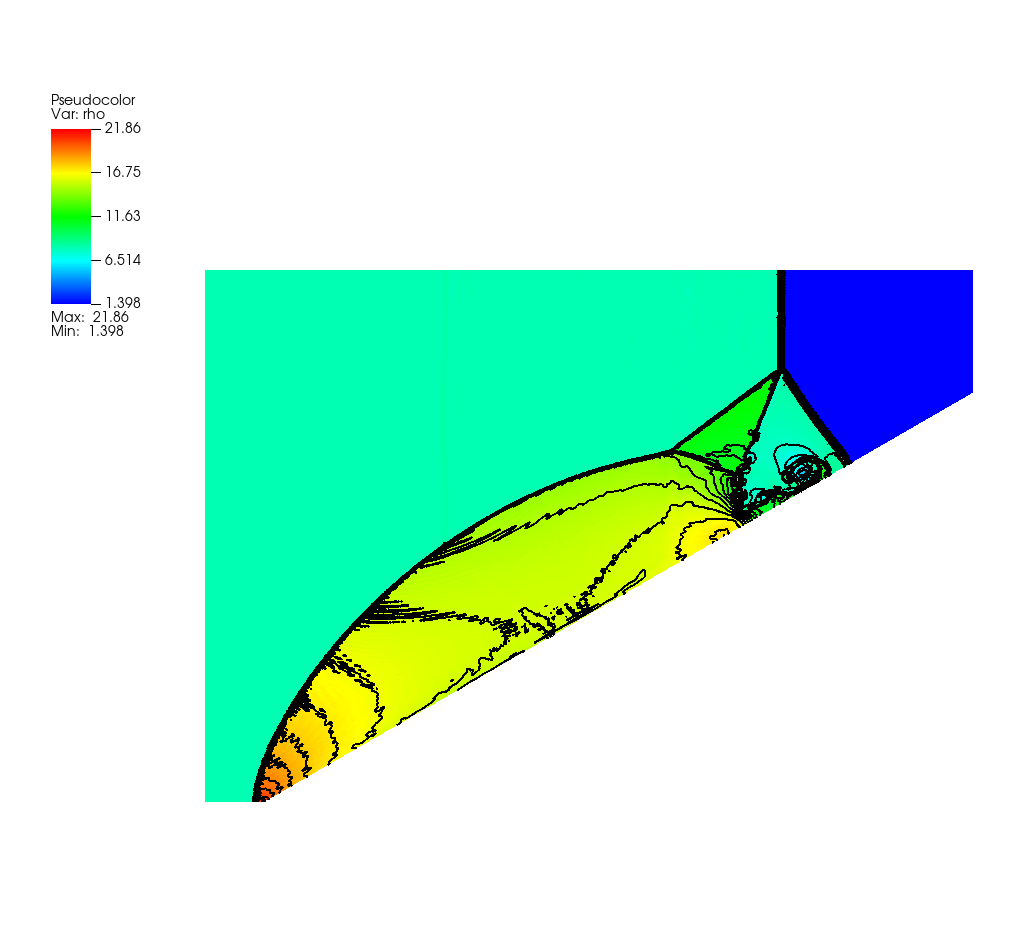}}
\subfigure[]{\includegraphics[width=0.45\textwidth]{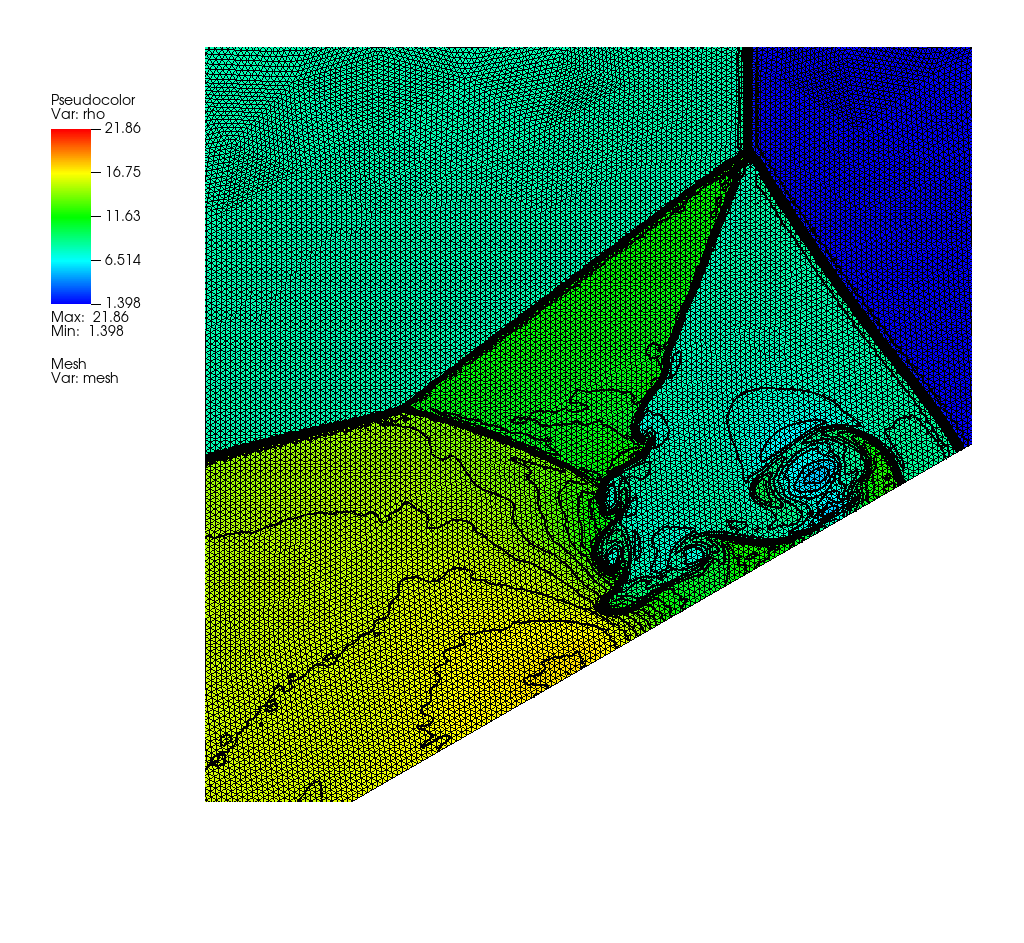}}
\caption{\label{fig:dmr}DMR test case. Density (point values) (a): Polygonal mesh, (b) Polygonal mesh and zoom, (c) Triangular mesh, (d) Triangular mesh and zoom. 30 isolines, CFL=30. The polygonal mesh is constructed from the triangular one by agglomeration.}
\end{center}
\end{figure}

\paragraph{DMR case.}
The final test case is the double Mach reflection problem. The domain consists of {the box $[-0.25, 3] \times [0, 2]$ from which a ramp of 30 degrees is subtracted starting at $x=0$}. 
The initial condition corresponds to a Mach 10 shock in a quiescent flow, where the pressure is chosen to be $p=1$ and the density is $\rho=\gamma$ with $\gamma=1.4$. The mesh has $1,555,600$ point DoFs and $776,501$ {triangles}.  This corresponds to a resolution of $h=1/N$ with $N\approx 350$. We also have used this mesh to build a polygonal. For $t_f=0.18$, a zoom of the density is represented in Figure \ref{fig:dmr}, as well as the computational meshes. This allows to get an idea of the width of the contact line and the shocks.
We notice the appearance of a roll-up structure on the slip line. This is more pronounced for the triangular case, simply because the elements are smaller (the polygons are roughly speaking twice as large as the triangles).

\section{Conclusion and perspectives}
We have presented a way to generalize the method of \cite{abgrall2023activefluxtriangularmeshes} to arbitrary polygonal control volumes using an approximation method inspired by the VEM technology. Here, we only use the VEM to replace the polynomial approximation of \cite{abgrall2023activefluxtriangularmeshes} by a polynomial-free approximation. More precisely, we can discretize the gradients using only the degrees of freedom, and not by using explicitly constructed basis functions. We have shown on several examples that the method preserves its stability and its accuracy, and may have surprisingly good results for the acoustic problem. {Despite the computational cost related to the setting of the VEM framework, which can nevertheless be carried out once and for all in the pre-processing stage, the novel scheme can handle arbitrarily shaped polygonal meshes.}

However, this is not the end of the story. In order to be able to compute solutions with shocks, we need to introduce, as usual, some nonlinear mechanism. This is done here, as in the above mentioned paper, by using either the MOOD paradigm or a bound preserving technique. We have relied on a rather crude implementation of the MOOD  method (in fact exactly the same as in \cite{abgrall2023activefluxtriangularmeshes}), but the results appear to be slightly less satisfactory. For this reason, we also have developed a Bound preserving method by blending a first order bound preserving algorithm with a high order one.For scalar problems, our formula are similar to those of \cite{Vilar_DGFV} and \cite{Duan_AF}. For the Euler case, \red{it turns} out that optimal blending factors can be obtained thanks to the  geometric interpretation  of the invariant domain, interpretation given in \cite{wu2023geometric}. We also intend to extend the method to higher than formally third order of accuracy following \cite{BarsukowAbgrall}. Mesh adaption will also be the topic of future work. In that direction, we notice \cite{Calhoun2023-ms} where an AMR strategy is adopted.
\section*{Acknowledgments}
The work of YL was supported by UZH Postdoc Grant, 2024 / Verf\"{u}gung Nr. FK-24-110 and SNFS grant 200020$\_$204917 ``Solving advection dominated problems with high order schemes with polygonal meshes: application to compressible and incompressible flow problems''. WB acknowledges research funding by the Italian Ministry of University and Research (MUR) in the framework of the PRIN 2022 project No. 2022N9BM3N.

\red{We also take this opportunity to thank the reviewers whose constructive criticisms have led to many improvement to this paper.}

\bibliographystyle{unsrt}
\bibliography{paper}
\appendix
\input{VEM_for_Dummies}
\input{N}

 \end{document}

%% file: command.tex
\usepackage{fullpage}

\usepackage{amsmath,amssymb,amsfonts,amsthm}
\usepackage{mathrsfs}
\usepackage{color,xcolor}
\usepackage{ifpdf}
\usepackage{psfrag}
\usepackage{graphicx,graphics,subfigure,epsfig}
\usepackage{epstopdf}
\usepackage{url}
\usepackage[colorlinks=true]{hyperref}
\usepackage{xspace}
\usepackage{bm}
\usepackage{ulem}
\usepackage{lineno}

\usepackage{extarrows}

\usepackage{csquotes}

\newtheorem{theorem}{Theorem}[section]

\newtheorem{remark}[theorem]{Remark}

\newcommand{\dpar}[2]{\dfrac{\partial #1}{\partial #2}}

\newcommand{\bbx}{{\mathbf {x}}}

\newcommand{\bbz}{{\mathbf {z}}}
\newcommand{\bbf}{{\mathbf {f}}}
\newcommand{\bbPhi}{{\mathbf {\Phi}}}
\newcommand{\R}{\mathbb R}
\newcommand{\N}{\mathbb N}
\newcommand{\K}{\mathbb K}
\renewcommand{\P}{\mathbb P}
\newcommand{\bbu}{{\mathbf {u}}}
\newcommand{\bbv}{{\mathbf {v}}}
\newcommand{\bbw}{{\mathbf {w}}}
\newcommand{\bbn}{{\mathbf {n}}}
\newcommand{\bsa}{{\bm {\alpha}}}
\newcommand{\bsb}{{\bm {\beta}}}
\newcommand{\hbbf}{\hat{\mathbf{f}}}
\newcommand{\lbbf}{\hat{\mathbf{f}}^{\rm LO}}
\newcommand{\Hbbf}{\hat{\mathbf{f}}^{\rm HO}}
\newcommand{\bbA}{\mathbf{A}}
\newcommand{\Id}{\text{Id }}
\newcommand{\bbN}{\mathbf{N}}
\newcommand{\bba}{\mathbf{a}}
\newcommand{\bbU}{\mathbf{U}}
\newcommand{\bbV}{\mathbf{V}}
\newcommand{\FF}{\mathcal F}
\newcommand{\dt}{\Delta t}
\newcommand{\sigmai}{{\sigma_i}}

\newcommand{\bbb}{\mathbf{b}}

%% file: intro.tex
\section{Introduction}
There exist many methods that \red{allow} to compute the solution of hyperbolic problems: finite volume including high order WENO ones, finite difference, continuous finite element (CFE), discontinuous Galerkin (dG) methods. A recent compilation is contained in \cite{Handbook1,Handbook2}. All these methods use a wide variety of meshes structures, so that one can wonder why still working on different ones. This is true, but the answer is not complete. Besides the many remaining problems (such as those related to structure preservation), we also believe that ease of implementation is important, as well as memory footprint. One way to reduce the memory footprint is to use a globally continuous representation of the solution, like what happens for CFE methods. Maybe more important is to develop methods that are friendly to mesh refinement. Mesh refinement is feasible for discontinuous finite element methods, but then one has to consider hanging nodes. However, $p$-refinement is very doable. The existence of hanging nodes complexifies the implementation, especially in 3D. In the case of CFE methods, mesh refinement can be very efficiently done but one has to change the mesh topology, see \cite{alauzet} for example, and $p$-refinement is complicated.

We do not believe there is a perfect solution, but we believe it is important to develop methods that have the potential of flexibility. Looking at the recent literature, some sort of compromise might be found by taking into account ideas coming from the Virtual Finite Element (VEM) and {"High Order Hybrid" (HHO)}  communities \cite{Acta,droniou1,droniou2}. The computational domain is discretized by polygons in 2D and polytopes in 3D. There, the solution can be globally continuous (but one can however relax this constraint if wished), and two types of degrees of freedom (DoFs) are introduced: some are point values on the boundary of the elements, and internal degrees in the form of moments are also considered. The functional representation is not made via polynomials (see appendix \ref{sec:VEM} for details), but it is polynomial on each edge or face of the elements. The polynomial degree can be edge/face dependent. This feature allows $h-p$ refinement  {without breaking the continuity constraint}. For example in 2D (to simplify), one can cut an edge in two or more, and represent the same polynomial function with sets of degree of freedom (DoF) attached to the new edges. The internal degrees of freedom (DoFs) are not affected by this. If an element is cut into several sub-elements, one can adapt this procedure, by also taking into account the special constraints of hyperbolic PDEs such as local conservation: this will be the topic of future researches.

Coming back to VEM type approximation, one can imagine a variational formulation, such as in \cite{daVeiga} for example. The problem is that one will have a mass matrix, as in classical finite element methods, that one will have to invert. Here, we chose to have diagonal mass matrix, as for dG, and a solution (but with discontinuous elements \ldots) can be found in \cite{Boscheri}.
One possible inspiration to get simultaneously a globally continuous approximation and diagonal mass matrix can be found in the so-called Active Flux methods.

The Active Flux method was initially introduced in \cite{AF1,AF2,AF3,Maeng,He} for solving hyperbolic problems on triangular unstructured meshes. The numerical solution is approximated employing the DoFs of the quadratic polynomials, which are all lying on the boundary of the elements, supplemented by another DoF: the average of the solution. This leads to a potentially third order accurate method, with an approximation that is globally continuous. The time integration relies on a back-tracing of characteristics. The extension of this scheme to square elements, in a finite difference fashion, has then been {done} in \cite{AF4,AF5}. For these original (fully-discrete) Active Flux methods, the evolution operators for point values are critical. Exact evolution operators based on characteristic methods have been developed for linear hyperbolic equations (see, e.g., \cite{EymannRoe2013,Roe2015,Barsukow,Calhoun2023-ms}). For nonlinear systems, several approximate evolution operators have been introduced, including those for Burgers' equation \cite{AF1,AF2,AF5}, the 1D compressible Euler equations \cite{AF2,AF4,AF5}, and hyperbolic balance laws \cite{AF5,Barsukow2024-tc}.

For nonlinear systems, constructing exact or approximate evolution operators becomes significantly more complex, particularly in multiple spatial dimensions. In \cite{Abgrall_AF}, a different, though very close, point of view is proposed. This work introduces a streamlined approach for evolving point values and combining different formulations (conservative and non-conservative) of a nonlinear hyperbolic system. The updates for both cell averages and point values are formulated in a semi-discrete form, enabling the use of standard Runge--Kutta time stepping algorithms. These routes have been followed in \cite{BarsukowAbgrall}, giving several versions of higher order approximations in one spatial dimension. Later on, the problem of {design limiters for nonlinear problems} has received attention: besides the {\textit{a-posteriori} limiting} technique already employed in \cite{Abgrall_AF} and other papers \cite{Liu2024,Abgrall2024_WBAF}, a direct convex limiting method is being considered, see \cite{Duan_AF, Abgrall_BP_PAMPA}. More recently, the semi-discrete Active Flux method introduced in \cite{Abgrall_AF} has been extended to multidimensional settings. In \cite{abgrall2023activefluxtriangularmeshes}, the method was extended to triangular meshes for hyperbolic conservation laws, where the DoFs are again given by Lagrange point values on the boundary of elements and one average in the element. Both conservative and non-conservative formulations are expressed in terms of conserved variables. In \cite{Duan_AF}, a related problem was studied on Cartesian meshes and this semi-discrete Active Flux method is referred to as a generalized Active Flux scheme. Further, in \cite{Liu_PAMPA_SW2D}, the semi-discrete Active Flux method has been extended to hyperbolic balance laws on triangular meshes. Here, the conservative formulation, given by conserved variables, is used to update the average DoF while the non-conservative formulation given by equilibrium variables updates the point value DoFs. This new hybrid procedure, combining conserved and other variables (e.g., primitive and equilibrium variables), is named the PAMPA (Point-Average-Moment PolynomiAl-interpreted) scheme in \cite{Liu_PAMPA_SW2D}.

In this work, we aim to extend the semi-discrete Active Flux method from \cite{abgrall2023activefluxtriangularmeshes} to more general polygonal meshes and incorporate the MOOD stabilisation and  the convex limiting technique from \cite{Abgrall_BP_PAMPA} to guarantee the invariant domain preserving property in the spirit of \cite{Zhang_MP,Guermond_IDP3,Guermond_IDP,Guermond_IDP2,Hajduk_MCL,Kuzmin_MCL1,Kuzmin_BP, wu2023geometric}. Both methods will  described and compared. The resulting scheme remains named PAMPA, but without knowing the explicit polynomial basis functions.


The format of this paper is as follows. First, we recall the approximation spaces and how we construct the polygonal meshes. Then, we describe the schemes, the high order one and the low order one that we need for the limiting strategies. Then we describe in details the {\textit{a-posteriori} limiting} strategy. Next, we introduce the monolithic convex limiting approach to blend the high- and low-order schemes. Finally, we present a set of numerical examples using triangle, quadrangle and general polygonal meshes, validating the accuracy and the robustness of the novel schemes, and allowing a fair comparison. A conclusion follows. In the appendixes, we recall the essential of VEM approximation.

%% file: pente.tex
\begin{table}[!h]
\begin{center}
\begin{tabular}{|c||cc||cc||cc|}\hline
\multicolumn{7}{|c|}{Average values}\\
\hline
$h$                      &$L^\infty$              & slope   &$L^1$                   & slope  &$L^2$ & slope \\ 
\hline
  $0.224 $            &  $0.334 $             &-            & $9.591\,10^{-3}$ &-          &$3.490\,10^{-2}$&-\\
   $0.112\,10^{-1 }$& $0.197 $              & $0.75$ &$3.261\,10^{-3}$ &$1.55$ &$1.509\,10^{-2}$& $1.71 $ \\
  $ 5.60\,10^{-2}$ & $ 5.418\,10^{-2}$& $1.86$ &$6.385\,10^{-4}$ &$2.35$ &$3.447\,10^{-3}$&$2.65$ \\
  $ 2.80\,10^{-2}$ &  $8.693\,10^{-3}$& $2.63$ &$9.086\,10^{-5}$ &$2.81$ &$5.161\,10^{-4}$&$2.72$ \\
  $ 1.40\,10^{-3}$ &  $ 1.144\,10^{-3}$& $2.92$&$1.169\,10^{-5}$ &$2.95$ &$6.716\,10^{-5}$&$2.64 $\\
      \hline \hline \multicolumn{7}{|c|}{Point values}\\
      \hline
      $h$                &$L^\infty$              & slope     &$L^1$                & slope &$L^2$                  & slope \\ 
\hline
 $0.224 $             &  $0.530$               &-             &$9.591\,10^{-3}$&-          & $3.259\,10^{-2}$& -\\
 $  0.112\,10^{-1}$ & $0.228$                &$0.75$   &$3.261\,10^{-3}$&$1.55$&$1.296\,10^{-2}$ & $1.20$\\
 $  5.60\,10^{-2} $& $5.656\,10^{-2}$  &$1.86$   &$6.385\,10^{-4}$&$2.35$&$2.872\,10^{-3}$ & $2.12$\\
 $  2.80\,10^{-2 }$& $8.881\,10^{-3}$  &$2.63$   &$9.086\,10^{-5}$&$2.81$ &$4.284\,10^{-4}$ & $2.73$\\
 $     1.40\,10^{-3}$&$1.190\,10^{-3}$ &$ 2.92$  &$1.169\,10^{-6}$&$2.95$&$5.680\,10^{-5}$. & $2.94$ \\
      \hline
      \end{tabular}
      \end{center}
      \caption{\label{error:convection:tri}Errors for the average and point values, triangular mesh, rotation problem \eqref{scalar:test}, 1 rotation.}
   \end{table}
   \begin{table}[!h]
\begin{center}
   \begin{tabular}{|c||cc||cc||cc||}\hline
   \multicolumn{7}{|c|}{Average values}\\ \hline
$h$                     &$L^\infty$            & slope     &$L^1$                   & slope   &$L^2$                 & slope \\ 
\hline
$0.224 $             & $0.333  $           &              -&$9.591\; 10^{-3} $&-            & $3.490\;10^{-2} $&-\\
$0.112  $            & $0.197 $            &  $0.75$  &$3.261\; 10^{-3} $& $1.55$ & $1.509\; 10^{-2} $& $1.20 $\\
$ 5.60\;10^{-2}$& $5.418\;10^{-2}$& $1.86$  &$6.385\; 10^{-4} $& $2.35$ & $3.447\; 10^{-3} $& $2.12 $\\
$2.80\;10^{-2}$ & $8.693\;10^{-3}$& $2.63$  &$9.086\; 10^{-5} $& $2.81$ & $5.161\; 10^{-4} $& $2.73 $ \\
$1.40\;10^{-2}$ & $1.144\;10^{-3}$& $2.92$  &$1.169\; 10^{-5} $& $2.96$ & $6.716\; 10^{-5} $ & $2.94$\\
      \hline
\hline
 $h $                    &$L^\infty$              & slope       &$L^1$ & slope &$L^2$ & slope \\ 
\hline\multicolumn{7}{|c|}{Point values}\\
      \hline
  $ 0.224           $  &  $ 0.530 $            &  -            & $7.098\;10^{-3} $&-             &  $3.259\; 10^{-2} $&-\\
  $ 0.112            $  & $0.227  $            & $0.75$   & $2.271\; 10^{-3} $& $1.64$ & $1.296\; 10^{-2 }$& $1.33 $\\
 $   5.60\;10^{-2}$ & $ 5.656\;10^{-2} $& $1.86 $ & $4.367\; 10^{-4} $& $2.37$ & $2.872\; 10^{-3} $& $2.17 $\\
  $  2.80\;10^{-2}$ & $  8.882\; 10^{-3 }$& $2.63$& $6.271\; 10^{-5} $& $2.80$ &  $4.284\; 10^{-4} $&  $2.74 $\\
  $1.40\;10^{-2}$   & $1.198\;10^{-3}$    & $2.92$&$8.488\;  10^{-6} $& $2.88$ & $5.680\; 10^{-5} $ & $2.91$\\
      \hline
      \end{tabular}
            \end{center}
            \caption{\label{error:convection:poly}Errors for the average and point values, {polygonal mesh}, rotation problem \eqref{scalar:test}, 1 rotation.}
       \end{table}

%% file: VEM_for_Dummies.tex
\section{VEM approximation: basic facts}\label{sec:VEM}
 
Any basis of $V_k(P)$ is virtual, meaning that it is not explicitly computed in closed form. Consequently, the evaluation of $v_h(\bbx)$ at some $\bbx\in P$ it is not straightforward. One way to proceed would be to evaluate $\pi(v_h)$, that is the $L^2$ projection of $v_h$ on $\P^k(P)$. Since we want to do it only using the degrees of freedom, it turns out that this is impossible in practice. But, as shown in \cite{hitch,vem}, it is possible to define a space $W_k(P)$ for which computing the $L^2$ projection is feasible.
It is constructed from $V_k(P)$ in two steps. First, we consider the approximation space $\widetilde{V}_k(P)$ given by
$$\widetilde{V}_k(P)=\{ v_h, v_h \text{ is continuous on }\partial P \text{ and } (v_h)_{\partial P}\in \P^k(\partial P); \text{ and }\Delta v_h\in \P^{k-2}(P)\}.$$
For $p\in \N$, let $M^\star_p(P)$ be the vector space generated by the scaled monomial of degree $p$ exactly,
$$m(\bbx)\in M^\star_p(P), \qquad m(\bbx)=\sum_{\bsa, |\bsa|=p} \beta_{\bsa} m_{\bsa}(\bbx).$$
Then, we consider $W_k(P)$, which is the subspace of $\widetilde{V}_k(P)$ defined by 
$$W_k(P)=\{ w_h\in \widetilde{V}_k(P), \langle w_h-\pi^\nabla w_h,q\rangle=0, \quad \forall q\in M^\star_{k-2}(P)\cup M^\star_k(P)\}.$$
In \cite{vem}, it is shown that $\dim V_k(P)=\dim W_k(P)$.

This approximation space is defined as follows. $w_h\in W_k(P)$ if and only if
 \begin{enumerate}
  \item $w_h$ is a polynomial of degree $k$ on each edge $e$ of $P$, that is $(w_h)_{|e}\in \P^k(e)$,
 \item $w_h$ is continuous on $\partial P$,
\item $\Delta w_h\in \P^{k-2}(P)$,
\item $\int_{P} w_h m_{\bsa} \; {\rm d}\bbx=\int_{P} \pi^\nabla w_h m_{\bsa} \; {\rm d}\bbx$ for $|\bsa|=k-1, k$.
\end{enumerate}
The degrees of freedom are the same as in $V_k(P)$:
\begin{enumerate}
 \item The value of $w_h$ on the vertices of $P$,
 \item On each edge of $P$, the value of $v_k$ at the $k-1$ internal points of the $k+1$ Gauss--Lobatto points on this edge,
 \item The moments up to order $k-2$ of $w_h$ in $P$,
 $$m_{\bsa}(w_h):=\dfrac{1}{\vert P\vert }\int_{P} w_h m_{\bsa}\; d\bbx, \qquad |\bsa|\leq k-2.$$
 \end{enumerate}

The $L^2$ projection of $w_h$ is computable. For any $\bsb$, if $\pi^0(w_h)=\sum_{\bsa, |\bsa|\leq k} s_{\bsa} m_{\bsa}$, we have
$$\langle \pi^0(w_h), m_{\bsb}\rangle=\sum_{\bsa, |\bsa|\leq k} s_{\bsa} \langle m_{\bsb}, m_{\bsa}\rangle=\langle w_h, m_{\bsb}\rangle.$$
The left hand side is computable since the inner product $\langle m_{\bsb}, m_{\bsa}\rangle$ only involves monomials. We need to look at the right hand side. If $|\bsa|\leq k-2$, $\langle w_h, m_{\bsb}\rangle=|P| m_{\bsb}(w_h)$ and if $|\bsa|=k-1$ or $k$, we have
$$\langle w_h, m_{\bsb}\rangle=\langle \pi^\nabla w_h, m_{\bsa}\rangle,$$ which is computable from the degrees of freedom only.

We note that if $w_h\in W_k(P)$, then $m_{\bsa}(\pi^0 w_h)=m_{\bsa}(w_h)$. Indeed
$$m_{\bsa}(\pi^0 w_h)=\frac{1}{|P|} \langle \pi^0 w_h ,m_{\bsa} \rangle =\frac{1}{|P|} \langle  w_h ,m_{\bsa} \rangle$$ by construction. This is not true for $\pi^\nabla$ in $V_k(P)$.
However, for $k\leq 2$, $V_k(P)=W_k(P)$

The last remark is that, since $\P^k(P)\subset V_k(P)$ and $\P^k(P)\subset W_k(P)$, the projections of $C^{k+1}(P)$ onto $V_k(P)$ and $W_k(P)$ defined by the degrees of freedom is $k+1$-th order accurate.

%% file: N.tex
\section{The N matrix}\label{append:N}
In this section, we show that $\bbN_\sigma=(\sum_{P,\sigma\in P}K^+_{\sigma})^{-1}$ {defined in \eqref{N}} has a meaning when
$$K_\sigma=A_xn^P_x+A_yn^P_y,\quad \forall~ \bbn^P=(n^P_x,n^P_y),$$ where $A_x\in\R^{m\times m}$ and $A_y\in\R^{m\times m}$ are the Jacobians of a symmetrizable system, the vectors $\bbn^P$ sum up to $0$,
$$\sum_{P}\bbn^P=\mathbf 0$$ and 
$$\bbN_\sigma^{-1}=\sum_{P} K^+_{\sigma}.$$

There exists $A_0$ a symmetric positive definite matrix such that $A_0A_x$ and $A_0A_y$ are symmetric. For the Euler equations, this is the Hessian of the entropy. From this we see that
$$\big (A_0K_{\sigma}\big )^+=A_0\big ( K_{\sigma}\big )^+.$$
This comes from
$$
A_0^{-1/2}\big (A_0 K_{\sigma}\big ) A_0^{-1/2}=A_0^{1/2}K_{\sigma}A_0^{-1/2}.$$

Second, we see that the eigenvectors of $K_\sigma$ are orthogonal for the metric defined by $A_0$,
$$\langle \bbu,\bbv\rangle_{A_0}=\bbu^T A_0\bbv$$

Last, we can split $\R^m$ as $\R^m=\bbU\oplus \bbV$ where $\bbU$ is the vector space generated by the vectors that are eigenvectors of $A_x$ and $A_y$, and $\bbV$ its orthogonal for the above metric. For the Euler equations, $\bbU$ is generated by the eigenvector associated to the transport of entropy.

We see that for any $\bbx\in \R^m$,  defining a orthonormal basis of $\bbU$ as $(\bbu_i)_{i=1,\ldots, m}$ and writing $\bbx=\sum_{i=1}^{m}\langle\bbx,\bbu_i\rangle \bbu_i +\bbv$, we have
$$K_\sigma^+\bbx = \sum_{i=1}^m\lambda_{i}^+\langle\bbx,\bbu_i\rangle \bbu_i+K_\sigma^+\bbv,$$
where $K_\sigma^+\bbv\in \bbV$.
Hence,
$$\big ( \sum_{P,\sigma\in P} K^+_{\sigma}\big )\bbx=
\sum_{i=1}^m\big (\sum_P  \lambda_{i}^+\big ) \langle\bbx,\bbu_i\rangle \bbu_i+ \sum_P K^+_{\sigma}\bbv .$$

Then we see that $\sum\limits_{P,\sigma\in P} K^+_{\sigma}$ is invertible on $\bbV$ because for any $\bbv$, 
$$\langle \bbv, \sum_{P,\sigma\in P} K^+_{\sigma}\bbv\rangle \geq 0$$
and if it were $0$ for some $\bbv\neq \mathbf{0}$ with $\bbv\in V$, then we would have for all $\sigma\in P$, 
$\langle \bbv,K^+_{\sigma}\bbv\rangle=0$. This implies that $\big ( K^+_{\sigma}\big)^{1/2}\bbv=0$, that is $\bbv\in \bbV$ is in the null space of $\big ( K^+_{\sigma}\big)^{1/2}$. But the null space of $\big (K^+_{\sigma}\big)^{1/2}$ is that of $K^+_{\sigma}$ so that 
$$K^+_{\sigma}\bbv=0,$$
that is $\bbv\neq\mathbf{0}$ would be a common eigenvector of all the $\K_\sigma$, this is impossible: $\bbv=0$. This shows that the restriction of $\sum_{P,\sigma\in P}K^+_{\sigma}$, that we still denote by $\sum_{P,\sigma\in P} K^+_{\sigma}$ is invertible on $\bbV$.

In the end we get
$$\bbN_{\sigma}K^+_{\sigma}\bbx=\sum_{i=1}^m \dfrac{\lambda_i^+}{\sum_P \lambda_i^+}\langle\bbx,\bbu_i\rangle \bbu_i +
\big (\sum_P K^+_{\sigma}\big )^{-1}K^+_{\sigma}\bbv.$$

If we consider $\text{sign}(K_{\sigma})$ instead, the proof is the same.

%% file: paper.bbl
\begin{thebibliography}{10}

\bibitem{AF1}
T.A. Eyman and P.L. Roe.
\newblock Active flux.
\newblock 49th AIAA Aerospace Science Meeting, 2011.

\bibitem{Abgrall_AF}
R.~Abgrall.
\newblock A combination of residual distribution and the active flux
  formulations or a new class of schemes that can combine several writings of
  the same hyperbolic problem: application to the 1d {Euler} equations.
\newblock {\em Commun. Appl. Math. Comput.}, 5(1):370--402, 2023.

\bibitem{abgrall2023activefluxtriangularmeshes}
R.~Abgrall, J.~Lin, and Y.~Liu.
\newblock Active flux for triangular meshes for compressible flows problems.
\newblock {\em Beijing Journal of Pure and Applied Mathematics}, 2025.
\newblock in press, also Arxiv preprint 2312.11271.

\bibitem{Abgrall_BP_PAMPA}
R.~Abgrall, M.~Jiao, Y.~Liu, and K.~Wu.
\newblock Bound preserving {P}oint-{A}verage-{M}oment
  {P}olynomi{A}l-interpreted ({PAMPA}) scheme: one-dimensional case.
\newblock {\em submitted}, 2024.
\newblock Arxiv: 2410.14292.

\bibitem{Handbook1}
R.~Abgrall and C.W. Shu, editors.
\newblock {\em {Handbook of Numerical Methods for Hyperbolic Problems: Basic
  and Fundamental Issues}}, volume~17 of {\em Handbook of Numerical Analysis}.
\newblock North-Holland, 2016.

\bibitem{Handbook2}
R.~Abgrall and C.W. Shu, editors.
\newblock {\em {Handbook of Numerical Methods for Hyperbolic Problems: Applied
  and Modern Issues}}, volume~18 of {\em Handbook of Numerical Analysis}.
\newblock North-Holland, 2016.

\bibitem{alauzet}
A.~Belme, A.~Dervieux, and F.~Alauzet.
\newblock Time accurate anisotropic goal-oriented mesh adaptation for unsteady
  flows.
\newblock {\em J. Comput. Phys.}, 231(19):6323--6348, 2012.

\bibitem{Acta}
Louren{\c{c}}o Beir{\~a}o~da Veiga, Franco Brezzi, L.~Donatella Marini, and
  Alessandro Russo.
\newblock The virtual element method.
\newblock {\em Acta Numerica}, 32:123--202, 2023.

\bibitem{droniou1}
Louren{\c{c}}o Beir{\~a}o~da Veiga, Franco Dassi, Daniele~A. Di~Pietro, and
  J{\'e}r{\^o}me Droniou.
\newblock Arbitrary-order pressure-robust {DDR} and {VEM} methods for the
  {Stokes} problem on polyhedral meshes.
\newblock {\em Comput. Methods Appl. Mech. Eng.}, 397:31, 2022.
\newblock Id/No 115061.

\bibitem{droniou2}
Daniele~A. Di~Pietro and J{\'e}r{\^o}me Droniou.
\newblock From {Finite} {Elements} to {Hybrid} {High}-{Order} methods.
\newblock Preprint, {arXiv}:2503.00425 [math.{NA}] (2025), 2025.

\bibitem{daVeiga}
L.~Beir{\~a}o~da Veiga, F.~Dassi, and S.~G{\'o}mez.
\newblock {SUPG}-stabilized time-{DG} finite and virtual elements for the
  time-dependent advection-diffusion equation.
\newblock {\em Comput. Methods Appl. Mech. Eng.}, 436:31, 2025.
\newblock Id/No 117722.

\bibitem{Boscheri}
Walter Boscheri and Giulia Bertaglia.
\newblock Local virtual element basis functions for space-time discontinuous
  {Galerkin} schemes on unstructured {Voronoi} meshes.
\newblock {\em Commun. Comput. Phys.}, 36(2):348--388, 2024.

\bibitem{AF2}
T.A. Eyman and P.L. Roe.
\newblock Active flux for systems.
\newblock 20 th AIAA Computationa Fluid Dynamics Conference, 2011.

\bibitem{AF3}
T.A. Eyman.
\newblock {\em Active flux}.
\newblock PhD thesis, University of Michigan, 2013.

\bibitem{Maeng}
Jungyeoul Maeng.
\newblock {\em On the Advective Component of Active Flux Schemes for Nonlinear
  Hyperbolic Conservation Laws}.
\newblock PhD thesis, {Applied and Interdisciplinary Mathematics, University of
  Michigan}, 2017.
\newblock {https://deepblue.lib.umich.edu/handle/2027.42/138695}.

\bibitem{He}
Fanchen He.
\newblock {\em Towards a New-generation Numerical Scheme for the Com- pressible
  Navier-Stokes Equations with the Active Flux Method}.
\newblock PhD thesis, {Applied and Interdisciplinary Mathematics, University of
  Michigan}, 2021.
\newblock https://deepblue.lib.umich.edu/handle/2027.42/169687.

\bibitem{AF4}
C.~Helzel, D.~Kerkmann, and L.~Scandurra.
\newblock A new {ADER} method inspired by the active flux method.
\newblock {\em Journal of Scientific Computing}, 80(3):35--61, 2019.

\bibitem{AF5}
W.~Barsukow.
\newblock The active flux scheme for nonlinear problems.
\newblock {\em J. Sci. Comput.}, 86(1):Paper No. 3, 34, 2021.

\bibitem{EymannRoe2013}
T.~A. Eymann and P.~L. Roe.
\newblock Multidimensional {A}ctive {F}lux schemes.
\newblock In American~Institute of~Aeronautics and Astronautics, editors, {\em
  21st AIAA Computational Fluid Dynamics Conference}, 2013.

\bibitem{Roe2015}
D.~Fan and P.~L. Roe.
\newblock Investigations of a new scheme for wave propagation.
\newblock In American~Institute of~Aeronautics and Astronautics, editors, {\em
  22nd AIAA Computational Fluid Dynamics Conference}, 2015.

\bibitem{Barsukow}
Wasilij Barsukow, Jonathan Hohm, Christian Klingenberg, and Philip~L. Roe.
\newblock The active flux scheme on {Cartesian} grids and its low {Mach} number
  limit.
\newblock {\em J. Sci. Comput.}, 81(1):594--622, 2019.

\bibitem{Calhoun2023-ms}
D.~Calhoun, E.~Chudzik, and C.~Helzel.
\newblock The {C}artesian grid {A}ctive {F}lux method with adaptive mesh
  refinement.
\newblock {\em J. Sci. Comput.}, 94:54, 2023.

\bibitem{Barsukow2024-tc}
W.~Barsukow and J.~P. Berberich.
\newblock A well-balanced {A}ctive {F}lux method for the shallow water
  equations with wetting and drying.
\newblock {\em Commun. Appl. Math. Comput.}, 6:2385--2430, 2024.

\bibitem{BarsukowAbgrall}
R.~Abgrall and W.~Barsukow.
\newblock Extensions of active flux to arbitrary order of accuracy.
\newblock {\em ESAIM, Math. Model. Numer. Anal.}, 57(2):991--1027, 2023.

\bibitem{Liu2024}
Y.~Liu and W.~Barsukow.
\newblock An arbitrarily high-order fully well-balanced hybrid finite
  element-finite volume method for a one-dimensional blood flow model.
\newblock {\em SIAM J. Sci. Comput.}, 47(4):a2041--a2073, 2025.

\bibitem{Abgrall2024_WBAF}
R.~Abgrall and Y.~Liu.
\newblock A new approach for designing well-balanced schemes for the shallow
  water equations: a combination of conservative and primitive formulations.
\newblock {\em SIAM J. Sci. Comput.}, 46:A3375--A3400, 2024.

\bibitem{Duan_AF}
J.~Duan, W.~Barsukow, and C.~Klingenberg.
\newblock Active flux methods for hyperbolic conservation laws---flux vector
  splitting and bound-preserving.
\newblock {\em SIAM Journal on Scientific Computing}, 2024.
\newblock Arxiv: 2411.00065.

\bibitem{Liu_PAMPA_SW2D}
Y.~Liu.
\newblock Well-balanced {P}oint-{A}verage-{M}oment {P}olynomi{A}l-interpreted
  ({PAMPA}) methods for shallow water equations on triangular meshes.
\newblock {\em arXiv preprint}, 2024.
\newblock arXiv: 2409.12606.

\bibitem{Zhang_MP}
X.~Zhang, Y.~Xia, and C.-W. Shu.
\newblock Maximum-principle-satisfying and positivity-preservinghigh order
  discontinuous galerkin schemes for conservation laws on triangular meshes.
\newblock {\em J. Sci. Comput.}, 50:29--62, 2012.

\bibitem{Guermond_IDP3}
J.-L. Guermond, B.~B.~Popov, and I.~Tomas.
\newblock Invariant domain preserving discretization independent schemes and
  convex limiting for hyperbolic systems.
\newblock {\em Comput. Method. Appl. M.}, 347:143--175, 2019.

\bibitem{Guermond_IDP}
J.-L. Guermond, M.~Nazarov, B.~Popov, and I.~Tomas.
\newblock Second-order invariant domain preserving approximation of the euler
  equations using convex limiting.
\newblock {\em SIAM J. Sci. Comput.}, 40:A3211--A3239, 2018.

\bibitem{Guermond_IDP2}
J.-L. Guermond and B.~Popov.
\newblock Invariant domains and first-order continuous finite
  elementapproximation for hyperbolic systems.
\newblock {\em SIAM J. Numer. Anal.}, 54:2466--2489, 2016.

\bibitem{Hajduk_MCL}
H.~Hajduk.
\newblock Monolithic convex limiting in discontinuous galerkin discretizations
  of hyperbolic conservation laws.
\newblock {\em Comput. Math. Appl.}, 87:120--138, 2021.

\bibitem{Kuzmin_MCL1}
D.~Kuzmin.
\newblock Monolithic convex limiting for continuous finite element
  discretizations of hyperbolic conservation laws.
\newblock {\em Comput. Method. Appl. M.}, 361:112804, 2020.

\bibitem{Kuzmin_BP}
D.~Kuzmin, M.~Quezada~de Luna, D.~Ketcheson, and J.~Gr{\"u}ll.
\newblock Bound-preserving flux limiting for high-order explicit {R}unge
  {K}utta time discretizations of hyperbolic conservation laws.
\newblock {\em J. Sci. Comput.}, 91:21, 2022.

\bibitem{wu2023geometric}
Kailiang Wu and Chi-Wang Shu.
\newblock Geometric quasilinearization framework for analysis and design of
  bound-preserving schemes.
\newblock {\em SIAM Review}, 65(4):1031--1073, 2023.

\bibitem{gmsh}
C.~Geuzaine and J.-F. Remacle.
\newblock Gmsh: a three-dimensional finite element mesh generator with built-in
  pre- and post-processing facilities.
\newblock {\em International Journal for Numerical Methods in Engineering},
  79(11):1309--1331, 209.

\bibitem{maire2011}
Pierre-Henri Maire.
\newblock A unified sub-cell force-based discretization for cell-centered
  lagrangian hydrodynamics on polygonal grids.
\newblock {\em International Journal for Numerical Methods in Fluids},
  65(11-12):1281--1294, 2011.

\bibitem{HTCLag_boscheri2024}
Walter Boscheri, Michael Dumbser, and Pierre-Henri Maire.
\newblock A new thermodynamically compatible finite volume scheme for
  lagrangian gas dynamics.
\newblock {\em SIAM Journal on Scientific Computing}, 46(4):A2224--A2247, 2024.

\bibitem{CWENO_BGK_boscheri2020}
Walter Boscheri and Giacomo Dimarco.
\newblock High order central weno-implicit-explicit runge kutta schemes for the
  bgk model on general polygonal meshes.
\newblock {\em Journal of Computational Physics}, 422:109766, 2020.

\bibitem{ader_fse_boscheri2020}
Walter Boscheri.
\newblock A space-time semi-lagrangian advection scheme on staggered voronoi
  meshes applied to free surface flows.
\newblock {\em Computers and Fluids}, 202:104503, 2020.

\bibitem{hitch}
L.~Beir{\~a}o~da Veiga, Franco Brezzi, L.~D. Marini, and A.~Russo.
\newblock The {Hitchhiker}'s guide to the virtual element method.
\newblock {\em Math. Models Methods Appl. Sci.}, 24(8):1541--1573, 2014.

\bibitem{orthVEM2017}
Stefano Berrone and Andrea Borio.
\newblock Orthogonal polynomials in badly shaped polygonal elements for the
  virtual element method.
\newblock {\em Finite Elements in Analysis and Design}, 129:14--31, 2017.

\bibitem{DeconinckMario}
H.~Deconinck and M.~Ricchiuto.
\newblock {\em Encyclopedia of Computational Mechanics}, chapter Residual
  distribution schemes: foundation and analysis.
\newblock John Wiley \& sons, 2007.
\newblock DOI: 10.1002/0470091355.ecm054.

\bibitem{Abgrall99}
R.~Abgrall.
\newblock {Toward the ultimate conservative scheme: Following the quest.}
\newblock {\em J. Comput. Phys.}, 167(2):277--315, 2001.

\bibitem{CiarletRaviart}
P.G. Ciarlet and P.A. Raviart.
\newblock General {Lagrange and Hermite Interpolation with Applications to
  Finite Element Methods}.
\newblock {\em Archive For Rational Mechanics and Application}, 46(3):177--199,
  1972.

\bibitem{chinois}
Long Chen and Jianguo Huang.
\newblock Some error analysis on virtual element methods.
\newblock {\em Calcolo}, 55(1):23, 2018.
\newblock Id/No 5.

\bibitem{AbgrallOeffnerLiuDG}
R{\'e}mi Abgrall, Philipp {\"O}ffner, and Yongle Liu.
\newblock Some new properties of the {PamPa} scheme.
\newblock {\em submitted}, 2025.
\newblock Preprint, {arXiv}:2508.17147 [math.{NA}] (2025).

\bibitem{WassilijPetrov}
Wasilij Barsukow.
\newblock Semi-discrete {Active} {Flux} as a {Petrov}-{Galerkin} method.
\newblock Preprint, {arXiv}:2508.15017 [math.{NA}] (2025), 2025.

\bibitem{GuermondPopovFast}
Jean-Luc Guermond and Bojan Popov.
\newblock Fast estimation from above of the maximum wave speed in the {Riemann}
  problem for the {Euler} equations.
\newblock {\em J. Comput. Phys.}, 321:908--926, 2016.

\bibitem{Vilar_DGFV}
F.~Vilar.
\newblock Local subcell monolithic {DG}/{FV} convex property preserving scheme
  on unstructured grids and entropy consideration.
\newblock {\em J. Comput. Phys.}, 521:113535, 2025.

\bibitem{GauthierLBM}
G.~Wissocq, Y.~Liu, and R.~Abgrall.
\newblock A positive- and bound-preserving vectorial lattice {B}oltzmann method
  in two dimensions.
\newblock {\em in preparation}, 2024.

\bibitem{zbMATH05380178}
Alexander Kurganov, Guergana Petrova, and Bojan Popov.
\newblock Adaptive semidiscrete central-upwind schemes for nonconvex hyperbolic
  conservation laws.
\newblock {\em SIAM J. Sci. Comput.}, 29(6):2381--2401, 2007.

\bibitem{lax1998}
Peter~D Lax and Xu-Dong Liu.
\newblock Solution of two-dimensional riemann problems of gas dynamics by
  positive schemes.
\newblock {\em SIAM Journal on Scientific Computing}, 19(2):319--340, 1998.

\bibitem{KLZ}
A.~Kurganov, Y.~Liu, and V.~Zeitlin.
\newblock Numerical dissipation switch for two-dimensional central-upwind
  schemes.
\newblock {\em ESAIM Mathematical Modelling and Numerical Analysis},
  55:713--734, 2021.

\bibitem{GKL}
N.~K. Grag, A.~Kurganov, and Y.~Liu.
\newblock Semi-discrete central-upwind {R}ankine-{H}ugoniot schemes for
  hyperbolic systems of conservation laws.
\newblock {\em Journal of Computational Physics}, 428:110078, 2021.

\bibitem{WDKL}
B.-S. Wang, W.~Don, A.~Kurganov, and Y.~Liu.
\newblock Fifth-order {A-WENO} schemes based on the adaptive diffusion
  central-upwind {R}ankine-{H}ugoniot fluxes.
\newblock {\em Communications on Applied Mathematics and Computation},
  5:295--314, 2023.

\bibitem{vem}
B.~Ahmed, A.~Alsaedi, F.~Brezzi, L.D. Marini, and A.~Russo.
\newblock Projectors for {Virtual Element Methods}.
\newblock {\em Comput. Math. Appl.}, 66(3), 2013.

\end{thebibliography}
